\documentclass[12pt,centertags]{amsart}
\usepackage{amssymb}
\usepackage[all]{xy}
\usepackage{amsmath,amstext,amsthm,a4,amssymb,amscd}
\usepackage[mathscr]{eucal}
\usepackage{mathrsfs}
\usepackage{epsf}
\usepackage{tikz}
\usepackage{extarrows}
\usetikzlibrary{matrix,arrows,decorations.pathmorphing}
\usepackage[a4paper,top=3cm,bottom=3cm,
left=2.2cm,right=2.2cm]
{geometry}
\numberwithin{equation}{section}
\allowdisplaybreaks[3]

\usepackage{color}
\usepackage{soul}

\newcommand{\field}[1]{\mathbb{#1}}
\newcommand{\C}{\field{C}}
\newcommand{\N}{\field{N}}
\newcommand{\Q}{\field{Q}}
\newcommand{\R}{\field{R}}
\newcommand{\Z}{\field{Z}}

\def\cC{\mathscr{C}}

\def\cU{\mathscr{U}}
\def\cE{\mathscr{E}}

\def\mA{\mathcal{A}}

\def\mD{\mathcal{D}}
\def\mE{\mathcal{E}}
\def\mF{\mathcal{F}}
\def\mG{\mathcal{G}}

\def\mL{\mathcal{L}}
\def\mM{\mathcal{M}}
\def\mN{\mathcal{N}}

\def\mU{\mathcal{U}}
\def\mV{\mathcal{V}}

\newcommand\mS{\mathcal{S}}

\def\Im{{\rm Im}}

\def\ra{\rangle}

\DeclareMathOperator{\End}{End}
\DeclareMathOperator{\Hom}{Hom}

\DeclareMathOperator{\rank}{rk}
\DeclareMathOperator{\Id}{Id}

\DeclareMathOperator{\tr}{Tr}

\DeclareMathOperator{\ind}{Ind}

\DeclareMathOperator{\td}{Td}

\DeclareMathOperator{\ch}{ch}

\newtheorem{thm}{Theorem}[section]
\newtheorem{lemma}[thm]{Lemma}
\newtheorem{prop}[thm]{Proposition}
\newtheorem{cor}[thm]{Corollary}

\theoremstyle{definition}
\newtheorem{rem}[thm]{Remark}
\theoremstyle{definition}
\newtheorem{defn}[thm]{Definition}

\newcommand{\be}{\begin{eqnarray}}
	\newcommand{\ee}{\end{eqnarray}}

\newcommand{\wi}{\widetilde}
\newcommand{\var}{\varepsilon}
\numberwithin{equation}{section}
\numberwithin{thm}{section}

\newcommand{\comment}[1]{}

\begin{document}
	
	\title[Localization formula for Bismut-Cheeger $\eta$-forms]
	{Localization formula for Bismut-Cheeger $\eta$-forms}
	
	\author{Bo LIU}
	\address{School of Mathematical Sciences,  Key Laboratory of 
		MEA (Ministry of Education) \& Shanghai Key Laboratory 
		of PMMP,  East China Normal University, Shanghai 
		200241, China}
	\email{bliu@math.ecnu.edu.cn}
	
	\author{Xiaonan MA}
	
	\address{Chern Institute of Mathematics \& LPMC, Nankai 
		University, Tianjin 300071, P.R. China 
	}  
	\email{xiaonan.ma@nankai.edu.cn}	
	
	\date{\today}
	
	\begin{abstract}
		We establish a localization formula for equivariant Bismut-Cheeger
		$\eta$-forms under a torus action, extending our previous work 
		\cite{LM20} on the localization of equivariant 
		$\eta$-invariants for a circle action. 
		One of the main technical ideas behind this generalization
		is to introduce a $g$-equivariant 
		differential K-theory model in the sense of
		 Liu \cite{Liu21}, which modifies Bunke-Schick's
		model \cite{Bunke2009,BunkeSchick13}. This construction
		enables us to derive a localization formula for pushforward maps 
		in differential K-theory. 
	\end{abstract}
	
	\maketitle
	
	
	\setcounter{section}{-1}
	

	\section{Introduction}\label{s00}  
	
	In our previous paper \cite{LM20}, we established a localization
	formula for equivariant $\eta$-invariants, which relates 
	an equivariant $\eta$-invariant on a closed $S^1$-spin$^c$ 
	manifold to the corresponding invariant on the fixed point set 
	of the $S^1$-action, using the differential K-theory which
	was introduced by Freed–Hopkins \cite{FH00}. 
	The present paper aims to generalize this 
	result to the family case with fiberwise torus action.	
	
	Let $G$ be a compact Lie group, and let $Z$ be a closed 
	odd-dimensional oriented $G$-manifold equipped with a 
	$G$-invariant Riemannian metric $g^{TZ}$. 
	An \textbf{equivariant geometric triple} $\underline{E} 
	= (E, h^E, \nabla^E)$ over $Z$ consists of a $\mathbb{Z}_2$-graded 
	complex vector bundle $E = E_+ \oplus E_-$ such that $E_{\pm}$ 
	are $G$-equivariant; a $G$-invariant Hermitian metric $h^E$ 
	on $E$ satisfying $E_+ \perp E_-$ and a $G$-invariant connection 
	$\nabla^E$ preserving both $h^E$ and the $\mathbb{Z}_2$-grading 
	of $E$. For an ungraded vector bundle $E$, we regard it as 
	$\mathbb{Z}_2$-graded with $E_- = 0$. The notation $E$ or 
	$\underline{E}$ is called \textit{ungraded} when we emphasize 
	the absence of $\mathbb{Z}_2$-grading.
	
	We assume $TZ$ admits a $G$-equivariant spin$^c$ structure, 
	meaning the $G$-action on $Z$ lifts naturally to its 
	associated spin$^c$ principal bundle. In particular, there 
	exists a $G$-equivariant complex line bundle $L_Z$ over $Z$ 
	satisfying the characteristic class relation
	$\omega_2(TZ) \equiv c_1(L_Z) \mod 2,	$
	where $\omega_2$ denotes the second Stiefel-Whitney class 
	and $c_1$ the first Chern class. 
	Given an equivariant geometric triple $\underline{L_Z}$
	associated with $L_Z$ over $Z$, we denote by $\bar{\eta}_g
	(\underline{TZ}, \underline{L_Z}, \underline{E})$ the 
	corresponding equivariant reduced $\eta$-invariant for 
	$g \in G$ (see \cite[Def. 1.2]{LM20} for details).
	
	We temporarily assume $G = S^1$. Let $Z^{S^1}$ be the fixed 
	point set of the circle action on $Z$. Each connected 
	component $Z_{\alpha}^{S^1}$ is then a closed $S^1$-equivariant 
	spin$^c$ manifold with an associated equivariant line bundle 
	$L_{\alpha}$ (see \cite[(1.33)]{LM20}). 
	Let $N_{\alpha}$ denote the normal bundle of $Z_{\alpha}^{S^1}$ 
	in $Z$, and let $\underline{L_{\alpha}}$ and $\underline{N_{\alpha}}$ 
	be the corresponding induced equivariant geometric triples over 
	$Z_{\alpha}^{S^1}$.	Set $A_{S^1}=\{g\in S^1: Z^g\neq Z^{S^1} \}$, 
	which is a finite subset of $S^1$ (see \cite[Prop. 1.1]{LM20}), 
	here $Z^g$ is the fixed point set of $g$-action.
	The following theorem is the main result in \cite{LM20}.
	
	\begin{thm}\cite[Th. 0.5]{LM20}\label{thm:0.01}
		There exist equivariant geometric triple 
		$\underline{\mE_{\alpha}}$ over $Z_{\alpha}^{S^1}$ 
		and $\chi_{\alpha}\in \Z[x]$, which depend on
		$\underline{N_{\alpha}}$, such that
		for any equivariant geometric triple $\underline{E}$ over $Z$,
		there exist $f, h\in \Z[x]$, such that for any $g\in
		S^1\backslash A_{S^1}$,
		\begin{align}\label{eq:0.03}
			\bar{\eta}_g\left(\underline{TZ}, \underline{L_Z}, 
			\underline{E}\right)-\sum_{\alpha}
			\bar{\eta}_g\left(\underline{TZ_{\alpha}^{S^1}}, 
			\underline{L_{\alpha}}, \underline{\mE_{\alpha}}\otimes 
			\underline{E}|_{Z_{\alpha}^{S^1}}\right)/\chi_{\alpha}(g)
			=\frac{f(g)}{h(g)}.
		\end{align}
		Note that $\underline{\mE_{\alpha}}/\chi_{\alpha}$ is 
		written as $\underline{\lambda_{-1}(N_{\alpha}^*)^{-1}_{\mN}}$
		in \cite{LM20} and constructed explicitly in 
		\cite[(2.100)]{LM20}.		
	\end{thm}
	
	The main purpose of this paper is to generalize Theorem
	\ref{thm:0.01} to the family case.
	
	Consider a $G$-equivariant submersion $\pi_Z: W \to B$ of closed 
	$G$-manifolds with closed fiber $Z$. The relative tangent bundle 
	$TZ = \ker(d\pi_Z: TW \to TB)$ 	forms an equivariant subbundle 
	of $TW$. We assume that $TZ$ is oriented and fix an orientation
	which corresponds to $o_Z\in H^0(W,\Z_2)$ by 
	\cite[Th. II.1.2]{LM89}.
	Let $g^{TZ}$ be a $G$-invariant Euclidean metric on $TZ$.
	We assume that the oriented vector bundle $TZ$ has a $G$-equivariant 
	spin$^c$ structure with associated equivariant complex line bundle 
	$L_Z$ over $W$.
	Let $\underline{L_Z}$ be an equivariant geometric triple
	associated with $L_Z$.
	Let $T_{\pi_Z}^HW$ be an equivariant horizontal subbundle of 
	$TW$ such that $TW=T_{\pi_Z}^HW\oplus TZ.$
	We set
	\begin{align}\label{eq:0.02}
		\underline{\pi_Z}:=(W, o_Z, T_{\pi_Z}^HW, g^{TZ}, \underline{L_Z}).
	\end{align}
	Let $\underline{E}$ be an equivariant geometric triple over $W$.
	We denote by
	\begin{align}\label{eq:1.02}
		\mF=(\underline{\pi_Z}, \underline{E})
	\end{align}
	an \textbf{equivariant geometric family} over $B$.	
	For an equivariant geometric family $\mF$, we can uniquely define
	the fiberwise Dirac operator  $D(\mF)$ with family index 
	$\ind(D(\mF))\in K_G^{\bullet}(B)$, the $G$-equivariant 
	$K^{\bullet}$-group of $B$.
	If the dimension of all fibers are even/odd, 
	$\ind(D(\mF))\in K_G^{0/1}(B)$.
	
	Let $\Omega^{\bullet}(B)$ be the space of smooth differential
	forms on $B$, and $\mathrm{Im}\,d=d\Omega^{\bullet}(B)$
	the space of exact forms on $B$.
	Given an equivariant geometric family $\mathcal{F}$ with trivial 
	$G$-action on $B$, the equivariant Bismut-Cheeger $\eta$-form, 
	the family extension of the equivariant $\eta$-invariant, can be 
	defined under two distinct conditions:
	\begin{enumerate}
		\item \textbf{(Bundle Condition):} When $\ker D(\mathcal{F})$ 
		forms a vector bundle over $B$, for each $g\in G$, the 
		equivariant Bismut-Cheeger $\eta$-form $\widetilde{\eta
		}_g(\mathcal{F})\in \Omega^{\bullet}(B)$ \cite{BC89,Liu17} is 
		well-defined (Definition \ref{defn:1.01}). In the special 
		case where $B$ is a point and $\dim Z$ is odd, $\widetilde{\eta
		}_g(\mathcal{F})$ reduces to half of the Atiyah-Patodi-Singer 
		$\eta$-invariant.
		
		\item \textbf{(Perturbation Condition):} When there exists 
		a perturbation operator $\mathcal{A}$ (Definition 
		\ref{defn:1.03}) making $D(\mathcal{F})+\mathcal{A}$ 
		fiberwise invertible (it is equivalent to $\ind(D(
		\mathcal{F}))=0\in K_G^{\bullet}(B)$ when the fiber has 
		positive dimension by Proposition \ref{prop:1.4}), 
		the equivariant Bismut-Cheeger $\eta$-form 
		$\widetilde{\eta}_g(\mathcal{F},\mathcal{A})\in \Omega^{\bullet}(B)/
		\mathrm{Im}\,d$ \cite{Liu21,MP97a,MP97b}  with perturbation 
		$\mA$ is well-defined modulo exact forms (Definition 
		\ref{defn:1.02}). If $B$ is a point and $\dim Z$ is odd,  
		$\mathcal{A}$ can be chosen so that $\widetilde{\eta}_g(
		\mathcal{F},\mathcal{A}) = \bar{\eta}_g(\underline{TZ},
		\underline{L_Z},\underline{E})$
		(cf. \cite[Remark 3.20]{Liu21}).
	\end{enumerate}
	
	For $G=\mathbb{T}^m$, $m$-dimensional torus, let $W^{\mathbb{T}}$ 
	be the fixed point set of the $\mathbb{T}^m$-action on $W$, 
	which is a totally geodesic submanifold of $W$.
	With the exception of Section \ref{s0304}, 
	we assume throughout that $W^{\mathbb{T}}$ is non-empty.
	The restriction $\pi_{Z^{\mathbb{T}}}:=\pi_Z|_{W^{\mathbb{T}}}:
	W^{\mathbb{T}}\to B$ forms a fiber bundle with closed fiber 
	$Z^{\mathbb{T}}$. 
	Assume that $B$ is connected and oriented.
	For each connected component $W^{\mathbb{T}}_{\alpha}$ of the 
	fixed point set, the relative tangent bundle $TZ_{\alpha}^{
		\mathbb{T}}$ corresponding to $\pi_{Z^{\mathbb{T}}_{\alpha}}: 
	W^{\mathbb{T}}_{\alpha} \to B$ inherits a canonical 
	$\mathbb{T}^m$-equivariant $\text{spin}^c$ structure, 
	with associated equivariant complex line bundle $L_{\alpha}$ 
	as specified in (\ref{eq:1.32a}). The geometric data 
	$\underline{\pi_Z}$ naturally induces corresponding data 
	$\underline{\pi_{Z^{\mathbb{T}}_{\alpha}}}$ on each 
	$W_{\alpha}^{\mathbb{T}}$ (see (\ref{eq:1.38b})). The embedding 
	$W_{\alpha}^{\mathbb{T}} \hookrightarrow W$ determines a normal 
	bundle $N_{\alpha}$ with equivariant complex structure, 
	yielding well-defined equivariant geometric triples 
	$\underline{L_{\alpha}}$ and $\underline{N_{\alpha}}$ over 
	$W_{\alpha}^{\mathbb{T}}$.
	
	Set $A_{\mathbb{T}}=\{g\in {\mathbb{T}}^m: W^g\neq 
	W^{\mathbb{T}} \}$. Then $A_{\mathbb{T}}$ is contained by 
	the union of a finite number of closed 
	subgroups of ${\mathbb{T}}^m$, each of positive codimension 
	(see Proposition \ref{prop:1.05}). 
	
	We expect a formula that replaces the $\eta$-invariants in the 
	left hand side of (\ref{eq:0.03}) by the corresponding 
	$\eta$-forms. Since the bundle condition is not easy to 
	verify, we seek an equivariant geometric family with
	vanishing family index to define the $\eta$-form. 
	
From the spin$^c$ version of the localization formula for the equivariant 
	family index in \cite[Lemma 2.7 and Prop. 2.10]{ASegal68}, 
	we have the following $\mathbb{T}^m$-equivariant 
	family version of \cite[Th. 0.1]{LM20}.
	
	\begin{thm}\label{thm:0.02}
	There exist equivariant geometric triple 
	$\underline{\mE_{\alpha}}$ over $W_{\alpha}^{\mathbb{T}}$ 
	and $\chi_{\alpha}\in \Z[x_1,\cdots,x_m]$, which depend on
	$\underline{N_{\alpha}}$, such that for any equivariant geometric 
	triple $\underline{E}$ over $W$, letting $\mathcal{F}_{\alpha}
	=\left(\underline{\pi_{Z^{\mathbb{T}}_{\alpha}}},  
	\underline{\mE_{\alpha}}\widehat{\otimes} \underline{E}|_{
		W_{\alpha}^{\mathbb{T}}}\right)$, 
		and for any $g\in	\mathbb{T}^m\backslash A_{\mathbb{T}}$, we have
	\begin{align}\label{eq:0.04}
	\ind (D(\mF))=\sum_{\alpha}\ind (D(\mF_{\alpha}))/\chi_{\alpha}^{}
	\in K_{\mathbb{T}^m}^{\bullet}(B)_{I(g)}.
	\end{align}	
	Here $K_{\mathbb{T}^m}^{\bullet}(B)_{I(g)}$ denotes the 
	localization of the $\mathbb{T}^m$-equivariant topological
	K-group at the prime ideal $I(g)$, which contains all
	characters of $\mathbb{T}^m$ vanishing at $g$ (see (\ref{eq:2.086})).
	\end{thm}
	
	For a proof of Theorem \ref{thm:0.02}, see 
	Corollary \ref{cor:3.5}.
	Note that there are many possible choices for $\mE_{\alpha}$ and 
$\chi_{\alpha}$ here, as long as $\mE_{\alpha}/\chi_{\alpha}$ is a geometric 
	realization of $\lambda_{-1}(N_{\alpha})^{-1}\in 
	K^0_{\mathbb{T}^m}(W^{\mathbb{T}})_{I(g)}$ (cf. \cite[Th. 1.1]{LMZ03} 
	for $m=1$). In Corollary \ref{cor:3.5}, $\underline{\mE_{\alpha}}$ 
	and $\chi_{\alpha}$ are constructed explicitly in
	(\ref{eq:2.092})-(\ref{eq:3.15}) and (\ref{eq:3.43}),
	which is the generalization of the choice in \cite{LM20}.
	If $m=1$, the constructions of $\underline{\mE_{\alpha}}$ 
	and $\chi_{\alpha}$ are the same as that in Theorem \ref{thm:0.01}.
	 In the reminder of the introduction,
	we always use this construction.

	For any finite dimensional $\mathbb{T}^m$-representation
	$M_{\chi}$ corresponding to the character $\chi$,
	we denote by $\underline{M_{\chi}}$ the equivariant 
	geometric family over $B$ given in (\ref{eq:2.43}),
	namely the trivial vector bundle with fiber $M_{\chi}$
	and with trivial metric and connection, and we write
	$\underline{M_{\chi}}\times_B\mF=(\underline{\pi_Z}, \pi_Z^*
	\underline{M_{\chi}}\widehat{\otimes} E)$. Then 
	we have
	\begin{align}\label{eq:0.05}
		\ind(D(\underline{M_{\chi}}\times_B\mF))
		=\chi\cdot \ind(D(\mF))\in K_{\mathbb{T}^m}^{\bullet}(B).
	\end{align}
	From (\ref{eq:0.04}) and (\ref{eq:0.05}),
	there exists character $\chi'$ such that
	for $\mF_M:=\underline{M_{\chi_0\chi'}}\times_B\mF$
	and $\mF_{M}^{\mathbb{T}}:=\sum_{\alpha} \underline{M_{\chi_0
			\chi'/\chi_{\alpha}}}
	\times_B \mF_{\alpha} $, we have
	\begin{align}\label{eq:0.06}
	\ind(D(\mF_M))=\ind(D(\mF_{M}^{\mathbb{T}}))\in K_{\mathbb{T}^m}^{\bullet}(B).
	\end{align}
	Here $\chi_0:=\prod_{\alpha}\chi_{\alpha}$
	and the sum of the equivariant geometric families is defined to
	 be their disjoint union.
	Let $\mF^{\mathrm{op}}$ denote the equivariant geometric 
	family which reverse the $\Z_2$-grading of $E$ in $\mF$.
	Then we have $\ind(D(\mF^{\mathrm{op}}))=-\ind(D(\mF))$.
	Thus from (\ref{eq:0.06}), we have
	\begin{align}\label{eq:0.07a}
	\ind(D(\mF_M+(\mF_{M}^{\mathbb{T}})^{\mathrm{op}}))=0
	\in K_{\mathbb{T}^m}^{\bullet}(B).
	\end{align}
	From Proposition \ref{prop:1.4}, Definition \ref{defn:1.02}
	and (\ref{eq:2.23}), there exists
	a perturbation $\mA$ such that 
	$\tilde{\eta}_g(\mF_M+(\mF_{M}^{\mathbb{T}})^{\mathrm{op}},
	\mA)\in \Omega^{\bullet}(B,\C)/\mathrm{Im}\,d$ is well defined.

	The main result of this paper is as follows.
	
	\begin{thm}\label{thm:0.03a}
	For any perturbation operator $\mA$ associated with
	$\mF_M+(\mF_{M}^{\mathbb{T}})^{\mathrm{op}}$,
	there exist $x\in K_{\mathbb{T}^m}^{\bullet}(B)$
	such that for any $g\in
	\mathbb{T}^m\backslash A_{\mathbb{T}}$, 
	modulo exact forms on $B$, we have
	\begin{align}
	\tilde{\eta}_g(\mF_M+(\mF_{M}^{\mathbb{T}})^{\mathrm{op}},
	\mA)=\ch_g(x).
	\end{align}
	\end{thm}
	
	In particular, we have the following corollary.
	
	\begin{cor}\label{thm:0.03}
	\begin{enumerate}
		\item (bundle condition) If $\ker D(\mathcal{F})$ and 
		$\ker D(\mathcal{F}_{\alpha})$ form vector bundles over $B$,
		there exist $f\in K^{\bullet}(B)[x_1,\cdots,x_m]$ and 
		$h\in \Z[x_1,\cdots,x_m]$, such that for any $g=(g_1,\cdots, 
		g_m)\in	\mathbb{T}^m\backslash A_{\mathbb{T}}$, modulo 
		exact forms on $B$,
		\begin{align}\label{eq:0.07}
			\widetilde{\eta}_g\left(\mathcal{F}\right)-\sum_{\alpha}
			\wi{\eta}_g\left(\mathcal{F}_{\alpha}\right)/\chi_{\alpha}(g)
			=\frac{\ch(f)(g_1,\cdots, g_m)}{h(g_1,\cdots, g_m)}.
		\end{align}
		\item (perturbation condition) If there exist perturbation 
		operators $\mA$ and $\mA_{\alpha}$ associated with $\mathcal{F}$
		and $\mathcal{F}_{\alpha}$ respectively, then there exist 
		$f\in K^{\bullet}(B)[x_1,\cdots,x_m]$ and $h\in \Z[x_1,\cdots,x_m]$, 
		such that for any $g=(g_1,\cdots, g_m)\in
		\mathbb{T}^m\backslash A_{\mathbb{T}}$, modulo exact forms on $B$,
		\begin{align}\label{eq:0.08}
			\widetilde{\eta}_g\left(\mathcal{F}, \mA\right)-\sum_{\alpha}
			\wi{\eta}_g\left(\mathcal{F}_{\alpha}, \mA_{\alpha}\right)
			/\chi_{\alpha}(g)
			=\frac{\ch(f)(g_1,\cdots, g_m)}{h(g_1,\cdots, g_m)}.
		\end{align}
		\end{enumerate}		
		Here $\ch: K^{0/1}(B)\to H^{\mathrm{even/odd}}(B,\R)$
		is the Chern character map.
	\end{cor}

If the dimension of all fibers of $\mathcal{F}$ is odd and 
$\ker D(\mF)$ forms a vector bundle over $B$,
the projection onto $\ker D(\mF)$ is a perturbation operator
associated with $\mF$ (cf. Definition \ref{defn:1.03}).
By Proposition \ref{prop:1.4}, we have 
$\ind(D(\mF))=0\in K_{\mathbb{T}^m}^1(B)$, which implies $\ind(D(\mF_M))$ 
vanishes. From (\ref{eq:0.06}), $\ind(D(\mF_M^{\mathbb{T}}))=0\in 
K_{\mathbb{T}^m}^1(B)$. Note that if the fixed point set
of $\mathbb{T}^m$-action on $W$ is non-empty, the dimension of any fiber 
of $\mathcal{F}_M^{\mathbb{T}}$ is non-zero.
So from Proposition \ref{prop:1.4} again, there exists 
a perturbation operator $\mA$ associated with 
$\mathcal{F}_M^{\mathbb{T}}$.

\begin{cor}\label{cor:0.5}
If the dimension of all fibers of $\mathcal{F}$ is odd and 
$\ker D(\mF)$ forms a vector bundle over $B$, for any 
perturbation operator $\mA$ associated with 
$\mathcal{F}_M^{\mathbb{T}}$,
there exist $f\in K^{\bullet}(B)[x_1,\cdots,x_m]$ and 
$h\in \Z[x_1,\cdots,x_m]$, such that for any $g=(g_1,\cdots, 
g_m)\in	\mathbb{T}^m\backslash A_{\mathbb{T}}$, modulo 
exact forms on $B$, we have
\begin{align}\label{eq:0.07b}
\widetilde{\eta}_g\left(\mathcal{F}\right)-
\wi{\eta}_g\left(\mathcal{F}_M^{\mathbb{T}},\mA \right)/(\chi_{0}\chi')(g)
=\frac{\ch(f)(g_1,\cdots, g_m)}{h(g_1,\cdots, g_m)}.
\end{align}
\end{cor}

\begin{rem}\label{rem:0.6}
	If $B$ is a point, then $K^0(B)=\Z$. The Chern character map
	on $K^0(B)$ is identical. In this case, we can choose $\mA$ in 
	Corollary \ref{cor:0.5} as the projections onto 
	$\ker D(\mF_M^{\mathbb{T}})$. By \cite[Remark 3.20]{Liu21},
	\begin{align}\label{eq:0.12a}
	\wi{\eta}_g\left(\mathcal{F}_M^{\mathbb{T}},\mA \right)
	=(\chi_0\chi')(g)\cdot\sum_{\alpha}\bar{\eta}_g
	(\underline{TZ_{\alpha}^{\mathbb{T}}},\underline{L_{\alpha}}, 
	\underline{\mE_{\alpha}}\widehat{\otimes} \underline{E}|_{
		Z_{\alpha}^{\mathbb{T}}})/\chi_{\alpha}(g).
	\end{align}
	Thus if  $m=1$, 
	Corollary \ref{cor:0.5} degenerates to Theorem \ref{thm:0.01}.
\end{rem}

We can also extend \cite[Th. 0.6]{LM20} to the family case.

\begin{thm}\label{thm:0.7}
If $W^{\mathbb{T}}=\emptyset$, $\ind(D(\mF))=0\in 
K_{\mathbb{T}^m}^{\bullet}(B)_{I(g)}$. 
Thus there exists character $\chi$ such that 
$\ind(D(\underline{M_{\chi}}\times_B\mF))=0
\in K_{\mathbb{T}^m}^{\bullet}(B)$. 
For a perturbation operator $\mA$ associated with 
$\underline{M_{\chi}}\times_B\mF$,
there exist 
$f\in K^{\bullet}(B)[x_1,\cdots,x_m]$ and $h\in \Z[x_1,\cdots,x_m]$, 
such that for any $g=(g_1,\cdots, g_m)\in
\mathbb{T}^m\backslash A_{\mathbb{T}}$, 
modulo exact forms on $B$,
\begin{align}\label{eq:0.08a}
	\widetilde{\eta}_g\left(\underline{M_{\chi}}\times_B\mF, \mA\right)
	=\frac{\ch(f)(g_1,\cdots, g_m)}{h(g_1,\cdots, g_m)}.
\end{align}
In particular, if $W^{\mathbb{T}}=\emptyset$ and $\ker D(\mF)$ forms 
a vector bundle over $B$, there exist 
$f\in K^{\bullet}(B)[x_1,\cdots,x_m]$ and $h\in \Z[x_1,\cdots,x_m]$, 
such that for any $g=(g_1,\cdots, g_m)\in
\mathbb{T}^m\backslash A_{\mathbb{T}}$, 
modulo exact forms on $B$,
\begin{align}\label{eq:0.08b}
	\widetilde{\eta}_g\left(\mathcal{F}\right)
	=\frac{\ch(f)(g_1,\cdots, g_m)}{h(g_1,\cdots, g_m)}.
\end{align}
\end{thm}

As in Remark \ref{rem:0.6}, if the dimension of all fibers
of $\mF$ is odd and $m=1$, (\ref{eq:0.08b}) degenerates to
\cite[Th. 0.6]{LM20}.

	In \cite{LM20}, Theorem \ref{thm:0.01} was established by 
	combining a commutative diagram for localized $g$-equivariant 
	differential K-theory for $g\in S^1\backslash A_{S^1}$,
	\begin{equation}\label{eq:0.09} 
		\begin{split}
			\xymatrix{
				\widehat{K}^0_g(Z^{S^1})_{I(g)} 
				\ar[rd]_{\widehat{f_{Y^{S^1}}}_!}
				&&\widehat{K}^0_g(Z)_{I(g)}
				\ar[ll]_{
					\left[\underline{\lambda_{-1}(N^*)},0\right]^{-1}
					\cup\,\hat{\iota}^*} 
				\ar[ld]^{\widehat{f_Y}_!}\\
				&\C/\Q_g,&
			}
		\end{split}
	\end{equation}
	with our extension of Goette's result \cite{Go00} on the comparison 
	between equivariant $\eta$-invariants and infinitesimal 
	equivariant $\eta$-invariants  \cite{LM22} here
	$\mathbb{Q}_g=\{P(g)/Q(g)\in \C : P, Q\in \mathbb{Z}[x],
	Q(g)\neq 0 \}$. For complete 
	details of this commutative diagram, see \cite{LM20}. 
	
	Since the comparison formula also works for equivariant 
	$\eta$-forms \cite{LM22}, the present work focuses on 
	generalizing the commutative diagram (\ref{eq:0.09}).
	In Theorem \ref{thm:1.08}, we extend (\ref{eq:0.09}) to 
	the following version for $g\in \mathbb{T}^m
	\backslash A_{\mathbb{T}}$:
	\begin{equation}\label{eq:0.10} 
		\begin{split}
			\xymatrix{
				\widehat{K}^0_g(W^{\mathbb{T}})_{I(g)} 
				\ar[rd]_{\widehat{\pi_{Z^{\mathbb{T}}}}_!}
				&&\widehat{K}^0_g(W)_{I(g)}
				\ar[ll]_{
					\left[\underline{\lambda_{-1}(N^*)},0\right]^{-1}
					\cup\,\hat{\iota}^*} 
				\ar[ld]^{\widehat{\pi_{Z}}_!}\\
				&\widehat{K}^{\bullet}_g(B)_{I(g)}.&
				}
		\end{split}
	\end{equation}
 If $B$ is a point, we have (see (\ref{eq:2.63}))
	\begin{align}\label{eq:0.11}
		\widehat{K}^1_g(B)_{I(g)}=\C/\Q_g.
	\end{align}
	
	The main difficulty in establishing (\ref{eq:0.10})
	is that the definition of the $g$-equivariant differential K-theory
	in \cite{LM20}, which orignates from Freed-Lott's model, 
	is hard to be used in this problem. It cannot manage the 
	$\widehat{K}^0$-group and $\widehat{K}^1$-group 
	in the same framework and does not have
	enough space to do analysis. 
	In \cite{Liu21}, we introduce another model for 
	equivariant differential K-theory, which modifies the model
	developped by Bunke-Schick \cite{Bunke2009} originally. 
	The differential K-theory defined by this model is equivalent to 
	the previous one, although the definition is more complicated. However,
	a $g$-equivariant modification of this model is just right for 
	solving this problem.
	Almost all analytic technical details, which are used in this paper,
	was discussed in \cite{Liu19,Liu21,LM22}. For closely related works, 
	cf.  \cite{GS21,HM20,Tang2021}.
	
	This paper is organized as follows. 
	In Section \ref{s01}, we introduce two types of equivariant 
	$\eta$-forms and review some of their key analytic properties. 
	In Section \ref{s02}, 
	we develop a $g$-equivariant version of the equivariant
	differential K-theory model in \cite{Liu21},
	which modifies Bunke-Schick's model in \cite{Bunke2009,BunkeSchick13}.
	In this framework, we construct pushforward maps along 
	submersions and embeddings and establish 
	the isomorphism between our model and the $g$-equivariant 
	version of Freed-Lott's model \cite{FreedLott10,LM20} for 
	differential K$^0$-ring. 
	In Section 3, utilizing the new model for 
	$g$-equivariant differential K-theory, we derive  the commutative 
	diagram (\ref{eq:0.10}) and prove Theorem \ref{thm:0.03a}
	by combining our result \cite{LM22}.
	Finally we deduce Corollaries \ref{thm:0.03}, \ref{cor:0.5}
	from Theorem \ref{thm:0.03a} in Section \ref{s0303} 
	and prove Theorem \ref{thm:0.7} in Section \ref{s0304}.
	
	\
	
	\noindent{\bf Notations}.
	We denote by $\lfloor x\rfloor$ the maximal integer not larger 
	than $x\in \R$.
	
	In the whole paper, if $\mA$, $\mA'$ are $\Z_2$-graded algebras
	we will note $\mA\widehat{\otimes}\mA'$ as the $\Z_2$-graded 
	tensor product as in \cite[\S 1.3]{BGV04}.
	If one of $\mA$ and $\mA'$ is ungraded, we understand 
	it as $\Z_2$-graded by taking its odd part as zero.
	
	We denote by $R(G)$ the representation ring of a compact Lie
	group $G$. As a finite dimensional virtual 
	representation of $G$ is characterized by its character,
	we will not distinguish an element in $R(G)$ and its character.

	For the fiber bundle $\pi_Z: W\rightarrow B$, we will often use 
	the integration of the differential forms along the oriented 
	fiber $Z$ in this paper. Since the fiber may be odd 
	dimensional, we must make precisely our sign conventions: 
	for any $\alpha\in \Omega^{\bullet}(B)$ and 
	$\beta\in \Omega^{\bullet}(W)$, 
	\begin{align}\label{eq:0.11a}
		\int_Z(\pi^*\alpha)\wedge\beta=\alpha\wedge \int_Z\beta.
	\end{align}		
	\
	
\noindent{\bf Acknowledgments}. 
We would like to thank Prof. Weiping Zhang for helpful discussions. 
BL is partially supported by NSFC No.12225105, National Key R\&D Program 
of China No.2024YFA1013200 and Science and Technology 
Commission of Shanghai Municipality No.22DZ2229014.
XM is partially supported by Nankai Zhide Foundation, NSFC No. 12621005,
Central Government-Guided Fund Project for Local Science 
and Technology Development in Tianjin (No.\ 25ZYZYJC00010).
	
	
	\section{Equivariant Bismut-Cheeger $\eta$-forms}\label{s01}
	
	In this section, we present some properties of equivariant 
	Bismut-Cheeger $\eta$-forms. In Section \ref{s0101}, we 
	describe the torus action on a closed Riemannian manifold.
	In Section \ref{s0102}, we introduce the equivariant 
	geometric family. In Section \ref{s0103a}, we explain 
	the orientation reversal trick for odd geometric family.
	In Section \ref{s0104a},
	we recall the definitions of 
	two types of equivariant Bismut-Cheeger $\eta$-forms under 
	different conditions and explain the relations between them. 
	In Section \ref{s0103},  we define the equivariant infinitesimal
	$\eta$-form with perturbation and compare it with the
	equivariant $\eta$-forms with perturbation, which is 
	the parallel version of the comparison formula in 
	\cite{LM22}. In Section \ref{s0104}, we revisit the embedding 
	formula for equivariant $\eta$-forms in \cite{Liu19}. For 
	practical application, we focus specifically on the case where 
	the small fiber bundle corresponds to the fixed point set of the 
	torus action on the large fiber bundle.
	
	
	\subsection{Torus action}\label{s0101}
	
	Let $W$ be a closed Riemannian manifold with a smooth
	$m$-dimensional torus
	$\mathbb{T}^m$-action preserving the Riemannian metric.
	Let $W^{\mathbb{T}}$ be the fixed point set of the torus action, 
	which is a totally geodesic submanifold of $W$.
	Let $N$ be the normal bundle of $W^{\mathbb{T}}$ in $W$.
	Then for $y\in W^{\mathbb{T}}$, $N_y$ is a 
	$\mathbb{T}^m$-representation which decomposes 
	into a direct sum of non-trivial representations
	of $\mathbb{T}^m$
	\begin{align}\label{eq:1.004}
		N_y=N(\beta_1)_{y}\oplus \cdots\oplus N(\beta_k)_y.
	\end{align}
	Here for $1\leq j\leq k$, the weight $\beta_j\in\Hom(
	\mathbb{T}^m, S^1)$ and $N(\beta_j)_y$ is an 
	even-dimensional real vector space which has a complex structure
	such that $g\in \mathbb{T}^m$ acts by multiplication by
	$\beta_j(g)\in S^1$.
	For $g=(g_1,\cdots, g_m)\in \mathbb{T}^m$, there exists 
	$(\ell_{j1},\cdots, \ell_{jm})\in \mathbb{Z}^m\backslash\{0\}$, such that 
	\begin{align}\label{eq:1.005}
		\beta_j(g)=g_1^{\ell_{j1}}\cdots g_m^{\ell_{jm}}.
	\end{align}
	For $g_s=e^{2\pi \mathrm{i}\theta_s}$, $1\leq s\leq m$, we have 
	\begin{align}\label{eq:1.003}
	\beta_j(e^{2\pi i\theta_1},\cdots,e^{2\pi \mathrm{i}\theta_m})
	=\exp\left(2\pi \mathrm{i}\sum_{s=1}^m\ell_{js}\theta_s\right).
	\end{align}
	
	Let $\{W_{\alpha}^{\mathbb{T}} \}_{\alpha\in \mathfrak{B}}$ 
	be the connected components of $W^{\mathbb{T}}$ 	
	and $N_{\alpha}$ be the normal bundle 	of $W^{\mathbb{T}}_{\alpha}$ 
	in $W$. 
	From (\ref{eq:1.004}) and the proof of \cite[Prop. 2.2]{
		Segal68}, $N_{\alpha}$ can be decomposed 
	into a direct sum of $\mathbb{T}^m$-equivariant vector bundles
	\begin{align}\label{eq:1.04a}
	N_{\alpha}=N_{\alpha}(\beta_1)\oplus \cdots\oplus N_{\alpha}(\beta_k).
	\end{align}	
	Here the weights $\beta_1,\cdots, \beta_k\in\Hom(
	\mathbb{T}^m, S^1)$ are distinct and $N_{\alpha}(
	\beta_j)$ admits a $\mathbb{T}^m$-equivariant complex structure
	such that $g\in \mathbb{T}^m$ acts by multiplication by 
	$\beta_j(g)\in S^1$. Thus $N_{\alpha}$ can be regarded as a 
	$\mathbb{T}^m$-equivariant complex vector bundle over 
	$W_{\alpha}^{\mathbb{T}}$.
	
	For $g\in \mathbb{T}^m$, we denote by $W^g$ the fixed point 
	set of $g$-action. 	
	Set 
	\begin{align}
		A_{\mathbb{T}}:=\{g\in {\mathbb{T}}^m: W^g\neq W^{\mathbb{T}}\}.
	\end{align}
	The following proposition is the generalization of 
	\cite[Prop. 1.1]{LM20}.
	
	\begin{prop}\label{prop:1.05}
		Let $W$ be a closed manifold with a smooth $\mathbb{T}^m$-action.
		Then there exists a finite number of
		closed subgroups $H_0,\cdots, H_t$ of $\mathbb{T}^m$,
		each of positive codimension, such that
		\begin{align}\label{eq:1.26} 
			A_{\mathbb{T}}\subset H_0\cup\cdots\cup H_t.
		\end{align}	
	\end{prop}
	\begin{proof}
		We may assume that the group action on $W$ is nontrivial.		
		Note that there exists $\var>0$ such that the exponential map of $W$
		\begin{align}\label{eq:1.5b}
			(y,Z)\in \mU_{\var}:=\{(y,Z)\in N_y: y\in W^{\mathbb{T}},  
			|Z|<\var \}\rightarrow \exp_y^W(Z)
		\end{align}
		is a diffeomorphism from $\mU_{\var}$ into a tubular 
		neighborhood $\mV_{\var}$ of $W^{\mathbb{T}}$ in $W$. 
		Thus for $(y,Z)\in \mU_{\var}$ 
		with $Z=(Z_1,\cdots, Z_k)\in N_{\alpha, y}$, where $Z_i\in 
		N(\beta_i)_y$, by (\ref{eq:1.003}) and (\ref{eq:1.04a}), we have
	\begin{align}\label{eq:1.9b} 
		g(y,Z)=\left(y, \left(e^{2\pi \mathrm{i} 
		\sum_{s=1}^m\ell_{1s}\theta_s}
		Z_1, \cdots, e^{2\pi \mathrm{i} \sum_{s=1}^m\ell_{ks}\theta_s}Z_k
		\right)\right).
		\end{align}
		Note that $\ell_{js}$ is locally constant on $W^{\mathbb{T}}$. 
		For $y\in W_{\alpha}^{\mathbb{T}}$,
		let $A_{\alpha,j}$ be the subset of $[0,1]^m$ which consists of
		the solutions of 
		\begin{align}
			\ell_{j1}\theta_1+\cdots+\ell_{jm}\theta_m\in \Z 
		\end{align}
		for $(\theta_1,\cdots, \theta_m)\in [0,1]^m$.
		Easy to see that $A_{\alpha,j}$ is the union of a finite number
		of closed submanifolds of $[0,1]^m$ and 
		$e^{2\pi \mathrm{i} A_{\alpha,j}}$ is a codimension
		1 subgroup of $\mathbb{T}^m$. Thus 
		\begin{align}\label{eq:1.10b} 
		\mU_{\var}^g=W^{\mathbb{T}}
		\quad \text{if}\quad
			(\theta_1,\cdots, \theta_m)\notin 
			\bigcup_{\alpha}\bigcup_j  A_{\alpha,j}.
		\end{align}
		
		Let $W_1:=W\backslash\mU_{\var/2}$.
		Then for any $x\in W_1$, the stabilizer $\mathbb{T}_x$ of $x$ is 
		a proper closed subgroup of $\mathbb{T}^m$
		whose codimension is positive.
		Set $N_x=T_xW_1/T_x(\mathbb{T}^m\cdot x)$, which is a linear 
		representation of $\mathbb{T}_x$.
		By the slice theorem, there exists a $\mathbb{T}^m$-equivariant 
		diffeomorphism from an equivariant open neighborhood $U_x$
		of the zero section in $\mathbb{T}^m\times_{\mathbb{T}_x}N_x$ to 
		an open neighborhood of $\mathbb{T}^m\cdot x$ in $W_1$, which 
		sends the zero section $\mathbb{T}^m/\mathbb{T}_x$ onto the orbit 
		$\mathbb{T}^m\cdot x$ by the map
		$g\in \mathbb{T}^m\mapsto g\cdot x$. Now on this neighborhood $U_x$,
		if $g\notin \mathbb{T}_x$, $U_x^g=\emptyset$. 
		By using the compactness of $W_1$  
		we obtain (\ref{eq:1.26}).
		
		The proof of Proposition \ref{prop:1.05} is complete.
	\end{proof}

	
	\subsection{Equivariant geometric family}\label{s0102}
	
	In the remainder of this section, unless otherwise stated, we assume that 
	$G$ is a compact Lie group.
	
	Let $\pi_Z:W\rightarrow B$ be a $G$-equivariant submersion of 
	closed $G$-equivariant Riemannian manifolds with closed fiber 
	$Z$ such that the relative tangent bundle $TZ$ over $W$ is 
	oriented and has a $G$-equivariant spin$^c$ structure.
	Let $\underline{E}:=(E, h^E, \nabla^E)$ be 
	an equivariant geometric triple over $W$ and
	\begin{align}\label{eq:1.04}
		\mF:=(\underline{\pi_Z}, \underline{E})
	\end{align}
	be an equivariant geometric family over $B$ defined in 
	(\ref{eq:1.02}). The opposite 
	family $\mF^{\mathrm{op}}$ of $\mF$ is defined by reversing 
	the $\Z_2$-grading of $E$. 
	For two equivariant geometric families $\mF$, $\mF'$ 
	over $B$, we can form their sum $\mF+\mF'$ over $B$
	as a new equivariant geometric family: the underlying fibration of 
	the sum is $\pi\sqcup \pi':W\sqcup W'\to B$,
	where $\sqcup$ is the disjoint union and the remaining 
	structures of $\mF+\mF'$ are induced in the obvious way.
	For two equivariant geometric triples $\underline{E}$ and
	$\underline{F}$, we can identify
	$(\underline{\pi_Z}, \underline{E\oplus F})$  with 
	$(\underline{\pi_Z}, \underline{E})+(\underline{\pi_Z}, \underline{F})$ 
	through the natural bijection. 
	 We say an equivariant geometric family 
	$\mF$ is even (resp. odd) if for any connected component of 
	the fiber, the dimension is even (resp. odd).
	For any equivariant geometric family $\mF$, it can be
	decomposed into $\mF=\mF^{\mathrm{even}}+\mF^{\mathrm{odd}}$
	uniquely, where $\mF^{\mathrm{even}}$ 
	(resp. $\mF^{\mathrm{odd}}$) is even (resp. odd).
	
	Let $T_{\pi_Z}^HW $ be a $G$-invariant horizontal 
	subbundle of $TW$.
	Let $g^{TB}$ be a $G$-invariant metric on $TB$.
	We equip 
	\begin{align}\label{eq:1.05}
		TW=T_{\pi_Z}^HW\oplus TZ
	\end{align} 
	with the $G$-invariant metric
	$g^{TW}=\pi_Z^*g^{TB}\oplus g^{TZ}$. Let $\nabla^{TW}$ be the 
	Levi-Civita connection on $(TW, g^{TW})$. Let 
	$P^{TZ}:TW\to TZ$ be the orthogonal projection associated with 
	(\ref{eq:1.05}).
	The connection $\nabla^{TZ}$ on $TZ$ is defined by (see e.g., 
	\cite[Def. 1.6]{Bi86}, \cite[(1.2)]{Bismut1997})
	\begin{align}\label{eq:1.15}
		\nabla^{TZ}:=P^{TZ}\nabla^{TW}P^{TZ},
	\end{align} 
	which restricts to the Levi-Civita connection of 
	$(TZ, g^{TZ})$ 
	on each fiber and only depends on ($T_{\pi_Z}^HW, g^{TZ}$)
	(see \cite[Th. 1.9]{Bi86}). Let $C(TZ)$ be the Clifford 
	algebra bundle of ($TZ, g^{TZ}$), whose fiber at $x\in W$
	is the Clifford algebra $C(T_xZ)$ of the Euclidean space
	$(T_xZ, g^{T_xZ})$. Let $\mathcal{S}(TZ, L_Z)$ 
	be the equivariant spinor bundle associated with the 
	equivariant spin$^c$ structure on $TZ$.
	This bundle carries a smooth action of $C(TZ)$, called the Clifford action.
	Let $\nabla^{\mathcal{S}_Z}$ be the connection on  
	$\mathcal{S}(TZ, L_Z)$ induced from $\nabla^{TZ}$ and 
	$\nabla^{L_Z}$. On $\mathcal{S}(TZ, L_Z)\widehat{\otimes}E$, 
	we set $\nabla^{\mathcal{S}_Z\widehat{\otimes} E}:=
	\nabla^{\mathcal{S}_Z}\widehat{\otimes} 1+1 \widehat{\otimes}
	\nabla^E$. For an equivariant geometric family $\mathcal{F}$,
	the fiberwise Dirac operator $D(\mathcal{F})$ associated with 
	$\mathcal{F}$ is defined by 
	\begin{align}\label{e01029}
		D(\mF):=\sum_ic(e_i)\nabla_{e_i}^{\mS_Z\widehat{\otimes} E},
	\end{align}
	where $\{e_i\}$ is a local orthonormal frame of $TZ$
	and $c$ is the Clifford action.
	Note that the definition of the fiberwise Dirac operator
	is independent of the choice of the local orthonormal frame.
	From (\ref{e01029}), the $G$-action commutes with $D(\mF)$.
	
	For an even (resp. odd) equivariant geometric family $\mF$,
	the classical construction of Atiyah-Singer
	assigns to this family its equivariant (analytic)
	index $\ind(D(\mF))\in K_G^0(B)$ (resp. $K_G^1(B)$)  
	(cf. \cite{AS69,AS71} or \cite[Appendix B]{Liu21}). 
	Here $K_G^{0/1}(B)$ denotes the $G$-equivariant 
	topological K$^{0/1}$-group of $B$. 
	Let $K_G^{\bullet}(B):=K_G^0(B)\oplus K_G^1(B)$.
	For any equivariant geometric family 
	$\mF=\mF^{\mathrm{even}}+\mF^{\mathrm{odd}}$, 
	we often write $\mathrm{Ind}(D(\mF))\in K_G^{\bullet}(B)$,
	meaning that the equivariant family index of $\mF^{\mathrm{even}}$ 
	(resp. $\mF^{\mathrm{odd}}$)
	lies in $K_G^0(B)$ (resp. $K_G^1(B)$).
	For the odd case, since the twisted bundle $E=E_+\oplus E_-$
	is $\Z_2$-graded here, we need to make precisely our notations.
	For $\mF_+=(\underline{\pi_Z}, \underline{E_+}\oplus 0)$ and
	$\mF_-=(\underline{\pi_Z}, 0\oplus\underline{E_-})$, we have
	$\mF=\mF_++\mF_-$ and define
	\begin{align}\label{eq:1.14}
	\ind(D(\mF)):=\ind(D(\mF_+))-\ind(D(\mF_-^{\mathrm{op}}))\in K_G^1(B).
	\end{align}
	Note that the definition of the usual Dirac operator is independent
	of the $\Z_2$-grading.
	Since the only difference between $\mF_-$ and $\mF_-^{\mathrm{op}}$
	is the $\Z_2$-grading, we may regard 
	\begin{align}\label{eq:1.15a}
	D(\mF)=D(\mF_+)\oplus (-D(\mF_-))
	\end{align}
	in this case.
However, this point of view is not entirely clear as a definition 
and may be confusing. We will realize this expression rigorously
as a ungraded version by reversing the orientation of the fiber of the 
$\mF_-$ part in Section \ref{s0103a}.

	By \cite[Def. 2.7, 2.8]{
		Segal68} (see also \cite[(3.1)]{LM20}), $K_G^1(B)$ is 
	characterized by the exact sequence 
	\begin{align}\label{eq:3.032} 
		0\rightarrow K_{G}^{1}(B)\rightarrow
		K_{G}^0(B\times S^1)\overset{i^*}  
		{\rightarrow}K_{G}^0(B)\rightarrow 0.
	\end{align}
	Here the $G$-action on $S^1$ is trivial and $i:B\to B\times S^1$, 
	$x\mapsto (x,*)$, is an embedding for $*\in S^1$ fixed. 
	For $g\in G$, we have the $\mathbb{Z}_2$-graded ring
	homomorphism
	\begin{align}\label{eq:1.28}
		\ch_g:K^{0/1}_G(B)\to H^{\mathrm{even/odd}}(B^g, \mathbb{C}).
	\end{align}
	For $y\in K^1_G(B)$, we see the definition of $\ch_g(y)$ in 
	\cite[(3.7)]{LM20}. 
	
	\subsection{Orientation reversal}\label{s0103a}
	
	In this subsection, we explain the orientation reversal 
	explanation for odd equivariant geometric family which is 
	originally used in \cite{Liu21}.
	
	Let $V$ be an Euclidean space of dimension $n$ with orientation
	$o$. 
	Then the complex Clifford
	algebra $C(V)$ is naturally $\Z_2$-graded from the
	$\N$-grading of the tensor algebra after reduction mod $2$.
	We write $C(V)=C_0(V)\oplus C_1(V)$.
	
	If $n=2m$ is even, the spinor $\mS(V)$ is the unique
	irreducible representation of $C(V)$. 
	In this subsection, we also write the spinor as $\mS_o(V)$
	and the Clifford action for $v\in C(V)$ as $c_o(v)$
	to specify the orientation $o$.
	Let $\tau_o$ denote the canonical $\Z_2$-grading of $\mS_o(V)$.
	Then $\mS_o(V)=\mS_{o,+}(V)\oplus \mS_{o,-}(V)$ with 
	respect to $\tau_o$. 
	
	If $n=2m-1$ is odd, 
	\begin{align}\label{eq:a.01}
		j: C(V)\to C_0(V\oplus \R), \quad v\mapsto v\cdot e
	\end{align}
	defines an isomorphism
	of algebras, where
	$e$ is the unit vector in $\R$. Let $\tilde{o}$ denote the 
	orientation of $V\oplus \R$ induced by the orientation $o$
	of $V$. Then
	$\mS_{\tilde{o},\pm}(V\oplus \R)$ are also the two inequivalent
	irreducible representations of $C(V)$ through the isomorphism $j$. 
	In this paper, we always take $\mS_o(V):=\mS_{\tilde{o},+}(V\oplus \R)$
	as the spinor for $C(V)$ as in \cite[\S 2.1]{Liu21}.
	If we reverse the orientation of $V$, 
	the orientation of $V\oplus \R$ is reversed and
	$\tau_{-\tilde{o}}=-\tau_{\tilde{o}}$. Thus we have
	$\mS_{-\tilde{o},+}(V\oplus \R)=\mS_{\tilde{o},-}(V\oplus \R)$.
	Set
	\begin{align}\label{eq:a.02}
		U:=c_{\tilde{o}}(e): 
		\mS_o(V)=\mS_{\tilde{o},+}(V\oplus \R)\longrightarrow 
		\mS_{-o}(V)=\mS_{\tilde{o},-}(V\oplus \R).
	\end{align}
	Then $U$ is an unitary isomorphism.
	By (\ref{eq:a.01}) and 
	(\ref{eq:a.02}), we have
	\begin{align}\label{eq:a.03}
		U c_o(v)=c_{\tilde{o}}(e)c_{\tilde{o}}(v)c_{\tilde{o}}(e)
		=-c_{\tilde{o}}(v)c_{\tilde{o}}(e)c_{\tilde{o}}(e)=-c_{-o}(v)U.
	\end{align} 
	
	From now on we consider an odd equivariant geometric family
	$\mF$ written as in (\ref{eq:1.02}). Then from (\ref{eq:a.02}), 
	we obtain an isomorphism
	\begin{align}\label{eq:a.04}
		U_{E}:\mS_{o_Z}(TZ, L_Z)\otimes E_-\to 
		\mS_{-o_Z}(TZ, L_Z)\otimes E_-
	\end{align}
	fiberwise by $U_{E}=U\otimes \mathrm{Id}_{L_Z^{1/2}\otimes E_-}$ locally.
	
	For $b\in B$, let $\mathscr{E}_b$ be the set of smooth sections 
	over $Z_b:=\pi^{-1}(b)$ of $\mathcal{S}(TZ, L_Z)\widehat{\otimes}
	E|_{Z_b}$. As in \cite{Bi86}, we regard $\mathscr{E}$ as an 
	infinite dimensional vector bundle over $B$ equipped
	with a fiberwise $L^2$-inner product.
	Let $\cE_{\pm}$ (resp. $\cE_{\pm}^{re}$) be the infinite dimensional
	bundle over $B$, the fiber of which is the smooth sections
	of $\mS_{o_Z}(TZ, L_Z)\otimes E_{\pm}$ (resp. $
	\mS_{-o_Z}(TZ, L_Z)\otimes E_{\pm}$) over $Z$.
	Then $U_E$ induces a fiberwise unitary map $\mU:\cE_-\to 
	\cE_-^{re}$.
	Let $\nabla^{\mS_Z\otimes E_-,o}$ denote the 
	natural connection on $\mS_{o_Z}(TZ, L_Z)\otimes E_-$.
	Then by (\ref{eq:a.04}), we have
	\begin{align}\label{eq:a.05}
		\mU \nabla^{\mS_Z\otimes E_-,o} =\nabla^{\mS_Z\otimes E_-,-o} \mU.
	\end{align}
	
	Set 
	\begin{align}
		\underline{-\pi_Z}:=(W, -o_Z, T_{\pi_Z}^HW, g^{TZ}, \underline{L_Z})
	\end{align}
	and 
	\begin{align}\label{eq:a.07}
		\mF_+:=(\underline{\pi_Z}, \underline{E_+}\oplus 0),\quad
		\mF_-:=(\underline{\pi_Z}, 0\oplus \underline{E_-}),\quad 
		\mF_-^{re}:=(\underline{-\pi_Z}, \underline{E_-}\oplus 0).
	\end{align}
	Then there is a natural bijective map $\mF_-\to \mF_-^{re}$
	by identification, which was used in the proof of
	\cite[Prop. 3.3]{Liu21}. Under this identification,
	from (\ref{eq:a.03}) and (\ref{eq:a.05}), we have
	\begin{align}\label{eq:a.08}
		D(\mF_-^{re})=-\mU D(\mF_-^{\mathrm{op}})\,\mU^{-1} :\cE_-^{re}
		\to \cE_-^{re}
	\end{align}
	and
	\begin{align}\label{eq:a.09}
		\ind (D(\mF_-^{re}))=-\ind (D(\mF_-^{\mathrm{op}}))
		\in K^1_G(B).
	\end{align}
	
	From (\ref{eq:a.08}), the Dirac operator $D(\mF)$
	in (\ref{eq:1.15a}) is explained as
	\begin{align}\label{eq:1.28a}
		(\mathrm{Id}\oplus \mU^{-1})D(\mF_++\mF_-^{re})
		(\mathrm{Id}\oplus \mU).
	\end{align}
	and
	$\ind(D(\mF))=
	\ind(
	D(\mF_++\mF_-^{re}))\in K^1_G(B)$ where $\mF_++\mF_-^{re}$
	is ungraded.

	\subsection{Equivariant $\eta$-form}\label{s0104a}
	
	In this subsection, we introduce the definitions of 
	two types of equivariant $\eta$-forms.
	
	For the
	infinite dimensional vector bundle $\mathscr{E}$ over $B$,
	whose fiber is the set of smooth sections 
	of $\mathcal{S}(TZ, L_Z)\widehat{\otimes}
	E$ over $Z$,
	let $\nabla^{\mathscr{E},u}$ be the connection on $\mathscr{E}$
defined by Bismut-Freed in \cite[Def. 1.3]{BF86}, 
which preserves the $L^2$-inner
	product on $\mathscr{E}$. If $U\in TB$, let $U^H\in T_{\pi_Z}^HW$ 
	be its horizontal lift in $T_{\pi_Z}^HW$ such that $\pi_{Z,*}(U^H)=U$.
	Let $\{f_p \}$ be a local frame of $TB$ and denote by $f^p$ 
	the dual of $f_p$.
	By \cite[(3.18)]{Bi86} (cf. also \cite[(1.23)]{LM22}), 
	the rescaled Bismut superconnection
	$\mathbb{B}_u:\cC^{\infty}(B, \Lambda(T^*B)\widehat{\otimes}
	\mathscr{E})\to \cC^{\infty}(B, \Lambda(T^*B)\widehat{\otimes}
	\mathscr{E})$, $u>0$, is defined by
	\begin{align}\label{eq:1.16a}
		\mathbb{B}_u:=\sqrt{u}D(\mathcal{F})+\nabla^{\mathscr{E},u}-
		\frac{1}{8\sqrt{u}}c\left(P^{TZ}[f_p^H, f_q^H] \right)f^p\wedge
		f^q\wedge.
	\end{align}
	Note that $\mathbb{B}_u^2$ is a second-order elliptic differential 
	operator along the fiber.
	
	We define the supertrace $\tr_s$ on a trace class operator
	in $\mathrm{End}(\mathscr{E})$ as in \cite[\S 1.3]{BGV04}.
	For $\sigma=\alpha\widehat{\otimes}A$ with $\alpha\in \Lambda
	(T^*B)$, $A\in \mathrm{End}(\mathscr{E})$, we define
	\begin{align}\label{e01056}
		\tr_s[\sigma]=\alpha\cdot\tr_s[A].
	\end{align}
	We denote by $\tr_s^{\mathrm{odd/even}}[\sigma]$ 
	the part of $\tr_s[\sigma]$
	which takes values in odd or even forms. 
	For $\alpha\in \Omega^j(B)$, set
	\begin{align}\label{1.04}
		\psi_B(\alpha)=\left\{
		\begin{array}{ll}
			\left(\frac{1}{2\pi \mathrm{i}}\right)^{\frac{j}{2}}
			\cdot \alpha, & \hbox{if $j$ is even;} \\
			\frac{1}{\sqrt{\pi}}\left(\frac{1}{2\pi \mathrm{i}}
			\right)^{\frac{j-1}{2}}\cdot \alpha, & \hbox{if $j$ 
				is odd.}
		\end{array}
		\right.
	\end{align}
	For $\beta\in \Omega^{\bullet}(B\times [0,1]_u)$, if
	$\beta=\beta_0+du\wedge \beta_1$, with $\beta_0, \beta_1\in 
	\Lambda (T^*B)$, set
	$\{\beta\}^{du}:=\beta_1.$
	
	\begin{defn}\label{defn:1.01} \cite[Def. 2.3]{Liu17}
		Assume that the $G$-action on $B$ is trivial and
		$\ker D(\mathcal{F})$ forms a vector bundle over $B$.
		For $g\in G$, the equivariant Bismut-Cheeger $\eta$-form 
		is defined by
		\begin{align}\label{eq:1.07}
			\tilde{\eta}_g(\mathcal{F}):=
			-\int_0^{\infty}\left\{\psi_{\R\times B}\left.\widetilde{\tr}
			\right.\left[g\exp\left(-\left(
			\mathbb{B}_{u}+du\wedge\frac{\partial}{\partial u}\right)^2
			\right)\right]\right\}^{du}du\in \Omega^{\bullet}(B,\C),
		\end{align}	
		where
		\begin{align}\label{i16}
			\widetilde{\tr}=
			\left\{
			\begin{array}{ll}
				\tr_s & \hbox{if $\mF$ is even;} \\
				\tr_s^{\mathrm{odd}} & \hbox{if $\mF$ is odd.}
			\end{array}
			\right.
		\end{align}
		When $\mF$ is even (resp. odd), $\tilde{\eta}_g(\mathcal{F})
		\in \Omega^{\mathrm{odd}}(B,\C)$ (resp. $\Omega^{
			\mathrm{even}}(B,\C)$).
	\end{defn}
	
	If $B$ is a point and $\dim Z$ is odd, for 
	$\bar{\eta}_g(\underline{TZ}, 
	\underline{L_Z},\underline{E})$ in \cite[(1.15)]{LM20},
	we have
	\begin{align}
		\tilde{\eta}_g(\mathcal{F})=\bar{\eta}_g(\underline{TZ}, 
		\underline{L_Z},\underline{E})-\frac{1}{2}\mathrm{Tr}|_{\ker
			D(\mathcal{F})}[g],
	\end{align}
	which is half of the equivariant APS $\eta$-invariant.
	
	If $\ker D(\mathcal{F})$ forms a vector bundle
	over $B$,
	the $L^2$-inner product on $\mathscr{E}$ induces a 
	$G$-invariant metric on $\ker D(\mathcal{F})$.
	Let $P^{\ker}: \mathscr{E}\to \ker D(\mathcal{F})$
	be the orthogonal projection associated with this $L^2$-inner
	product on $\mathscr{E}$. Let $\nabla^{\ker}:=P^{\ker}\circ
	\nabla^{\mathscr{E},u}\circ P^{\ker}$. Then
	$\underline{\ker D(\mathcal{F})}$ is an equivariant
	geometric triple over $B$. Furthermore, we have (cf. 
	\cite[(2.81)]{Liu17}, also \cite[(1.44), (1.49)]{LM22})
	\begin{align}\label{eq:3.73}
		d\tilde{\eta}_g(\mathcal{F})=\Bigg\{
		\begin{array}{ll}
			\int_{Z^g}\mathrm{Td}_g(\nabla^{TZ}, \nabla^{L})\mathrm{ch}_g
			(\underline{E})-\mathrm{ch}_g
			\left(\underline{\ker D(\mathcal{F})}\right),\hspace{10mm} 
			& \hbox{if $\mF$ is even;} \\
			\int_{Z^g}\mathrm{Td}_g(\nabla^{TZ}, \nabla^{L})\mathrm{ch}_g
			(\underline{E}),\hspace{10mm} & \hbox{if $\mF$ is odd.}
		\end{array}
		\Bigg.
	\end{align}
	For the definitions of equivariant characteristic forms
	$\mathrm{Td}_g(\cdot)$ and $\mathrm{ch}_g(\cdot)$,
	see e.g., \cite[(1.19), (1.21)]{LM20}.
	
	\
	
	In the following, we introduce another type of equivariant 
	$\eta$-forms by replacing the bundle condition by the 
	perturbation condition. 
	
	\begin{defn}\label{defn:1.03}
		A smooth family of equivariant
		bounded self-adjoint pseudodifferential operators
		 $\mA$ along the fiber 
		is called a \textbf{perturbation operator} associated
		with $\mathcal{F}$
		if $(D(\mathcal{F})+\mA)|_{Z_b}$ is invertible 
		for any $b\in B$, and, when $\mF$ is even,
		$\mA$ anti-commutes with the 
		$\Z_2$-grading of $\mathcal{S}(TZ, L_Z)\widehat{\otimes}E$.
	\end{defn}
	
Remark that when $\mF$ is odd, $D(\mF)$ here should be read as 
(\ref{eq:1.15a}). We can explain the perturbation operator
 rigorously using 
the orientation reversal in	Section \ref{s0103a}:
for a smooth family of equivariant
bounded self-adjoint pseudodifferential operators $\mA$ 
along the fiber, we say $\mA$ is a perturbation operator if 
$D(\mF_++\mF_-^{re})+\mA^{re}$ is invertible (cf. (\ref{eq:1.28a})), where
\begin{align}\label{eq:1.36b}
	\mA^{re}=(\mathrm{Id}\oplus \mU)\circ\mA\circ
	(\mathrm{Id}\oplus \mU^{-1}).
\end{align}
We note
that $\mathcal{A}$ does not preserve 
the $\mathbb{Z}_2$-grading in general for odd case,
although \cite[Prop. 3.3 (iii)]{Liu21} 
seems to suggest otherwise.
That $\mathbb{Z}_2$-graded preserving property 
should be understood as holding
only after the orientation reversal, as can be 
derived from the proof there. Accordingly, in the odd case, 
the equivariant spectral section of \cite[Def. 3.2]{Liu21} 
is likewise understood after the orientation reversal, 
and this is how it is used throughout \cite{Liu21}.
None of the results in \cite{Liu19,Liu21,Liu22} rely on the 
$\mathbb{Z}_2$-graded preserving property of $\mathcal{A}$ in the odd case.	 

	\begin{prop}\label{prop:1.4}\cite[Prop. 3.3 (i)]{Liu21}
	(cf. also \cite[Prop. 1]{MP97a} and \cite[
	Prop. 2]{MP97b})
	If there exists a perturbation operator associated with 
	$\mathcal{F}$, then $\ind(D(\mF))=0\in K_G^{\bullet}(B)$. Conversely, 
	if $\ind(D(\mF))=0\in K_G^{\bullet}(B)$ and at least one component
	of the fiber has nonzero dimension, there exists a 
	perturbation operator associated with $\mathcal{F}$.
\end{prop}
	
	For a perturbation operator $\mA$ associated with $\mF$,	
	we extend $\mA$ to $1\widehat{\otimes}\mA$ on 
	$\cC^{\infty}(B,\pi_Z^*\Lambda(T^*B)\widehat{\otimes}\cE)$
	as an odd element of the $\Z_2$-graded tensor product of 
	$\Z_2$-graded algebras such that
	\begin{align}\label{eq:1.23a}
		(\alpha\widehat{\otimes}1)(1\widehat{\otimes}\mA)=
		(-1)^{\deg \alpha}(1\widehat{\otimes}\mA)(\alpha\widehat{\otimes}1).
	\end{align} 
	We usually abbreviate $1\widehat{\otimes}\mA$ by $\mA$ when there is 
	no confusion.
	
	Let $\chi\in \cC^{\infty}(\R)$ be a cut-off function such that
	\begin{align}\label{eq:1.08}
		\chi(u)=
		\left\{
		\begin{aligned}
			&0,  &\hbox{if $u<1$;} \\
			&1,  &\hbox{if $u>2$.}
		\end{aligned}
		\right.
	\end{align}
	
	\begin{defn}\label{defn:1.02} \cite[Def. 3.12]{Liu21}
Assume that the $G$-action on $B$ is trivial. For $g\in G$, the equivariant 
$\eta$-form with perturbation operator $\mA$ is defined by
		\begin{multline}\label{eq:1.10}
			\tilde{\eta}_g(\mathcal{F},\mA):=
			-\int_0^{\infty}\left\{\psi_{\R\times B}\left.\widetilde{\tr}
			\right.\left[g\exp\left(-\left(
			\mathbb{B}_u+\sqrt{u}\chi(\sqrt{u})\mA+du\wedge\frac{\partial}{
				\partial u}\right)^2\right)\right]\right\}^{du}du
			\\
			\in \Omega^{\bullet}(B,\C).
		\end{multline}	 
		
	\end{defn}
	
	A priori, this definition depends on the choice of the cut-off function 
	$\chi$ in (\ref{eq:1.08}). However, choosing a different cut-off function 
	satisfying (\ref{eq:1.08}) changes the eta form by an exact 
	form \cite{Liu21}. Hence we  regard 
	\begin{align}
	\tilde{\eta}_g(\mathcal{F},\mA)\in \Omega^{\bullet}(B,\C)/\,
	d\Omega^{\bullet}(B,\C),
	\end{align}
	that is independent of the choice of the cut-off function.
	When $\mF$ is even (resp. odd), $\tilde{\eta}_g(\mF,\mA)\in
	\Omega^{\mathrm{odd}}(B,\C)/\,d\Omega^{\bullet}(B,\C)$ (resp. $\Omega^{
		\mathrm{even}}(B,\C)/\,d\Omega^{\bullet}(B,\C)$).
	
	If $B$ is a point and $\dim Z$ is odd, the projection
	$P^{\ker}$ is a perturbation operator. In this case, by 
	\cite[Remark 3.20]{Liu21}, $\tilde{\eta}_g(\mathcal{F},
	P^{\ker})=\bar{\eta}_g(\underline{TZ},
	\underline{L_Z},\underline{E})$, the equivariant reduced eta invariant.
	
	Moreover, as in (\ref{eq:3.73}), we have (cf. \cite[(3.66)]{Liu21})
	\begin{align}\label{eq:1.25a}
		d\,\widetilde{\eta}_g(\mathcal{F},\mA)=\int_{Z^g}
		\mathrm{Td}_g(\nabla^{TZ}, \nabla^{L})\mathrm{ch}_g
		(\underline{E}).
	\end{align}
	
	If $\mF$ is odd, under the orientation reversal,
	we have
    \begin{multline}\label{eq:a.12}
		\tilde{\eta}_g(\mathcal{F},\mA)
		=
		-\int_0^{\infty}\Bigg\{\psi_{\R\times B}\cdot
		\tr|_{\cE_+\oplus \cE_-^{re}}\Bigg[g\exp\Bigg(-\Bigg(
	(\mathrm{Id}\oplus \mU)\mathbb{B}_u(\mathrm{Id}\oplus \mU^{-1})
		\Bigg.\\
	\Bigg.+\sqrt{u}\,\chi(\sqrt{u})\mA^{re}+du\wedge\frac{\partial}{\partial u}
	\Bigg)^2\Bigg)\Bigg]\Bigg\}^{du}\,du.
	\end{multline}
	In fact, for a trace class operator $M$ on $\cE_+\oplus \cE_-$
	along the fiber, we have
	\begin{multline}\label{eq:a.13}
	\tr|_{\cE_+\oplus \cE_-^{re}}[(\mathrm{Id}\oplus \mU)M 
	(\mathrm{Id}\oplus \mU^{-1})]
	=\tr[M|_{\cE_+}]+\tr[\mU M\mU^{-1}|_{\cE_-^{re}}]
	\\
	=\tr[M|_{\cE_+}]+\tr[\mU M|_{\cE_-}\mU^{-1}].
	\end{multline}
	From (\ref{eq:a.03}), since $n$ is odd,
	\begin{align}\label{eq:a.14}
		\tr[\mU c_o(e_1)\cdots c_o(e_n)\mU^{-1}]=
		-\tr[c_{-o}(e_1)\cdots c_{-o}(e_n)].
	\end{align}
	Thus (\ref{eq:a.12}) follows directly from Definition \ref{defn:1.02}, 
	(\ref{eq:a.13}) and (\ref{eq:a.14}).

	Until now we introduce two types of equivariant $\eta$-forms
	under different conditions. If both conditions hold, 
	these two types of equivariant $\eta$-forms are related.
	The following proposition compares to 
	\cite[Lemma 3.14]{Liu21} and \cite[Prop. 2.3, 
	Lemma 3.5]{Liu22}.
	
	\begin{prop}\label{prop:1.5}
	Assume that the $G$-action on $B$ is trivial.
	If there exists a perturbation operator $\mA$ associated with
	$\mathcal{F}$ and $\ker D(\mathcal{F})$ forms a vector bundle
	over $B$, there exists $x_1\in K_G^1(B)$ or $x_2\in K_G^0(B)$
	such that for any $g\in G$,
	\begin{align}
	\tilde{\eta}_g(\mathcal{F},\mA)=\Bigg\{
	\begin{array}{ll}
	\tilde{\eta}_g(\mathcal{F})+\ch_g(x_1)\ \mathrm{mod}\ 
	d\,\Omega^{\bullet}(B),\	&  \hbox{if $\mF$ is even;} \\
		\tilde{\eta}_g(\mathcal{F})+\frac{1}{2}\ch_g(\underline{\ker
	D(\mathcal{F})})+\ch_g(x_2)\ \mathrm{mod}\ d\,\Omega^{\bullet}(B),
	\  &	\hbox{if $\mF$ is odd.}
		\end{array}
		\Bigg.
	\end{align}
	\end{prop}
	
	\
	
	Let $T_{\pi_Z}^{H'}W$ be another $G$-invariant horizontal
	subbundle of $TW$ and
	$g^{TZ'}$ be another $G$-invariant metric on $TZ$.
	Let $\nabla^{TZ'}$ be the connection on $TZ$ defined in 
	(\ref{eq:1.15}) associated with $(T_{\pi_Z}^{H'}W, g^{TZ'})$.
	Let $(h^{L_Z'}, \nabla^{L_Z'})$ and $(h^{E'}, \nabla^{E'})$ be 
	another group of metrics and connections on $L_Z$ and $E$
	such that
	$\underline{L_Z'}=(L_Z, h^{L_Z'}, \nabla^{L_Z'})$ and 
	$\underline{E'}=(E, h^{E'}, \nabla^{E'})$ are 
	equivariant geometric triples over $W$.
	Set $\mathcal{F}'=(W, o_Z, T_{\pi}^{H'}W, g^{TZ'},
	\underline{L_Z'},
	\underline{E'})$.
	
	If there exists a perturbation operator $\mA$ associated with 
	$\mathcal{F}$, then it implies
	 $\ind (D(\mathcal{F}))=0\in K_G^{\bullet}(B)$.
	Since the equivariant family index is a topological invariant,
	we have $\ind (D(\mathcal{F}'))=0\in K_G^{\bullet}(B)$ too. 
	Thus from Proposition \ref{prop:1.4}, if at least one component 
	of the fiber has nonzero dimension, the perturbation operator 
	associated with $\mathcal{F}'$ exists.
	The following proposition was established in \cite[Th. 1.2]{Liu21}, 
	which is a generalization of \cite[Th. 0.1]{DZ98}. 
	
	\begin{prop}\label{prop:2.01}
		Let $\mathcal{A}$ and $\mA'$ be perturbation operators
		associated with $\mathcal{F}$ and $\mathcal{F}'$ respectively.
		Then there exists $x\in K_G^{\bullet}(B)$ such that
		for any $g\in G$, modulo exact forms on $B$, 
		\begin{multline}\label{D07}
			\tilde{\eta}_{g}(\mathcal{F}', \mA')-\tilde{\eta}_{g}(
			\mathcal{F}, \mA)=\int_{Z^g}\widetilde{\mathrm{Td}}_g( 
			\nabla^{TZ}, \nabla^{L_Z}, \nabla^{TZ'}, \nabla^{L_Z'})
			\, \ch_g(\underline{E})
			\\
			+\int_{Z^g}\mathrm{Td}_g(\nabla^{TZ'}, \nabla^{L_Z'} )\, 
			\widetilde{\ch}_g(\underline{E}, \underline{E'})
			+\ch_g(x).
		\end{multline} 
		For the definitions of equivariant Chern-Simons form
		$\widetilde{\mathrm{Td}}_g(\cdot)$ and $\widetilde{\ch}_g
		(\cdot)$, see e.g., \cite[(1.29)]{LM20}.
	\end{prop}
	
	
	Conversely, we have the following proposition, which is 
	the $g$-version of \cite[(4.29), (4.30)]{Liu21}. 	
	
	\begin{prop}\label{prop:1.11}
	For any $x\in K_G^{\bullet}(B)$, there exists an equivariant
	geometric family $\mF$ and perturbation operators $\mA_1$ and $
	\mA_2$, such that for any $g\in G$,
	\begin{align}
	\tilde{\eta}_{g}(\mathcal{F}, \mA_2)-\tilde{\eta}_{g}(\mathcal{F}, \mA_1)
	=\ch_g(x)\in \Omega^{\bullet}(B)/\,d\Omega^{\bullet}(B).
	\end{align}
	\end{prop}


	\subsection{Equivariant infinitesimal $\eta$-forms}\label{s0103}
	
	In this subsection, we
	introduce the equivariant infinitesimal $\eta$-form with 
	perturbation and
	 use the notations in Section
	\ref{s0104a}. We assume that 
	the $G$-action on $B$ is trivial.
	
	Let $\mathcal{F}$ be an equivariant geometric family 
	defined in (\ref{eq:1.04}). 
	Let $|\cdot|$ be an $\mathrm{Ad}_G$-invariant norm on the Lie algebra
	$\mathfrak{g}$ of $G$.
	For $K\in \mathfrak{g}$,
	let 
	$K^Z(x)=\left.\frac{\partial}{\partial t}\right|_{t=0}
	e^{tK}\cdot x\quad \text{for}\ x\in W$
	be the induced vector field 
	on $W$ and 
	let $\mL_{K}$ be 
	the corresponding Lie derivative given by
	$\mL_Ks=\left.\frac{\partial}{\partial t}\right|_{t=0}
	\left(e^{-tK}. s\right)$
	for any smooth section $s$ of a vector bundle over $W$.
	
	Let $Z(g)\subset G$ be the centralizer of $g\in G$ with Lie
	algebra $\mathfrak{z}(g)$. Then in the sense of the adjoint
	action, $\mathfrak{z}(g)=\{K\in \mathfrak{g}: 
	\mathrm{Ad}_gK=K \}$.
	Set
	\begin{align}
		\mathbb{B}_{K,u}'=\mathbb{B}_u+\sqrt{u}\chi(\sqrt{u})\mA+
		\frac{c(K^Z)}{4\sqrt{u}}.
	\end{align}
	
	\begin{defn}\label{defn:1.07} 
		Assume that there exists a perturbation operator $\mA$ 
		associated with $\mathcal{F}$.
		For $g\in G$, $K\in \mathfrak{z}(g)$ with $|K|$ small enough, 
		the equivariant infinitesimal 
		$\eta$-form with perturbation operator $\mA$ is defined by
		\begin{align}\label{eq:1.33}
		\tilde{\eta}_{g,K}(\mathcal{F},\mA):=
		-\int_0^{\infty}\left\{\psi_{\R\times B}\left.\widetilde{\tr}
		\right.\left[g\exp\left(-\left(\mathbb{B}_{K,u}'+du\wedge\frac{
			\partial}{\partial u}\right)^2-\mathcal{L}_K\right)\right]
		\right\}^{du}du.
		\end{align}	 
	\end{defn}
	
	The following proposition is a parallel version of 
	\cite[Th. 0.1]{LM22}. 
	
	\begin{prop}\label{prop:1.08} 
	For any $g\in G$, there exists $\beta>0$ such that if
	$K\in \mathfrak{z}(g)$ with $|K|<\beta$,  
	$\tilde{\eta}_{g,K}(\mathcal{F},\mA)$ is well-defined, and
		\begin{align}
			d\, \tilde{\eta}_{g,K}(\mathcal{F},\mA)=
			\int_{Z^g}\mathrm{Td}_{g,K}(\nabla^{TZ},\nabla^{L_Z})\, 
			\ch_{g,K}(\underline{E}).
		\end{align}
		Moreover, for fixed $K\in \mathfrak{z}(g)$, 
		$\wi{\eta}_{g,zK}(\mathcal{F},\mA)$ is analytic in
		$z\in \C$ for $|zK|<\beta$, that is,
		\begin{align}
			\frac{\partial}{\partial \bar{z}}
			\wi{\eta}_{g,zK}(\mathcal{F},\mA)=0.	
		\end{align} 
		For the definitions of $\mathrm{Td}_{g,K}(\cdot)$ and 
		$\ch_{g,K}(\cdot)$, see e.g., \cite[(1.49), (1.51)]{LM20}.
		Just like 
		$\wi{\eta}_{g}(\mathcal{F},\mA)$ in Definition \ref{defn:1.02},
		if we change another cut-off function satisfying (\ref{eq:1.08}), 
		$\wi{\eta}_{g,zK}(\mathcal{F},\mA)$ also differs from 
		its original by an exact form. We therefore also 
		regard $\wi{\eta}_{g,zK}(\mathcal{F},\mA)$
		as an element of $\Omega^{\bullet}(B,\C)/\,d\Omega^{\bullet}(B,\C)$,
		that is independent of the choice of the cut-off function.	 
	\end{prop}
	\begin{proof}
	Comparing with the proof of \cite[Th. 2.2]{LM22}, 
	for a fixed cut-off function we 
	only need to prove that for $u\in [1,2]$,
	\begin{align}\label{eq:1.36a}
		\left|\left\{\left.\widetilde{\tr}\right.\left[g
		\exp\left(-\left(\mathbb{B}_{K,u}'+du\wedge\frac{\partial}{
	\partial u}\right)^2-\mathcal{L}_K\right)\right]\right\}^{du} \right|
	\end{align}
		is uniformly bounded. This can be proved following the same 
		argument in \cite[\S 5.1]{LM22} by replacing the weighted 
		Sobolev norms there by the usual Sobolev norms along the fiber. 
		Remark that whenever we use the finite propagation speed
		method in \cite{LM22}, we have
		the restriction $t\in (0,1]$. So although 
		$\mA$ is a pseudodifferential operator which
		may destroy the finite propagation speed, 
		all localization argument in \cite{LM22} holds for 
		$\mathbb{B}_{K,u}'$ since $\chi(t)=0$
	for $t\in (0,1]$\footnote{We thank Prof. Shu Shen for reminding us 
	that pseudodifferential operators may affect the validity of the 
	localization.}.
		
	For the independence of the cut-off function, the argument is
	 the same as for $\wi{\eta}_{g}(\mathcal{F},\mA)$ in \cite{Liu21}. 
		For completeness, we outline the proof as follows.
		Let $\chi_0$ and $\chi_1$ be two cut-off functions satisfying
		(\ref{eq:1.08}). Set 
		\begin{align}
			\mathbb{B}_{K,u}'(s):=\mathbb{B}_u+\sqrt{u}\left((1-s)\chi_0+
			s\chi_1\right)(\sqrt{u})\mA+
			\frac{c(K^Z)}{4\sqrt{u}}
		\end{align}
		and 
		\begin{align}
		r(u,s):=\psi_{\R\times B}\left.\widetilde{\tr}
		\right.\left[g\exp\left(-\left(\mathbb{B}_{K,u}'(s)+du\wedge\frac{
			\partial}{\partial u}
		+ds\wedge\frac{
			\partial}{\partial s}\right)^2-\mathcal{L}_K\right)\right].
		\end{align}
		Then $\wi{\eta}_{g,zK,i}(\mathcal{F},\mA)=-\int_0^{\infty}
		r(u,s)^{du}|_{s=i}\, du$ is the equivariant infinitesimal 		
		$\eta$-form with perturbation $\mA$	associated with $\chi_i$
		for $i=0, 1$. Here we write
		\begin{align}
		r(u,s)=r(u,s)^0+du\wedge r(u,s)^{du}+ds\wedge r(u,s)^{ds}+du\wedge
		ds\wedge r(u,s)^{duds}
		\end{align}
		where $r(u,s)^{\bullet}\in \Omega^{\bullet}(B,\C)$.
		Since for $u\notin [1,2]$, $r(u,s)$ does not depend on $s$,
		we see that $r(u,s)^{ds}=r(u,s)^{duds}=0$ for $u\notin [1,2]$.
		Following the same way as \cite[(2.95)]{Liu17}, we have
	\begin{align}
	\wi{\eta}_{g,zK,1}(\mathcal{F},\mA)-\wi{\eta}_{g,zK,0}(\mathcal{F},\mA)
		=d \left(\int_{u\in [1,2]}\int_{s\in [0,1]} r(u,s)^{duds}du
		ds\right).
	\end{align}
	\end{proof}
	%
	
	Let $\vartheta_K\in T^*Z$ be the $1$-form which is dual to 
	$K^Z$ by the metric $g^{TZ}$. 
	We identify $\vartheta_K$ to a 1-form on $W$ which vanishes 
	on $T_{\pi_Z}^HW$.
	The following proposition is the parallel version of the 
	main result in \cite{LM22}.

	\begin{prop}\label{prop:1.09}	
		Assume that there exists a perturbation operator $\mA$ 
		associated with $\mathcal{F}$.
		For $g\in G$ and $K_{0}\in \mathfrak{z}(g)$, 
		there exists  $\beta>0$ such that for any $K=tK_{0}$,  $K\neq 0$ 
		and $-\beta<t<\beta$,
		modulo exact forms on $B$, we have
		\begin{align}\label{eq:1.35a}
			\tilde{\eta}_{g,K}(\mathcal{F}, \mA)=\tilde{\eta}_{ge^K}(
			\mathcal{F},\mA)+\mM_{g,K},
		\end{align}
		where $\mM_{g,K}$ is a well-defined integral defined by
		\begin{align}
			\mM_{g,K}=-\int_0^{+\infty}\int_{Z^g}\frac{\vartheta_K}{2 i\pi v}
			\exp\left(\frac{d\vartheta_K-2i\pi |K^Z|^2}{2 i\pi v}
			\right)\mathrm{Td}_{g,K}(\nabla^{TZ}, \nabla^{L_Z})
			\ch_{g,K}(\underline{E})\frac{dv}{v},
		\end{align}
		and $t^{\lfloor(\dim W^g+1)/2\rfloor}
		\mM_{g,tK_0}$ is real analytic in $t\in \R$, $|t|<\beta$. 
		Moreover, we have
		\begin{align}\label{eq:0.12}
			d\mM_{g,K}=\int_{Z^{g}}
			\mathrm{Td}_{g,K}(\nabla^{TZ}, \nabla^{L_Z})
			\ch_{g,K}(\underline{E})-\int_{Z^{ge^K}}
			\mathrm{Td}_{ge^{K}}(\nabla^{TZ}, \nabla^{L_Z})
			\ch_{ge^K}(\underline{E}).
		\end{align}	
		In particular, if $G=\mathbb{T}^m$, for $g\in 
		\mathbb{T}^m\backslash A_{\mathbb{T}}$, we have $\mM_{g,K}=0$.
	\end{prop}
	\begin{proof}
		We replace $\mathbb{B}_t$ in \cite[\S 4]{LM22} by 
		$\mathbb{B}_t+\sqrt{t}\chi(\sqrt{t})\mA$. Then \cite[Th. 4.3]{LM22} 
		in this version naturally holds. 
		Using the argument in 
		\cite[\S 5.11]{LM22}, similarly as the proof of Proposition 
		\ref{prop:1.08}, we can obtain \cite[Th. 4.2]{LM22} in 
		this version. Following the process in \cite[\S 4.2]{LM22}, 
		we obtain (\ref{eq:1.35a}).
		The equation (\ref{eq:0.12}) is from \cite[Prop. 3.2]{LM22}.
		
		If $G=\mathbb{T}^m$, $\mathfrak{z}(g)=\mathfrak{g}$.
		If $g\in \mathbb{T}^m\backslash A_{\mathbb{T}}$, we have
		$W^g=W^{\mathbb{T}}$. Thus for $g\in \mathbb{T}^m\backslash 
		A_{\mathbb{T}}$, $\vartheta_K=0$ on $W^g$ for any 
		$K\in \mathfrak{g}$. So $\mM_{g,K}=0$. 
		
		The proof of Proposition \ref{prop:1.09} is complete.		
	\end{proof}
	
	
	\subsection{Embedding formula for equivariant $\eta$-forms}\label{s0104}
	
	Let $W$ be a $G$-equivariant closed Riemannian manifold.
	For a finite dimensional complex $G$-representation $M$, 
	we consider the $G$-action on $W\times M$ given by 
	\begin{align}\label{eq:2.86b} 
		g(b,u)=(gb,gu), \quad b\in W,\ u\in M.
	\end{align}
	Thus $W\times M\rightarrow W$ is an equivariant vector bundle over $W$.  
	We denote this equivariant vector bundle by $E_M$. 
	Similarly, for finite dimensional complex virtual $G$-representation
	$M=M_+-M_-$, we obtain an equivariant $\mathbb{Z}_2$-graded 
	vector bundle $E_M=E_{M_+}\oplus E_{M_-}$.
	We denote by
	\begin{align}\label{eq:2.43}
		\underline{M}=(E_M, h^M, d)
	\end{align}
	the equivariant geometric
	triple over $W$ with trivial metric $h^M$ and trivial connection $d$. 
	
	In the following, we take $G=\mathbb{T}^m$ and use the notations
	in Section \ref{s0101}.
	Let $\iota: W^{\mathbb{T}}\hookrightarrow W$ be the 
	equivariant Riemannian embedding. We consider the normal bundle
	$N_{\alpha}$ of $W_{\alpha}^{\mathbb{T}}$ in $W$ as an equivariant
	complex vector bundle and denote by $N_{\alpha}^{\R}$ the 
	underlying real vector bundle of $N_{\alpha}$. Let 
	$C(N_{\alpha}^{\R})$ be the Clifford algebra bundle of 
	$(N_{\alpha}, g^{N^{\R}})$, where $g^{N^{\R}}$ is the induced 
	metric on $N_{\alpha}^{\R}$. Then $\Lambda(\overline{N}_{\alpha}^*)$ 
	is a $C(N_{\alpha}^{\R})$-Clifford module. Namely, if 
	$u\in N_{\alpha}$, let $u^*\in \overline{N}_{\alpha}^*$ be 
	the metric dual of $u$. The Clifford action on  
	$\Lambda(\overline{N}_{\alpha}^*)$ is defined by
	\begin{align}\label{eq:1.37a}
		c(u)=\sqrt{2}u^*\wedge,\quad c(\bar{u})=-\sqrt{2}i_{\bar{u}},
	\end{align}
	for any $u\in N_{\alpha}$. Here $i_{\cdot}$ is the interior
	product on differential forms.
	
	For an equivariant ungraded geometric triple $\underline{\mu}$
	over $W^{\mathbb{T}}$, from the arguments in
	\cite[\S 1b]{BZ93}, \cite[\S 3.3]{Liu19} 
	and \cite[\S 1.4]{LM20}, we can construct the geometric version
	of the equivariant Atiyah-Hirzebruch direct image of 
	$\underline{\mu}$ as an equivariant geometric triple
	$\underline{\xi}=\iota_*\underline{\mu}$ over $W$.
	For the completeness, we repeat the construction explicitly
	as follows.
	
	For any $\delta>0$, set $\mU_{\alpha,\delta}=\{(y,Z)\in 
	N_{\alpha}^{\R}: y\in W_{\alpha}^{\mathbb{T}}, |Z|<\delta \}$. 
	Then there exists $\varepsilon_0>0$ such that the exponential map 
	$(y,Z)\in \mU_{\alpha,2\varepsilon_0}
	\to \exp_y^W(Z)$ is a diffeomorphism between $\mU_{\alpha,
		2\varepsilon_0}$ and an open $\mathbb{T}^m$-equivariant 
	tubular neighbourhood of $W_{\alpha}^{\mathbb{T}}$ in $W$.
	Without confusion, we will regard $\mU_{\alpha, 2\varepsilon_0}$
	as this neighbourhood of $W_{\alpha}^{\mathbb{T}}$ in $W$
	via this identification. We choose $\varepsilon_0>0$ small
	enough such that for any $\alpha\neq \beta\in \mathfrak{B}$,
	$\overline{\mU_{\alpha,
			2\varepsilon_0}}\cap \overline{\mU_{\beta,
			2\varepsilon_0}}=\emptyset$.
	
	For $Z\in N_{\alpha}^{\R}$, let $\tilde{c}(Z)\in \End(\Lambda(
	N_{\alpha}^*))$ be the transpose of $c(Z)$ on $\Lambda(
	\overline{N}_{\alpha}^*)$ in (\ref{eq:1.37a}). 
	In particular, for $u\in N_{\alpha}$,
	$\tilde{c}(u)=\sqrt{2}i_u$, $\tilde{c}(\bar{u})
	=-\sqrt{2}\bar{u}^*\wedge$.
	For the projection 
	$\pi_{\alpha}:N_{\alpha}\to W_{\alpha}^{\mathbb{T}}$, let 
	$\pi_{\alpha}^*(\Lambda(N_{\alpha}^*))$ be the pull back bundle
	of $\Lambda(N_{\alpha}^*)$ over $N_{\alpha}$. Let $\tilde{c}(Z)
	:\pi_{\alpha}^*(\Lambda^{\mathrm{even/odd}}(N_{\alpha}^*))|_Z
	\to \pi_{\alpha}^*(\Lambda^{\mathrm{odd/even}}(N_{\alpha}^*))|_Z$
	denote the corresponding pull back isomorphism at $Z\in 
	N_{\alpha}^{\R}$ for $Z\neq 0$.
	
	As $\mathbb{T}^m$ acts trivially on $W_{\alpha}^{\mathbb{T}}$,
	applying \cite[Cor. I.9.9]{LM89} for each weight part,
	there exists a $\mathbb{T}^m$-equivariant complex vector 
	bundle $F_{\alpha}$ over $W_{\alpha}^{\mathbb{T}}$ such that
	$\Lambda^{\mathrm{even}}(N_{\alpha}^*)\otimes \mu_{\alpha}
	\oplus F_{\alpha}$ with $\mu_{\alpha}=\mu|_{W_{\alpha}^{
			\mathbb{T}}}$, is a trivializable
	$\mathbb{T}^m$-equivariant complex vector bundle over 
	$W_{\alpha}^{\mathbb{T}}$ with the $\mathbb{T}^m$-equivariant
	trivialization map 
	$\varphi_{\alpha}: E_{M_{\alpha}}|_{W_{\alpha}^{\mathbb{T}}}
	\to \Lambda^{\mathrm{even}}(N_{\alpha}^*)\otimes \mu_{\alpha}
	\oplus F_{\alpha}$, where $M_{\alpha}$ is a finite dimensional
	$\mathbb{T}^m$-representation.
	Then
	\begin{align}
		\mathrm{i}\cdot\tilde{c}(Z)\oplus \pi_{\alpha}^*\Id_{F_{\alpha}}:
		\pi_{\alpha}^*(\Lambda^{\mathrm{even}}(N_{\alpha}^*)
		\otimes \mu_{\alpha}
		\oplus F_{\alpha})|_Z\to \pi_{\alpha}^*(\Lambda^{\mathrm{odd}}
		(N_{\alpha}^*)\otimes \mu_{\alpha}
		\oplus F_{\alpha})|_Z
	\end{align}
	induces a $\mathbb{T}^m$-equivariant isomorphism between 
	two $\mathbb{T}^m$-equivariant complex vector bundle over 
	$\overline{\mU_{\alpha,2\varepsilon_0}}\backslash W_{\alpha}^{
		\mathbb{T}}$. By adding trivial $\mathbb{T}^m$-equivariant 
	complex vector bundles for the part
	$F_{\alpha}$, we can assume that $M_{\alpha}=M_{\beta}=M$
	for any $\alpha\neq \beta\in \mathfrak{B}$. Now the identification
	$\pi_{\alpha}^*(\Lambda^{\mathrm{even}}(N_{\alpha}^*)\otimes 
	\mu_{\alpha}\oplus F_{\alpha})$ (resp. $\pi_{\alpha}^*(
	\Lambda^{\mathrm{odd}}(N_{\alpha}^*)\otimes \mu_{\alpha}
	\oplus F_{\alpha})$) with $(W\backslash \cup_{\alpha}\mU_{
		\alpha,\varepsilon_0})\times M$ on $\mU_{\alpha,2\varepsilon_0}
	\backslash \mU_{\alpha, \varepsilon_0}$ via the map $\varphi_{\alpha}$
	(resp. $(\mathrm{i}\cdot\tilde{c}(Z)\oplus 
	\pi_{\alpha}^*\Id_{F_{\alpha}})\circ \varphi_{\alpha}$) defines
	a $\mathbb{T}^m$-equivariant vector bundle $\xi_{+}$ 
	(resp. $\xi_-$) over $W$. Then the direct image of $\mu$ by $\iota$
	is given by 
	\begin{align}
		\iota_*\mu=\xi_+-\xi_-.
	\end{align}
	From the construction, we have
	\begin{align}\label{eq:1.42}
		\xi_{\pm}|_{W_{\alpha}^{\mathbb{T}}}=
		\Lambda^{\mathrm{even/odd}}(N_{\alpha}^*)\otimes \mu_{\alpha}
		\oplus F_{\alpha}.
	\end{align}
	
	Take $\underline{F_{\alpha}}=(F_{\alpha}, h^{F_{\alpha}},
	\nabla^{F_{\alpha}})$ as an equivariant ungraded geometric
	triple over $W_{\alpha}^{\mathbb{T}}$ associated with 
	$F_{\alpha}$. By a partition of unity argument, we can 
	construct a $\mathbb{T}^m$-invariant metric $h^{\xi}=h^{\xi_+}
	\oplus h^{\xi_-}$ and a $\mathbb{T}^m$-invariant $\Z_2$-graded 
	Hermitian
	connection $\nabla^{\xi}=\nabla^{\xi_+}\oplus \nabla^{\xi_-}$
	on $\xi=\xi_+\oplus \xi_-$, such that
	\begin{align}\label{eq:1.43}
		\underline{\xi_{\pm}}|_{W_{\alpha}^{\mathbb{T}}}=
		\underline{\Lambda^{\mathrm{even/odd}}(N_{\alpha}^*)}
		\otimes \underline{\mu_{\alpha}}\oplus \underline{F_{\alpha}}.
	\end{align}
	Here $\underline{\Lambda^{\mathrm{even/odd}}(N_{\alpha}^*)}$
	is the equivariant geometric triple with induced metrics and 
	connections. We denote by 
	\begin{align}
		\iota_*\underline{\mu}=\underline{\xi}
	\end{align}
	the geometric version of the
	direct image of $\underline{\mu}$.
	
	Remark that even if without the geometric structure, the 
	construction of $\xi$ is not unique.
	
	\begin{lemma}\label{lemma:1.11}
		Let $\xi_0$ and $\xi_1$ be two different
		representations of $\iota_*\mu$ such that for $i=0,1$,
		\begin{align}\label{eq:1.45}
			\xi_{i,\pm}|_{W_{\alpha}^{\mathbb{T}}}=
			\Lambda^{\mathrm{even/odd}}(N_{\alpha}^*)
			\otimes \mu_{\alpha}\oplus F_{i,\alpha}.
		\end{align}
		Then there exist finite dimensional $\mathbb{T}^m$-representations
		$M_0$, $M_1$ and isomorphisms
		\begin{align}\label{eq:1.39a}
		\Phi_+: \xi_{0,+}\oplus E_{M_0}\to \xi_{1,+}\oplus E_{M_1},
		\quad \Phi_-: \xi_{0,-}\oplus E_{M_0}\to \xi_{1,-}\oplus E_{M_1}
		\end{align}
		of equivariant complex vector bundles such that
		for any $\alpha\in \mathfrak{B}$,
		\begin{multline}\label{eq:1.47}
			\Phi_{\pm}|_{W_{\alpha}^{\mathbb{T}}}=
			\Id\oplus \Psi_{\alpha}: \Lambda^{\mathrm{even/odd}}(N_{\alpha}^*)
			\otimes \mu_{\alpha}\oplus (F_{0,\alpha}\oplus E_{M_0}|_{
				W_{\alpha}^{\mathbb{T}}})
			\\
			\longrightarrow \Lambda^{\mathrm{even/odd}}(N_{\alpha}^*)
			\otimes \mu_{\alpha}\oplus (F_{1,\alpha}\oplus E_{M_1}|_{
				W_{\alpha}^{\mathbb{T}}}),
		\end{multline}
	where $\Psi_{\alpha}:F_{0,\alpha}\oplus E_{M_0}|_{W_{\alpha}^{
	\mathbb{T}}}\to F_{1,\alpha}\oplus E_{M_1}|_{W_{\alpha}^{\mathbb{T}}}$
	is an isomorphism of equivariant complex vector bundles over 
	$W_{\alpha}^{\mathbb{T}}$.
	\end{lemma}
	\begin{proof}
		From the construction of $\iota_*\mu$, there exist 
		finite dimensional $\mathbb{T}^m$-representations
		$M_0'$ and $M_1'$ such that $\xi_{i,+}$ is equivariant isomorphic
		to $E_{M_i'}$ for $i=0,1$. Set
		\begin{align}\label{eq:1.48}
			\begin{split}
				&\xi_{0,+}':=\xi_{0,+}\oplus 
				\xi_{1,+},\quad \xi_{0,-}':=\xi_{0,-}\oplus 
				\xi_{1,+},
				\\
				&\xi_{1,+}':=\xi_{1,+}\oplus 
				\xi_{0,+},\quad \xi_{1,-}':=\xi_{1,-}\oplus 
				\xi_{0,+}.
			\end{split}
		\end{align}
		Then the corresponding vector bundle $F_{0,\alpha}'$ 
		(resp. $F_{1,\alpha}'$) associated with $\xi_0'$ (resp. $\xi_1'$)
		is $F_{0,\alpha}\oplus E_{M_1'}|_{W_{\alpha}^{\mathbb{T}}}$
		(resp. $F_{1,\alpha}\oplus E_{M_0'}|_{W_{\alpha}^{\mathbb{T}}}$).
		Note that in this case, $F_{0,\alpha}'$ is not necessary isomorphic 
		to $F_{1,\alpha}'$ as equivariant vector bundles.
		From (\ref{eq:1.48}), we have
		\begin{align}
			\Lambda^{\mathrm{even}}(N_{\alpha}^*)
			\otimes \mu_{\alpha}\oplus F_{0,\alpha}'\simeq 
			\Lambda^{\mathrm{even}}(N_{\alpha}^*)
			\otimes \mu_{\alpha}\oplus F_{1,\alpha}'
		\end{align}
		as equivariant vector bundles over $W_{\alpha}^{\mathbb{T}}$.
		So we have $[F_{0,\alpha}']=[F_{1,\alpha}']\in 
		K_{\mathbb{T}^m}^0(W_{\alpha}^{\mathbb{T}})$,
		 which means that there exists
		finite dimensional $\mathbb{T}^m$-representation $M_{\alpha}'$
		such that $F_{0,\alpha}'\oplus E_{M_{\alpha}'}|_{W_{\alpha}^{
				\mathbb{T}}}\simeq F_{1,\alpha}'\oplus E_{M_{\alpha}'}|_{
			W_{\alpha}^{\mathbb{T}}}$.
		Set $M':=\oplus_{\alpha\in \mathfrak{B}}M_{\alpha}'$. Set
		\begin{align}
			\begin{split}
				&F_{0,\alpha}'':=F_{0,\alpha}'\oplus E_{M'}|_{W_{\alpha}^{
						\mathbb{T}}}=F_{0,\alpha}\oplus 
				E_{M_1'\oplus M'}|_{W_{\alpha}^{\mathbb{T}}},
				\\
				&F_{1,\alpha}'':=F_{1,\alpha}'\oplus E_{M'}|_{W_{\alpha}^{
						\mathbb{T}}}=F_{1,\alpha}\oplus 
				E_{M_0'\oplus M'}|_{W_{\alpha}^{\mathbb{T}}}.
			\end{split}
		\end{align}
		Then we have $F_{0,\alpha}''\simeq F_{1,\alpha}''$.
		From the construction of $\iota_*\underline{\mu}$,
		we obtain this lemma for $M_0=M_1'\oplus M'$ and $M_1=M_0'
		\oplus M'$.
		
		The proof of Lemma \ref{lemma:1.11} is complete.
	\end{proof}
	
	\begin{prop}\label{prop:1.13}
		The embedding $\iota: W^{\mathbb{T}}\to W$ induces a 
		group homomorphism $\iota_!:K_{\mathbb{T}^m}^0(W^{\mathbb{T}})
		\to K_{\mathbb{T}^m}^0(W)$ by $[\mu]\mapsto [\xi]:= 
		[\xi_+]-[\xi_-]$.
	\end{prop}
	\begin{proof}
		From (\ref{eq:1.39a}), if $\xi_0$ and $\xi_1$ are two 
		representations of $\iota_*\mu$, we have
		\begin{multline}\label{eq:1.51}
			[\xi_{0,+}]-[\xi_{0,-}]=[\xi_{0,+}\oplus E_{M_0}]-
			[\xi_{0,-}\oplus E_{M_0}]=[\xi_{1,+}\oplus E_{M_1}]
			-[\xi_{1,-}\oplus E_{M_1}]
			\\
			=[\xi_{1,+}]-[\xi_{1,-}]\in K_{\mathbb{T}^m}^0(W).
		\end{multline}
		
		For equivariant complex vector bundles $\mu$ and $\mu'$
		over $W^{\mathbb{T}}$, from the construction of the
		direct image, 
		$\iota_*\mu\oplus\iota_*\mu'$
		is also a representation of $\iota_*(\mu\oplus \mu')$.
		So from (\ref{eq:1.51}), we have
		\begin{align}\label{eq:1.52}
			[\iota_*(\mu\oplus \mu')]=[\iota_*\mu]+[\iota_*\mu']
			\in K_{\mathbb{T}^m}^0(W).
		\end{align}
		
		For equivariant complex vector bundles $\mu$ and $\mu'$
		over $W^{\mathbb{T}}$, if $[\mu]=[\mu']\in 
		K_{\mathbb{T}^m}^0(W^{\mathbb{T}})$,
		there exists equivariant complex vector bundle 
		$\mu''$ such that $\mu\oplus \mu''$ is isomorphic to $\mu'\oplus
		\mu''$. From (\ref{eq:1.51}) and (\ref{eq:1.52}), we have
		\begin{align}
			[\iota_*\mu]=[\iota_*(\mu\oplus \mu'')]-[\iota_*\mu'']
			=[\iota_*(\mu'\oplus \mu'')]-[\iota_*\mu'']
			=[\iota_*\mu']\in K^0_{\mathbb{T}^m}(W).
		\end{align}
		Thus $\iota$ induces a well-defined map 
		$\iota_!:K_{\mathbb{T}^m}^0(W^{\mathbb{T}})\to 
		K_{\mathbb{T}^m}^0(W)$. From (\ref{eq:1.52}),
		we see that $\iota_!$ is a group homomorphism.
		
		The proof of Proposition \ref{prop:1.13} is complete.
	\end{proof}

	In the reminder of this subsection we assume that
	$\pi_Z:W\to B$ is a $\mathbb{T}^m$-equivariant fiber bundle 
	over a connected closed manifold $B$
	with closed oriented fiber $Z$. 
	We use the notations in Section \ref{s0101} and
	assume that 
	the group action on $B$ is trivial. In this case $W^{\mathbb{T}}$ 
	is a $\mathbb{T}^m$-equivariant fiber bundle over  $B$
	with closed oriented fiber $Z^{\mathbb{T}}$.
	We assume that $TZ$ has an equivariant spin$^c$ structure.
	Let $\underline{\pi_Z}$ be a family of $\mathbb{T}^m$-equivariant
	geometric data described in (\ref{eq:0.02}).
	
	Let 
	$\pi_{Z_{\alpha}^{\mathbb{T}}}:=\pi_Z|_{W^{\mathbb{T}}_{\alpha}}:
	W_{\alpha}^{\mathbb{T}}\to B$ be the restriction of $\pi_Z$ on 
	$W^{\mathbb{T}}_{\alpha}$, which is 
	a $\mathbb{T}^m$-equivariant fiber bundle with
	fiber $Z_{\alpha}^{\mathbb{T}}$.
	Set 
	\begin{align}\label{eq:1.32a}
		L_{Z, \alpha}^{\mathbb{T}}:= L_Z|_{W^{\mathbb{T}}_{\alpha}}\otimes
		(\det N_{\alpha})^{-1}.
	\end{align}
	Then  the relative tangent bundle $TZ_{\alpha}^{\mathbb{T}}$
	has an equivariant spin$^c$ structure with equivariant complex line 
bundle $L_{Z, \alpha}^{\mathbb{T}}$ 
(cf. \cite[(1.33)]{LM20}, \cite[(1.47)]{LMZ03}). 
	Let $L_{Z}^{\mathbb{T}}$ be the line bundle over $W^{\mathbb{T}}$ 
	such that $L_{Z}^{\mathbb{T}}|_{W^{\mathbb{T}}_{\alpha}}
	=L_{Z,\alpha}^{\mathbb{T}}$.
	Let $\underline{N}$ and $\underline{L_{Z}^{\mathbb{T}}}$
	be the induced equivariant geometric triples over
	$W^{\mathbb{T}}_{}$. 
	Set 
	\begin{align}\label{eq:1.38b}
		\underline{\pi_{Z^\mathbb{T}}}=(W^{\mathbb{T}}, o_{Z^{\mathbb{T}}}, 
		T_{\pi_Z}^HW|_{W^{\mathbb{T}}},
		g^{TZ^{\mathbb{T}}},
		\underline{L_{Z}^{\mathbb{T}}}),
	\end{align}
	where $g^{TZ^{\mathbb{T}}}$
	is induced from $g^{TZ}$ and $o_{Z^{\mathbb{T}}}$ is the 
	orientation of $TZ^{\mathbb{T}}$ induced from that
	of $TZ$ (cf. \cite[Prop. 2.1]{Liu17}, \cite[Lemma 4.1]{LMZ002}).
	Note that $TZ|_{W^{\mathbb{T}}_{\alpha}}=TZ^{\mathbb{T}}_{\alpha}
	\oplus N_{\alpha}^{\R}$ as real vector bundles.
	Since $W^{\mathbb{T}}_{\alpha}$
	is totally geodesic in $W$, we have
	$\nabla^{TZ}|_{W_{\alpha}^{\mathbb{T}}}=
	\nabla^{TZ_{\alpha}^{\mathbb{T}}}\oplus \nabla^{N_{\alpha}}$,
	where  $\nabla^{TZ_{\alpha}^{\mathbb{T}}}$ is the connection on $
	TZ_{\alpha}^{\mathbb{T}}$ defined by (\ref{eq:1.15}).

	Note that in \cite[\S 3.3]{Liu19} and
	\cite[\S 1.4]{LM20}, we simply assume that $\mu$ is ungraded.
	Then $\xi$ is $\mathbb{Z}_2$-graded with $\xi=\xi_{+}\oplus \xi_{-}$.
	In fact, if $\mu=\mu_+\oplus \mu_-$ is $\Z_2$-graded,
	letting $\underline{\xi^+}=\iota_* \underline{\mu_+}$
	and  $\underline{\xi^-}=\iota_* \underline{\mu_-}$,
	we can regard $\underline{\xi}=\underline{\xi^+\oplus \xi^-}$ as 
	the direct image of $\underline{\mu}$, where the metrics and the 
	connections are defined by the direct sum. In this case
	the $\Z_2$-grading of $\xi$ is
	\begin{align}\label{eq:1.32}
		\xi_+=\xi_+^+\oplus \xi_-^-,\quad \xi_-=\xi_+^-\oplus \xi_-^+.
	\end{align}
	
	Set 
	\begin{align}\label{eq:1.80}
		\mathcal{F}_{W^{\mathbb{T}}}=(\underline{\pi_{Z^\mathbb{T}}},
		\underline{\mu}),
		\quad \mathcal{F}_{W}=(\underline{\pi_Z},
		\underline{\xi}).
	\end{align}
	By \cite[\S 3.3]{Liu19}, the geometric data of $\mathcal{F}_{
		W^{\mathbb{T}}}$ and $\mathcal{F}_{W}$ satisfy the fundamental
	assumption in \cite[\S 3.2]{Liu19}, which is the family version of
	\cite[(1.10)]{BZ93}.
	From \cite[Prop. 3.6]{Liu19}, 
	if $\mathrm{Ind}(D(\mathcal{F}_{
		W^{\mathbb{T}}}))=0\in K_{\mathbb{T}^m}^{\bullet}(B)$, then 
	$\mathrm{Ind}(D(\mathcal{F}_{
		W}))=0\in K_{\mathbb{T}^m}^{\bullet}(B)$.

	\begin{prop}\cite[Th. 3.7]{Liu19}\label{prop:1.10}
		Let $\mathcal{A}$ and $\mathcal{A}'$ be perturbation operators
		associated with $\mathcal{F}_W$ and $\mathcal{F}_{W^{\mathbb{T}}}$
		respectively. Then there exists $x\in K_{\mathbb{T}^m}^{\bullet}(B)$ 
		such that for any $g\in \mathbb{T}^m\backslash A_{\mathbb{T}}$,
		modulo exact forms on $B$,
		\begin{align}
			\widetilde{\eta}_g(\mathcal{F}_{W}, \mA)-
			\widetilde{\eta}_g(\mathcal{F}_{W^{\mathbb{T}}}, \mA')
			=\ch_g(x).
		\end{align}
	\end{prop}
	
	Remark that in \cite[Th. 3.7]{Liu19}, 
	$x\in K_{\mathbb{T}^m}^{\bullet}(B)$
	is an equivariant higher spectral flow depending on $W^g$
	for $g\in \mathbb{T}^m$. Since for $g\in \mathbb{T}^m\backslash
	A_{\mathbb{T}}$, $W^g=W^{\mathbb{T}}$ does not change, thus
	$x$ does not depend on $g\in \mathbb{T}^m\backslash
	A_{\mathbb{T}}$. From the same argument as 
	\cite[Remark 1.6]{LM20}, the equivariant
	Bismut-Zhang current in \cite[Th. 3.7]{Liu19} vanishes in 
	this case.
	
	Note that in \cite[Th. 3.7]{Liu19} we simply assume that
	$\mu$ is ungraded. But for $\Z_2$-graded case, by (\ref{eq:a.12}) 
	and (\ref{eq:1.32}), the proof is almost the same. 
	

	\section{Equivariant Differential K-theory}\label{s02}
	
	In this section, we introduce two models of $g$-equivariant
	differential K-theory for a fixed $g\in G$ and discuss some 
	properties of them. In Section \ref{s0201}, we recall the 
	$g$-equivariant extension of the Freed-Lott's model for 
	differential K$^0$-theory in \cite{LM20} and define the 
	pushforward map for $g$-equivariant differential 
	K$^0$-groups associated with the embedding from the fixed point 
	set of the torus action to the original manifold.
	In Section \ref{s0202}, we develop another $g$-equivariant 
	differential K-theory model based on
	equivalence classes of equivariant geometric families and 
	differential forms.
	We also explain the isomorphism between these two models
	for $g$-equivariant differential K$^0$-rings.
	In Section \ref{s0203},  we discuss the pushforward map along 
	the submersion in the
	geometric family framework.
	
	
	\subsection{$g$-equivariant differential K-theory: vector bundle 
		version}\label{s0201}
	
	Let $B$ be a closed manifold.
	Let $G$ be a compact Lie group, which acts on $B$ smoothly.
	In this section, we do not assume that the group action on $B$ 
	is trivial.
	
	The following definition is an equivariant extension of
	the definition of differential K-group in
	\cite[Def. 2.16]{FreedLott10}.
	
	\begin{defn}\cite[Def. 2.14]{LM20}\label{defn:2.12} 
		For $g\in G$, a cycle for the $g$-equivariant version of the 
		differential K-group
		of $B$ is a 
		pair $(\underline{E}, \phi)$ where 
		$\underline{E}=(E, h^E, \nabla^E)$
		is a $G$-equivariant ungraded geometric triple on $B$ and $\phi$
		is an element in $\Omega^{\mathrm{odd}}(B^g, \C)/\Im \, d$.
		We say two cycles $(\underline{E_1}, \phi_1)$ and $(\underline{E_2}, 
		\phi_2)$ are equivalent if there exist a $G$-equivariant ungraded
		geometric triple $\underline{E_3}=(E_3, h^{E_3}, \nabla^{E_3})$ 
		and a $G$-equivariant vector bundle isomorphism over $B$
		\begin{align}\label{eq:2.082}
			\Phi:E_1\oplus E_3\rightarrow E_2\oplus E_3
		\end{align}
		such that
		\begin{align}\label{eq:2.083}
		\widetilde{\ch}_g\left(\underline{E_1}\oplus\underline{E_3}, 
		\Phi^{*} \left(\underline{E_2}\oplus\underline{E_3}\right)\right)
		=\phi_2-\phi_1.
		\end{align} 
		We define the sum of the cycles in the natural way and denote by	
		$\widehat{K}_{FL,g}^0(B)$  
		the Grothendieck group of equivalence classes of cycles.
	\end{defn}
	
	Denote by $[\underline{E}, \phi]\in \widehat{K}^0_{FL,g}(B)$
	the equivalence class of a cycle $(\underline{E}, \phi)$.
	For $[\underline{E}, \phi], [\underline{F}, \psi]\in 
	\widehat{K}^0_{FL,g}(B)$, set
	\begin{align}\label{eq:2.084}
		[\underline{E}, \phi]\cup[\underline{F}, \psi]
		=[\underline{E}\otimes \underline{F}, 
		\ch_g(\underline{E})\wedge\psi+\phi\wedge\ch_g(\underline{F})
		-d\phi\wedge\psi  ].
	\end{align} 
	From (\ref{eq:2.084}), 
	the element $1:=[\underline{\C},0]$ is a unity. Here
	$\underline{\C}$ denotes the trivial line bundle over $B$ with 
	the trivial metric and connection.
	From \cite[Lemma 2.15]{LM20},
	the product (\ref{eq:2.084}) on $\widehat{K}^0_{FL,g}(B)$ is 
	well-defined, 
	associative and commutative. 
	Thus $(\widehat{K}^0_{FL,g}(B), +, \cup)$ is a commutative ring
	with unity $1$.
	
	By \cite[Prop. 3.1]{LM20}, we have the following exact sequence:
	\begin{align}\label{eq:1.11}
		K_G^{1}(B)\overset{\footnotesize{\ch_{g}}}{\longrightarrow} 
		\Omega^{\mathrm{odd}}(B^g, \C)/\Im\, d
		\overset{\footnotesize{a}}{\longrightarrow} \widehat{K}_{FL,g}^0(B) 
		\overset{\footnotesize{I}}{\longrightarrow} K_G^0(B)
		\longrightarrow 0.
	\end{align}
	Here the homomorphisms $a$ and $I$ are defined by 
	\begin{align}
		a(\phi)=[\underline{0}, -\phi],\quad I([\underline{E}, \phi])=[E],
	\end{align}
	for $\phi\in \Omega^{\mathrm{odd}}(B^g, \C)/\Im\, d$ and
	$[\underline{E}, \phi]\in \widehat{K}_{FL,g}^0(B)$.
	Remark that in \cite[Prop. 3.1]{LM20}, (\ref{eq:1.11}) is stated
	for $G=S^1$. But the proof there holds for any compact Lie group. 
	For the homomorphism 
	$a$, the definition here differs from that in \cite[(3.12)]{LM20} by 
	a minus sign, which does not affect the validity of (\ref{eq:1.11}). 
	We adopt this definition because we want (\ref{eq:3.052}) later on to 
	be a commutative diagram and (\ref{d033}) to carry no minus sign, 
	in agreement with the corresponding definitions 
	in \cite{Bunke2009,BunkeSchick13,Liu21}.
	
	We can also represent $\widehat{K}_{FL,g}^0(B)$ using 
	$\mathbb{Z}_2$-graded
	equivariant vector bundles. For $\mathbb{Z}_2$-graded 
	equivariant complex vector bundles
	$E=E_+\oplus E_-$ and $E'=E'_+\oplus E'_-$, 
	$[E_+]-[E_-]=[E'_+]-[E'_-]\in
	K^0_G(B)$ if and only if there exist equivariant complex 
	vector bundles
	$F$, $F'$ over $B$ and equivariant vector bundle isomorphisms  
	\begin{align}
	E_+\oplus F\simeq E'_+\oplus F',\quad E_-\oplus F\simeq E'_-\oplus F'.
	\end{align}
	
	\begin{prop}\label{prop:2.02}
	The $g$-equivariant differential K$^0$-group $\widehat{K}_{FL,g}^0(B)$ 
		is generated by the cycles $(\underline{E},
		\phi)$, where $\underline{E}$ is an equivariant geometric triple and 
		$\phi \in \Omega^{\mathrm{odd}}(B^g, \C)/\Im\, d$. Two cycles 
	$(\underline{E},\phi)$ and $(\underline{E'}, \phi')$ are equivalent if
		there exist geometric triples $\underline{F}$ and $\underline{F'}$ 
		such that $\underline{F_+}=\underline{F_-}$, $\underline{F'_+}=
		\underline{F'_-}$ and an equivariant $\mathbb{Z}_2$-graded
		vector bundle isomorphism 
		\begin{align}\label{eq:2.07}
			\Phi: E\oplus F\to E'\oplus F'
		\end{align}
		such that
		\begin{align}\label{eq:2.08}
			\widetilde{\ch}_g\left(\underline{E}\oplus\underline{F}, 
			\Phi^{*} 
			\left(\underline{E'}\oplus\underline{F'}\right)\right)
			=\phi_2-\phi_1.
		\end{align}	
	\end{prop}
	\begin{proof}
		From the proof of \cite[Prop. 3.1]{LM20},
		the group generated by the cycles $(\underline{E},
		\phi)$ satisfying the equivalence relations (\ref{eq:2.07}) and 
		(\ref{eq:2.08}) verifies the long exact sequence (\ref{eq:1.11})
		for $I(\underline{E},\phi)=[E_+]-[E_-]$. So it is isomorphic 
		to $\widehat{K}_{FL,g}^0(B)$ from the five lemma.
		
		The proof of Proposition \ref{prop:2.02} is complete.
	\end{proof}
	
	Using this $\mathbb{Z}_2$-graded version, denoting by 
	$E^{\mathrm{op}}$ the $\mathbb{Z}_2$-graded vector bundle
	reversing the $\mathbb{Z}_2$-grading of $E$, we have
	$-[\underline{E}, \phi]=[\underline{E^{\mathrm{op}}}, -\phi]$.
	So any element in $\widehat{K}_{FL,g}^0(B)$ can be represented
	as a $\mathbb{Z}_2$-graded cycle $(\underline{E},
	\phi)$.
	
	Let $M$ be a finite dimensional complex $G$-presentation.
	For the equivariant trivial vector bundle $E_M$ over $B$ 
	defined by (\ref{eq:2.86b}),
	$E\mapsto E_M\widehat{\otimes} E$ 
	endows the $G$-equivariant K-group $K_{G}^0(B)$
	with the structure of an $R(G)$-module.
	From the exact sequence in (\ref{eq:3.032}),
	we can equip a structure of an $R(G)$-module on $K_G^1(B)$.
	Let $\chi_M$ be the character for the representation $M$.
	From \eqref{eq:2.083}, since 
	$\ch_g(\underline{M})=\chi_{M}(g)$, constant on $B^g$,
	\begin{align}\label{eq:2.09}
		(\underline{E}, \phi)\mapsto (\underline{M}
		\widehat{\otimes} \underline{E},  
		\chi_{M}(g)\cdot\phi)
	\end{align}
	makes $\widehat{K}_{FL,g}^0(B)$ an 
	$R(G)$-module.  In particular, $\phi\mapsto \chi_{M}(g)\cdot\phi$
	makes $\Omega^{\bullet}(B^g,\C)/\mathrm{Im}\,d$ an $R(G)$-module.

	For $G=\mathbb{T}^m$, if $g\in \mathbb{T}^m$ such that 
	$B^g=B^{\mathbb{T}}$, the normal bundle $N$ has an equivariant 
	complex structure, thus an equivariant spin$^c$ structure.
	We can define the pushforward map for $\widehat{K}_{FL,g}^0$ 
	along the embedding 
	$\iota: B^{\mathbb{T}}\to B$ using the equivariant version of 
	the direct image (cf. \cite[Def. 4.14]{FreedLott10} and
	\cite[Def. 3.3]{LM20}) as follows. 
	We use the notations in Section \ref{s0104}.

	\begin{defn}\label{defn:2.02} For $g\in \mathbb{T}^m$ such 
		that $B^g=B^{\mathbb{T}}$,
		we define the push-forword map
		\begin{align}
			\hat{\iota}_{!}:\widehat{K}_{FL,g}^0(B^{\mathbb{T}})\to 
			\widehat{K}_{FL,g}^0(B)
		\end{align}
	by
	\begin{align}\label{eq:1.026}
	\hat{\iota}_{!}\left(\left(\underline{\mu},\phi\right)\right):=
	\left(\iota_*\underline{\mu}, 
	\ch_g\left(\underline{\Lambda^{}(N^*)}\right)\wedge\phi \right).
	\end{align}
	\end{defn}
	
	As remarked below Proposition \ref{prop:1.10}, the Bismut-Zhang current
	term in \cite[Def. 4.14]{FreedLott10} vanishes.
	The following proposition is the equivariant generalization of
	a special case in
	\cite[Prop. 4.18]{FreedLott10}. 
	
	\begin{prop}\label{prop:1.04}
		The pushforward map $\hat{\iota}_{!}:\widehat{K}_{FL,g}^0(
		B^{\mathbb{T}})\to \widehat{K}_{FL,g}^0(B)$ in Definition 
		\ref{defn:2.02} 
		is a well-defined homomorphism of $R(G)$-modules.
	\end{prop}
	\begin{proof}
		For a cycle $(\underline{\mu}, 0)$ of $\widehat{K}_{FL,g}^0
		(B^{\mathbb{T}})$ where $\underline{\mu}$ is ungraded,
		we claim that if $\underline{\xi_{0}}$ and $\underline{\xi_{1}}$ 
		are two geometric versions of the direct image of $\underline{\mu}$
		over $B$, then 
		\begin{align}\label{eq:2.11}
			[\underline{\xi_0}, 0]=[\underline{\xi_1}, 0]\in 
			\widehat{K}_{FL,g}^0(B).
		\end{align} 
		In fact, from Lemma \ref{lemma:1.11}, there exist 
		representations $M_0$, $M_1$ and isomorphisms
		\begin{align}
			\Phi_{+}:\xi_{0,+}\oplus E_{M_0}\to 
			\xi_{1,+}\oplus E_{M_1},\quad 
			\Phi_{-}:\xi_{0,-}\oplus E_{M_0}\to 
			\xi_{1,-}\oplus E_{M_1}
		\end{align}
		satisfying (\ref{eq:1.47}).
		Set
		\begin{align}
			\underline{\tilde{\xi}_{i,+}}:=\underline{\xi_{i,+}}
			\oplus \underline{M_i},\quad
			\underline{\tilde{\xi}_{i,-}}:=\underline{\xi_{i,-}}
			\oplus \underline{M_i}
		\end{align}
		for $i=0,1$. Note that by 
		(\ref{eq:1.43}), for $i=0,1$, we have
		\begin{align}
			\underline{\xi_{i,\pm}}|_{B_{\alpha}^{\mathbb{T}}}=
			\underline{\Lambda^{\mathrm{even/odd}}(N_{\alpha}^*)}
			\otimes \underline{\mu_{\alpha}}\oplus \underline{F_{i,\alpha}}.
		\end{align}
		Thus by (\ref{eq:1.47}), 
		\begin{multline}
			\widetilde{\ch}_g\left(\underline{\tilde{\xi}_0},
			\Phi^{-1}\underline{\tilde{\xi}_1}\right)=\widetilde{\ch}_g
			\left(\underline{F_{0,\alpha}}\oplus \underline{M_0}|_{
				B_{\alpha}^{\mathbb{T}}}, \Psi^{-1}
			(\underline{F_{1,\alpha}}\oplus \underline{M_1}|_{
				B_{\alpha}^{\mathbb{T}}})\right)
			\\
			-\widetilde{\ch}_g\left(\underline{F_{0,\alpha}}\oplus 
			\underline{M_0}|_{B_{\alpha}^{\mathbb{T}}}, \Psi^{-1}
			(\underline{F_{1,\alpha}}\oplus \underline{M_1}|_{
				B_{\alpha}^{\mathbb{T}}})\right)=0.
		\end{multline}
		From Proposition \ref{prop:2.02}, we obtain (\ref{eq:2.11}).
		
		For two ungraded geometric triples $\underline{\mu}$
		and $\underline{\mu'}$, $\iota_*\underline{\mu}\oplus\iota_*
		\underline{\mu'}$ is a geometric version of the direct image 
		of $\underline{\mu}\oplus\underline{\mu'}$.
		So from (\ref{eq:2.11}), as in (\ref{eq:1.52}), we have
		\begin{align}\label{eq:2.13}
			[\iota_*\underline{\mu},\phi]+
			[\iota_*\underline{\mu'},\phi']=	
			[\iota_*\underline{\mu}\oplus\iota_*
			\underline{\mu'},\phi+\phi']=
			[\iota_*(\underline{\mu}\oplus\underline{\mu'}), \phi
			+\phi']\in \widehat{K}_{FL,g}^0(B).
		\end{align}

		From Definition \ref{defn:2.12}, if $(\underline{\mu}, \phi)$ 
		is equivalent to $(\underline{\mu'}, \phi')$ in 
		$\widehat{K}_{FL,g}^0(B^{\mathbb{T}})$, 
		then there exists equivariant ungraded geometric 
		triple $\underline{\mu''}$ and an equivariant isomorphism 
		$\Phi: \mu\oplus \mu''\to \mu'\oplus \mu''$ over $B^{\mathbb{T}}$ 
		such that $	\widetilde{\ch}_g\left(\underline{\mu}\oplus
		\underline{\mu''}, \Phi^{-1} 
		\left(\underline{\mu'}\oplus\underline{\mu''}\right)\right)
		=\phi'-\phi.$
		From the construction of the direct image,
		we can choose $\underline{\xi}=\iota_*(\underline{\mu}\oplus 
		\underline{\mu''})$, $\underline{\xi'}=\iota_*(\underline{\mu'}\oplus 
		\underline{\mu''})$ and an equivariant isomorphism 
		$\Phi':\xi\to \xi'$ such that
		restricted on $B^{\mathbb{T}}_{\alpha}$, $\Phi'$ is identity on 
		$F_{\alpha}$ and equal to $\mathrm{Id}\otimes \Phi$ on 
		$\Lambda(N_{\alpha}^*)\otimes (\mu\oplus \mu'')$ associated with 
		(\ref{eq:1.42}). Thus we have $\widetilde{\ch}_g(\underline{\xi}, 
		(\Phi')^{-1}\underline{\xi'})
		=\ch_g(\underline{\Lambda(N^*)})\wedge(\phi'-\phi)$. 
		So from Definition \ref{defn:2.12} and (\ref{eq:2.11}),
		\begin{align}\label{eq:2.13a}
			[\iota_*(\underline{\mu}\oplus \underline{\mu''}), 
			\ch_g(\underline{\Lambda(N^*)})\wedge\phi]
			=[\iota_*(\underline{\mu'}\oplus \underline{\mu''}), \ch_g(
			\underline{\Lambda(N^*)})\wedge\phi']\in \widehat{K}_{FL,g}^0(B).
		\end{align}
		From (\ref{eq:2.13}) and (\ref{eq:2.13a}), we see that 
		$[\iota_*\underline{\mu}, \ch_g(\underline{\Lambda(N^*)})
		\wedge\phi]=[\iota_*\underline{\mu'}, \ch_g(\underline{
			\Lambda(N^*)})\wedge\phi']$. 
		So the map in Definition \ref{defn:2.02} is well-defined.
		
		By (\ref{eq:2.13}),  $\hat{\iota}_!$ is a group homomorphism.
		For a finite dimensional $G$-representation $M$.
		$\underline{M}\widehat{\otimes}\underline{\xi} $
		is a representation of $\iota_*(\underline{M}|_{B_{
				\alpha}^{\mathbb{T}}}\otimes \underline{\mu})$ by replacing 
		$\underline{F_{\alpha}}$ by 
		$\underline{M}|_{B_{\alpha}^{\mathbb{T}}}\otimes 
		\underline{F_{\alpha}}$ in the construction of the direct image. 
		So $\hat{\iota}_!$ is a homomorphism of $R(G)$-modules.
		
		The proof of Proposition \ref{prop:1.04} is complete.	
	\end{proof}


	\subsection{$g$-equivariant differential K-theory: geometric 
		family version}\label{s0202}
	In this subsection we use the notations and the conventions 
	in \cite{Liu21}.

	Let $\mF$ and $\mF'$ be two  equivariant geometric families over $B$.
	We mark the objects associated with $\mF'$ by $'$.
	We say 
	$\mF$ is isomorphic to $\mF'$ if there exist equivariant
	diffeomorphism
	$f:W\to W'$ with $\pi'\circ f=\pi$ and equivariant
	$\Z_2$-graded vector bundle isomorphism $f_E:E\to f^*E'$ such that 
	the orientation of the relative 
	tangent bundle, the equivariant spin$^c$ structure of the relative 
	tangent bundle, the horizontal subbundle and all geometric structures
	are preserved by these maps. 
	For the details, see \cite[Def. 2.2]{Liu21} (cf. also
	\cite[\S 2.1.7]{Bunke2009}).
	Recall that
	we say an equivariant geometric family $\mF$ is even (resp. odd) if
	for any connected component of the fiber, the dimension is even 
	(resp. odd).
	Let $\mathrm{F}_G^0(B)$ (resp. $\mathrm{F}_G^1(B)$) be the set of 
	equivalence classes of even (resp. odd) equivariant geometric families 
	over $B$. 
Let $\mathrm{F}_G^{\bullet}(B)=\mathrm{F}_G^0(B)\oplus \mathrm{F}_G^1(B)$.
	It is an additive abelian semigroup.
	
	Since isomorphic geometric families have the same family index,  
	the family index of the fiberwise Dirac operator 
	induces a semi-group homomorphism
	\begin{align}\label{eq:1.03}
		\ind: \mathrm{F}_G^{\bullet}(B)\to K_G^{\bullet}(B).
	\end{align}
	This homomorphism is surjective due to \cite[Prop. 2.6]{Liu21}
	(cf. also \cite[Lemma 2.18]{Bunke2009} and 
	\cite[\S 2.5.5]{BunkeSchick13}).
	
	The following definition extends \cite[Def. 4.1]{Liu21}, 
	\cite[Def. 2.1]{Bunke2009} and 
	\cite[Def. 2.4]{BunkeSchick13} to our $g$-version.
	
	\begin{defn}
		For $g\in G$, a family cycle for a $g$-equivariant 
		differential K-group over $B$ is a pair $(\mF, \phi)$, where 
	$\mF\in \mathrm{F}_G^{\bullet}(B)$ and 
	$\phi\in \Omega^{\bullet}(B^g,\C)/\Im \,d$.
		A family cycle $(\mF, \phi)$ is called even (resp. odd) 
		if $\mF$ is even (resp. odd) and $\phi\in \Omega^{\mathrm{odd}}(
		B^g,\C)/\Im \,d$ (resp. $\phi\in \Omega^{\mathrm{even}}(B^g,\C)/
		\Im \,d$).
		Two family cycles ($\mF, \phi$) and ($\mF', \phi'$) are called 
		isomorphic if $\mF$ and $\mF'$ are isomorphic
		and $\phi=\phi'$. Let $\widehat{\mathrm{IC}}_g^0(B)$ (resp. 
		$\widehat{\mathrm{IC}}_g^1(B)$) denote the set of isomorphic classes 
		of even (resp. odd) family cycles over $B$ with a natural abelian 
	semi-group structure: $(\mF,\phi)+(\mF',\phi')=(\mF+\mF',\phi+\phi')$.
	\end{defn}
	
	Remark that for two equivariant geometric families $\mF$ and $\mF'$,
	if $\ind(D(\mF))=\ind(D(\mF'))$, then 
	$\mF+\mF^{'\mathrm{op}}\in \ker (\ind)$.
	Thus from Proposition \ref{prop:1.4}, if $\mF+\mF^{'\mathrm{op}}$ 
	has non-zero dimensional fiber,
	there exists a perturbation operator $\mathcal{A}$ associated with
	$\mF+\mF^{'\mathrm{op}}$. Let $G'$ be the subgroup of $G$ generated
	by $g\in G$. Restricted on $B^g$, the geometric family 
	$(\mF+\mF^{'\mathrm{op}})|_{B^g}$ is $G'$-equivariant and the 
	$G'$-action
	on $B^g$ is trivial. Since $g\in G'$, for a perturbation
	operator $\mA$ associated with $\mF+\mF^{'\mathrm{op}}$, the 
	equivariant $\eta$-form
	$\widetilde{\eta}_g((\mF+\mF^{'\mathrm{op}})|_{B^g},
	\mA|_{B^g})$ is well-defined. Then our $g$-version of 
	\cite[Def. 4.3]{Liu21}, \cite[Def. 2.10]{Bunke2009} and 
	\cite[Def. 2.15]{BunkeSchick13} is follows.
	
	\begin{defn}\label{defn:2.03}
		We call two family cycles ($\mF, \phi$) and ($\mF', \phi'$) paired 
		if $\ind(D(\mF))=\ind(D(\mF'))$, and there exists
		a perturbation operator $\mA$ associated with 
		$\mF+\mF^{'\mathrm{op}}$ 
		such that
		\begin{align}
			\phi-\phi'= \widetilde{\eta}_g((\mF+\mF^{'\mathrm{op}})|_{B^g},
			\mA|_{B^g}).
		\end{align}
	\end{defn}
	
	The following lemma is the $g$-version of 
	\cite[Lemma 4.4]{Liu21} (cf. also 
	\cite[Lemmas 2.11, 2.12]{Bunke2009} and 
	\cite[Lemmas 2.16, 2.17]{BunkeSchick13}). For the reader's convenience,
     we state the detailed proof here using our notation.
	
	\begin{lemma}\label{lemma:2,7}
		The relation ``paired" 
		is symmetric, reflexive and compatible with the semigroup 
		structure on $\widehat{\mathrm{IC}}_g^*(B)$. 
	\end{lemma}
	\begin{proof}
		The compatible
		properties can be easily derived from Definitions
		\ref{defn:1.02} and \ref{defn:2.03}. 
		
		For the symmetric property, we only need to prove that
		if $\mA$ is a perturbation operator associated with
		$\mF$, then there exists a perturbation operator $\mA'$
		associated with $\mF^{\mathrm{op}}$ such that on $B^g$,
		\begin{align}
		\widetilde{\eta}_g(\mF^{\mathrm{op}},
		\mA')=-\widetilde{\eta}_g(\mF,
		\mA).
		\end{align}
	In fact, let $\mA^{\mathrm{op}}:\cE^{\mathrm{op}}\to \cE^{\mathrm{op}}$ 
	denote the operator which is 
	the same as $\mA$ on the underlying space and associated
	with the opposite $\Z_2$-grading. 
	By Definition \ref{defn:1.01} and the definition of the supertrace,
	we have
			\begin{align}\label{eq:2.22}
			\widetilde{\eta}_g(\mF^{\mathrm{op}},
			\mA^{\mathrm{op}})=-\widetilde{\eta}_g(\mF,
			\mA).
		\end{align}
		
		For the reflexive property,		
		We only need to prove the following statement: 
		for any equivariant geometric family $\mF$, there
		exists a perturbation operator $\mA$ associated with
		$\mF+\mF^{\mathrm{op}}$ such that on $B^g$,
		\begin{align}\label{eq:2.23}
			\tilde{\eta}_g(\mF+\mF^{\mathrm{op}}, \mA)=0.
		\end{align}
		From (\ref{eq:2.22}), if we can construct $\mA$ such that
		$\mA=\mA^{\mathrm{op}}$, then (\ref{eq:2.23}) holds.
		We will discuss the construction of such
		$\mA$ in three cases. 
		
		If $\mF$ is even and at least one component of the 
fiber has nonzero dimension, for the $\Z_2$-grading $\tau$ of 
$\mF$,
we have $\tau D(\mF)=-D(\mF)\tau$. On $\cE\oplus \cE$,
the Dirac operator and the $\Z_2$-grading are
\begin{align}
\mD:=D(\mF+\mF^{\mathrm{op}})=\left(\begin{array}{cc}
	D(\mF) & 0 \\
	0 & D(\mF)
\end{array}\right),\quad 
\hat{\tau}=\left(\begin{array}{cc}
	\tau & 0 \\
	0 & -\tau
\end{array}\right).
\end{align}
Set
\begin{align}
J=\left(\begin{array}{cc}
	0 & \tau \\
	\tau & 0
\end{array}\right).
\end{align}
Then $J \mD=-\mD J$.
We choose a real-valued even function $f\in \cC_0^{\infty}(\R)$
such that $f(0)\neq 0$. Then $\mA:=Jf(\mD)$ is
an equivariant self-adjoint smoothing operator along the fiber. 
Since $f$ is even, 
\begin{align}\label{eq:2.26a}
(\mD+\mA)^2=\mD^2+ \mD
J f(\mD)+J f(\mD)\mD
+\mA^2=\mD^2+f(\mD)^2
\end{align} 
is invertible and $\mA \hat{\tau}=-\hat{\tau}\mA$. 
So $\mA$ is a perturbation operator
associated with $\mF+\mF^{\mathrm{op}}$.
Let $S: \mF+\mF^{\mathrm{op}}\to (\mF+\mF^{\mathrm{op}})^{\mathrm{op}}
=\mF^{\mathrm{op}}+ \mF$, $(a,b)\mapsto 
(b,a)$ be the exchange map. Then
\begin{align}
\mA^{\mathrm{op}}=S\mA S=\mA.
\end{align}

If $\mF$ is odd, we use the notation in Section \ref{s0103a}.
Denote by
\begin{align}
	\begin{split}
&\mF_+=(\underline{\pi_Z}, \underline{E_+}\oplus 0), \quad
\ \mF_-=(\underline{\pi_Z}, 0\oplus\underline{E_-}), \quad
\quad\ \mF_-^+:=(\underline{\pi_Z}, \underline{E_-}\oplus 0), \\
&\mF_+^-:=(\underline{\pi_Z}, 0\oplus\underline{E_+}), \quad
\mF_+^{re}:=(\underline{-\pi_Z}, \underline{E_+}\oplus 0), \quad
\mF_-^{re}:=(\underline{-\pi_Z}, \underline{E_-}\oplus 0).
\end{split}
\end{align}
Then we have $\mF=\mF_++\mF_-$, $\mF^{\mathrm{op}}=\mF_-^++\mF_+^-$.
Denote by $\cE_{\pm}$ (resp. $\cE_{\pm}^{re}$) the infinite
dimensional bundle associated with $\mF_{\pm}$ (resp. 
$\mF_{\pm}^{re}$) over $B$.
Then the families of Dirac operators
\begin{align}
D(\mF_{\pm}):\cE_{\pm}\to \cE_{\pm},\quad
D(\mF_{\pm}^{re}):\cE_{\pm}^{re}\to \cE_{\pm}^{re}.
\end{align}
Denote by $\mU_{\pm}:\cE_{\pm}\to \cE_{\pm}^{re}$
the orientation reversing map in Section \ref{s0103a}.
Then by (\ref{eq:a.08}), we have
\begin{align}\label{eq:2.30a}
	\mU_+^{-1}D(\mF_+^{re})\,\mU_+=-D(\mF_+),\quad
		\mU_-^{-1}D(\mF_-^{re})\,\mU_-=-D(\mF_-).
\end{align} 
Denote by $\mU:=\mU_+\oplus \mU_-:\cE=\cE_+\oplus \cE_-\to \cE_+^{re}\oplus
\cE_-^{re}$. 
Set
\begin{align}
\mG:=\mF+\mF^{\mathrm{op}}=\mF_++\mF_-+\mF_-^++\mF_+^-,\quad
\mG^{re}:=\mF_++\mF_-^{re}+\mF_-^++\mF_+^{re}.
\end{align}
Let $\mathscr{G}$ (resp. $\mathscr{G}^{re}$) denote the infinite dimensional
bundle associated with $\mG$ (resp. $\mG^{re}$) over $B$. Let 
\begin{align}
P: \mathscr{G}^{re}\to (\cE_+\oplus \cE_-)\oplus (\cE_+^{re}+\cE_-^{re})
\end{align}
denote the obvious block permutation. And denote by 
$\mV:=(\mathrm{Id}_{\cE}\oplus \mU^{-1})\circ P: \mathscr{G}^{re}\to 
\cE\oplus \cE$. Then on $\cE\oplus \cE$, by (\ref{eq:2.30a}),
we have
\begin{align}
\mD:=\left( \begin{array}{cc}
	D(\mF) & 0 \\
	0 & -D(\mF)
\end{array}\right)=\mV D(\mG^{re}) \mV^{-1}.
\end{align}
Denote by $J:\cE\oplus \cE\to \cE\oplus \cE$, $(a,b)\mapsto (b,a)$
the exchange map. Then we have $J^*=J$, $J^2=\mathrm{Id}$ and
$J\mD=-\mD J$. Let $f\in \cC^{\infty}(\R)$ be the even function
introduced in the even case. Then $\mD+Jf(\mD)$ is invertible.
Thus $\mA^{re}:=\mV^{-1}Jf(\mD)\mV$ is a perturbation operator
associated with $\mG^{re}$. By the discussion in (\ref{eq:1.36b}),
$\mA:=\cU^{-1}\mA^{re}\cU$ is a perturbation operator
associated with $\mG$, where $\cU=\mathrm{Id}_{\cE_+}\oplus \mU_-
\oplus \mathrm{Id}_{\cE_-}\oplus \mU_+:\mathscr{G}\to \mathscr{G}^{re}$. 
From the construction of the perturbation $\mA$, 
we see that $\mA^{\mathrm{op}}=\mA$.

If the fibers are zero-dimensional, an equivariant geometric 
family is simply a finite-dimensional $\Z_2$-graded geometric triple
$\underline{E}=\underline{E_+}\oplus \underline{E_-}$.
In this case, $\underline{E^{\mathrm{op}}}
=\underline{E_-^{\mathrm{op}}}\oplus \underline{E_+^{\mathrm{op}}}$. 
Thus for $\mF+\mF^{\mathrm{op}}$, the positive part 
is $H_+=\underline{E_+}\oplus \underline{E_-}$, the negative part
is $H_-=\underline{E_-}\oplus \underline{E_+}$.
Let $C:H_+\to H_-$, $(a,b)\mapsto (b,a)$ be the exchange map.
Then on $H_+\oplus H_-$,
\begin{align}\label{eq:2.34}
\mA:=\left(
\begin{array}{cc}
	0 & C \\
	C & 0
\end{array}\right)
\end{align}
anticommutes with the $\Z_2$-grading $\tau$ of $H_+\oplus H_-$, which 
is $+\mathrm{Id}$ on $H_+$ and $-\mathrm{Id}$ on $H_-$.
Moreover, we have $\mA^*=\mA$, $\mA^2=\mathrm{Id}$.
Since $D(\mathcal{F}+\mF^{\mathrm{op}})=0$, 
we see that $\mA$ is a perturbation operator associated with
$\mF+\mF^{\mathrm{op}}$. From (\ref{eq:2.34}), it is obvious that
$\mA^{\mathrm{op}}=\mA$.

The proof of Lemma \ref{lemma:2,7} is complete.
	\end{proof}

	Note that the ``paired" relation is usually not transitive.
	It is extended to an equivalence relation ``$\sim$'' by
	\begin{multline}\label{eq:2.21}
(\mF,\phi)\sim(\mF',\phi') \iff \text{there exist cycles } 
\{(\mF_i,\phi_i)\}_{0\leq i\leq r}\	\text{such that for any }\  \\ 
				1\leq i\leq r,
	(\mF_i,\phi_i)\ \text{is paired with }(\mF_{i-1},\phi_{i-1}),
(\mF_0,\phi_0)=(\mF,\phi)\ \text{and }(\mF_r,\phi_r)=(\mF',\phi').
	\end{multline}
	The following definition extends \cite[Def. 4.5]{Liu21},
	\cite[Def. 2.14]{Bunke2009} and 
	\cite[Def. 2.19]{BunkeSchick13} to our $g$-version.
	
	\begin{defn}\label{d051}
		Let $\sim$ denote the equivalence relation generated 
		by the relation ``paired" in (\ref{eq:2.21}).  
		The equivariant differential K-group 
		$\widehat{K}_g^0(B)$ (resp. $\widehat{K}_g^{1}(B)$) is the group 
		completion of the abelian semigroup $\widehat{\mathrm{IC}}_g^{
			0}(B)/\sim$ (resp. $\widehat{\mathrm{IC}}_g^{
			1}(B)/\sim$).
	\end{defn}
	
	The non-equivariant version of Definition \ref{d051} is originally 
	explained in \cite{Bunke2009}, which is generalized to the 
	orbifold case in \cite{BunkeSchick13} using the language of stack.
	In \cite{Liu21}, the first author introduces 
	a modified version of the model in \cite{Bunke2009}
	using the eta form in Definition \ref{defn:1.02} and the 
	Dai-Zhang higher spectral flow \cite{DZ98}
	and extended the model to the equivariant case when 
	the compact Lie group action has finite stabilizers only.
	That equivariant model can also be regarded as an 
	analytic model of differential K-theory for orbifolds. 
	
	For $(\mF, \phi)\in \widehat{\mathrm{IC}}_g^{0/1}(B)$, we denote by 
	$[\mF, \phi]\in \widehat{K}_g^{0/1}(B)$ the corresponding class in 
	$g$-equivariant differential K-group. 
	As explained in \cite[(4.18)]{Liu21} (cf. also 
	\cite[Lemmas 2.15, 2.16]{Bunke2009} and  
	\cite[Lemmas 2.20, 2.21]{BunkeSchick13}),
	every element of $\widehat{K}_g^{\bullet}(B)$ can be represented in the 
	form $[\mF, \phi]$. Furthermore, we have	
	$-[\mF, \phi]=[\mF^{\mathrm{op}}, -\phi]$. 
	The following lemma is the $g$-equivariant version of
	\cite[Lemma 2.17]{Bunke2009} and 
	\cite[Lemma 2.22]{BunkeSchick13}.
	
	\begin{lemma}\label{lemma:2.09}
	If $[\mathcal{F}, \phi]=
	[\mathcal{F}', \phi']$, then there exists a family cycle
	$(\mathcal{F}'', \phi'')$ such that $(\mathcal{F}, \phi)
	+(\mathcal{F}'', \phi'')$ is paired with $(\mathcal{F}', \phi')
	+(\mathcal{F}'', \phi'')$.
	\end{lemma}
	\begin{proof}
If $[\mathcal{F}, \phi]=[\mathcal{F}', \phi']$, 
there exists $(\widetilde{\mF}, \tilde{\phi})$ such that
 $(\mathcal{F}, \phi)+(\widetilde{\mF}, \tilde{\phi})\sim 
 (\mathcal{F}', \phi')+(\widetilde{\mF}, \tilde{\phi})$. 
 By (\ref{eq:2.21}), there exist cycles $\{(\mF_i,\phi_i)\}_{0\leq i\leq r}$
such that for any
$1\leq i\leq r$,
$(\mF_i,\phi_i)$ is paired with $(\mF_{i-1},\phi_{i-1})$,
$(\mF_0,\phi_0)=(\mF+\widetilde{\mF},\phi+\tilde{\phi})$ and $(\mF_r,\phi_r)
=(\mF'+\widetilde{\mF},\phi'+\tilde{\phi})$.
From Definition \ref{defn:2.03}, for $1\leq i\leq r$,
there exists a perturbation operator $\mA_i$ associated with 
$\mF_i+\mF_{i-1}^{\mathrm{op}}$ such that
\begin{align}
\phi_i-\phi_{i-1}=\widetilde{\eta}_g((\mF_i+\mF_{i-1}^{\mathrm{op}})|_{B^g},
\mA_i|_{B^g}).
\end{align}
Thus by Definition \ref{defn:1.02}, we have
\begin{align}\label{eq:2.37b}
\phi'-\phi=\sum_{i=1}^r(\phi_i-\phi_{i-1})
=\widetilde{\eta}_g\big((\mF_r+\mF_{r-1}^{\mathrm{op}}+\cdots 
+ \mF_1+\mF_{0}^{\mathrm{op}})|_{B^g},
\mA_r\oplus\cdots\oplus \mA_1|_{B^g}\big).
\end{align}
From Definition \ref{defn:2.03} and (\ref{eq:2.37b}), we see that
$(\mF_r+\sum_{i=1}^{r-1}\mF_i, \phi') $ is paired with 
$(\mF_0+\sum_{i=1}^{r-1}\mF_i, \phi) $.
So we can take $\mF''=\widetilde{F}+\sum_{i=1}^{r-1}\mF_i$, $\phi''=0$. 

The proof of Lemma \ref{lemma:2.09} is complete.
	\end{proof}	
	
	For equivariant geometric families $\mF$ and $\mF'$, we form 
	their product $\mF \times_B \mF'$ over $B$ as follows.
	The total space of the underlying fibration of $\mF \times_B \mF'$ 
	is $ W\times_B W':=\{(w,w')\in W\times W' : \pi(w)=\pi'(w') \}$, 
	which is a fiber bundle over $B$ with fiber $Z\times Z'$. 
	Let $\mathrm{pr}_W:W\times_B W'\rightarrow W$ and 
	$\mathrm{pr}_{W'}:W\times_B W'\rightarrow W'$
	be the obvious projections.
	Let $\mathrm{pr}_W^*E \widehat{\otimes} \mathrm{pr}_{W'}^*E'$
	be the equivariant complex $\Z_2$-graded vector bundle over
	$ W\times_B W'$. 
	The remaining structures of $\mF\times_B \mF'$ are induced 
	in the obvious way. 
	
	For $[\mF,\phi]\in \widehat{K}_g^i(B)$ and $[\mF',\phi']\in 
	\widehat{K}_g^{\bullet}(B)$, $i=0,1$, as in \cite[Def. 4.1]{Bunke2009}, 
	\cite[Def. 3.22]{BunkeSchick13} and \cite[(4.41)]{Liu21},
	we define 
	\begin{multline}\label{eq:1.12}
		[\mF,\phi]\cup[\mF',\phi']:=\left[\mF\times_B\mF', 
		(-1)^i\int_{Z^g}\td_g(\nabla^{TZ},\nabla^{L_Z})\,
		\ch_g(\underline{E})\wedge \phi'\right.
		\\
		\left.+\phi\wedge \int_{Z^{'g}}\td_g(\nabla^{TZ'},\nabla^{L_{Z'}})\,
		\ch_g(\underline{E'})-(-1)^id\phi\wedge \phi'\right].
	\end{multline}
	By \cite[Prop. 4.11]{Liu21} (cf. also 
	\cite[Prop. 4.2]{Bunke2009} and 
	\cite[Prop. 3.23]{BunkeSchick13}), 
	this product is well-defined and commutative.
	So $\widehat{K}_g^{\bullet}(B)$ is a unital graded commutative ring
	with unit $[\underline{\C}, 0]$. Here we consider $\underline{\C}$ 
	as an equivariant geometric family in which
	the fiber is a point. 
	Furthermore, the product (\ref{eq:1.12}) is associative
	(cf. \cite[Prop. 4.11]{Liu21} and \cite[\S 4.1.4]{Bunke2009}).
	
	Let $f: B_1\rightarrow B_2$ be a $G$-equivariant smooth map.
	We define the pullback map
	\begin{align}\label{eq:1.13}
		\hat{f}^*:\widehat{K}_{g}^{\bullet}(B_2)\rightarrow \widehat{K}_{g}^{\bullet}(B_1)
	\end{align}
	as in \cite[(4.39)]{Liu21} (cf. also \cite[\S 2.3.2]{Bunke2009}
	and \cite[\S 2.4.2]{BunkeSchick13}). The product in (\ref{eq:1.12})
	is natural with respect to this pullback map.

	\begin{thm}\label{d031}(Compare with \cite[Th. 4.12]{Liu21}, 
		\cite[\S 3, \S 4]{Bunke2009}
		and \cite[\S 2.5]{BunkeSchick13}) Fix $g\in G$.
		The $g$-equivariant differential K-theory $\widehat{K}_{g}$ is a 
		contravariant functor $B\rightarrow \widehat{K}_{g}^{\bullet}(B)$
		from the category of compact smooth $G$-manifolds
		to unital $\Z_2$-graded commutative rings together with
		well-defined transformations
		
		\begin{enumerate}
			\item $R:\widehat{K}_{g}^{0/1}(B)\rightarrow 
			Z^{\mathrm{even/odd}}(B^g, \R)$ (curvature);
			\item $I: \widehat{K}^{0/1}_g(B)\rightarrow K^{0/1}_G(B)$ 
			(underlying $\mathrm{K}_G$-group);
			\item $a: \Omega^{\mathrm{odd/even}}(B^g, \R)/\Im 
			\,d\rightarrow \widehat{K}_g^{0/1}(B)$ (action of forms),
		\end{enumerate}	
		where  $Z^{\mathrm{even/odd}}(B^g, \R)$ denotes the vector space of
		even/odd degree closed
		differential forms, such that	
		\begin{enumerate}
		\item the following diagram commutes		
		\begin{center}\label{d032}
		\begin{tikzpicture}[>=angle 90]
				\matrix(a)[matrix of math nodes,
				row sep=2em, column sep=2.5em,
				text height=1.5ex, text depth=0.25ex]
				{\widehat{K}^{0/1}_g(B) & K^{0/1}_G(B)\\
				Z^{\mathrm{even/odd}}(B^g, \R)& H_{}^{\mathrm{even/odd}}(
					B^g, \R);\\};
			\path[->](a-1-1) edge node[above]{\footnotesize{$I$}} (a-1-2);
		\path[->](a-1-2) edge node[right]{\footnotesize{$\ch_{g}$}} (a-2-2);
			\path[->](a-1-1) edge node[right]{\footnotesize{$R$}} (a-2-1);
		\path[->](a-2-1) edge node[above]{\tiny{de Rham}} (a-2-2);
		\end{tikzpicture}
		\end{center}
			
			\item the exterior differential operator
			\begin{align}\label{d033}
				d=R\circ a;
			\end{align}
			
			\item  for $x,y\in \widehat{K}_g^*(B)$ and $\alpha\in 
			\Omega_{}^*(B^g,\R)/\Im\, d$, we have
			\begin{align}\label{d203} 
				R(x\cup y)=R(x)\wedge R(y),\quad I(x\cup y)=I(x)\cdot I(y), 
				\quad a(\alpha)\cup x=a(\alpha\wedge R(x));
			\end{align}
			
		\item the following sequence is exact:
		\begin{align}\label{d034}
			K_G^{1/0}(B)\overset{\footnotesize{\ch_{g}}}{\longrightarrow} 
			\Omega^{\mathrm{odd/even}}(B^g, \C)/\Im\, d
			\overset{\footnotesize{a}}{\longrightarrow} \widehat{K}_g^{0/1}(B) 
			\overset{\footnotesize{I}}{\longrightarrow} K_G^{0/1}(B)
			\longrightarrow 0.
		\end{align}	
		\end{enumerate}	
	\end{thm}
	\begin{proof}
		We define the natural transformations $R$, $I$ and $a$
		by
		\begin{align}\label{d058}
			R([\mF, \phi]):=\int_{Z^g}\mathrm{Td}_g(\nabla^{TZ}, 
			\nabla^{L})\mathrm{ch}_g(\underline{E})-d\phi,
		\end{align}
		\begin{align}\label{d054}
			I([\mF, \phi]):=\ind(D(\mF)),
		\end{align}
		and
		\begin{align}\label{d056}
			a(\phi):=[\emptyset, -\phi],
		\end{align}
		where $\emptyset$ denotes the empty geometric family.
		From (\ref{eq:1.25a}), Definitions \ref{defn:2.03} and \ref{d051},
		$R$, $I$ and $a$ are well-defined. 
		The commutativeness of the diagram in (1) follows directly from 
		the equivariant family index theorem.
		The properties (2) and (3) follow from straight calculations 
		by using the definitions.
		
		At last, the proof of  the exactness of (\ref{d034}) is almost the 
		same as the proof of \cite[(4.45)]{Liu21}.
		For completeness, we recall some details.
		
		Since (\ref{eq:1.03}) is surjective, so is $I$.	
		Next, it is obvious that $I\circ a=0$. 
		We need to prove that $\ker I\subseteq \mathrm{Im} a$.
		For $[\mF, \rho]\in \widehat{K}_g^{\bullet}(B)$,  $I([\mF, \rho])=0$ 
		implies $\ind(D(\mF))=0\in K_G^{\bullet}(B)$. 
		We can take a representative $\mF$ such that at least one component 
of the fiber has the nonzero dimension (cf. \cite[Example 2.5]{Liu21}).
By Proposition \ref{prop:1.4}, there exists a perturbation operator 
		$\mA$ with respect to  $D(\mF)$.
		By Definition \ref{defn:2.03} and (\ref{d056}), we have
		$
		[\mF, \rho]=a(\widetilde{\eta}_g(\mF, \mA)-\rho).
		$
		Finally, we will prove that $\ker a=\mathrm{Im} \ch_g$.
	By Proposition \ref{prop:1.11}, Definition \ref{defn:2.03} and 
	(\ref{d056}), for any $x\in K_G^{\bullet}(B)$, there exists
		equivariant geometric family $\mF$ and perturbations
		$\mA_1$ and $\mA_2$, such that 
		\begin{align}\label{d197}
		a\circ \ch_{g}(x)=(\emptyset, \tilde{\eta}_g(\mF, \mA_1)
		-\tilde{\eta}_g(\mF, \mA_2))=[\mF, 0]-[\mF,0]=0.
		\end{align}
		On the other hand,
		if $a(\rho)=0$, for any equivariant geometric family $\mF$ with
		a perturbation operator $\mA$, by Definition \ref{defn:2.03} 
		and (\ref{d056}), we have
		$[\mF, \tilde{\eta}_g(\mF, \mA)-\rho]=a(\rho)=0
		=[\mF, \tilde{\eta}_g(\mF, \mA)].$
		By Lemma \ref{lemma:2.09}, there exists another family cycle 
		$(\mF',\rho')$, 
		such that $(\mF+\mF', \rho'+\tilde{\eta}_g(\mF, \mA)-\rho)$
		 is paired with $(\mF+\mF', \rho'+\tilde{\eta}_g(\mF, \mA))$.
		By Definition \ref{defn:2.03}, there exists a
		perturbation operator $\mA'$ associated 
		$(\mF+\mF')+(\mF+\mF')^{\mathrm{op}}$ such that
		$\rho=\tilde{\eta}_g\left((\mF+\mF')+(\mF+\mF')^{\mathrm{op}},
		\mA' \right)$. From (\ref{eq:2.23}), there exists
		perturbation operator $\mA''$ associated 
		$(\mF+\mF')+(\mF+\mF')^{\mathrm{op}}$ such that
		$\tilde{\eta}_g\left((\mF+\mF')+(\mF+\mF')^{\mathrm{op}},
		\mA'' \right)=0$. So by Proposition \ref{prop:2.01}, 
		 there exists $x\in K_G^{\bullet}(B)$, such that
	\begin{align}\label{eq:2.45}
		\rho=\tilde{\eta}_g\left((\mF+\mF')+(\mF+\mF')^{\mathrm{op}},
		\mA' \right)
	-\tilde{\eta}_g\left((\mF+\mF')+(\mF+\mF')^{\mathrm{op}},\mA'' \right)
		=\ch_g(x).
	\end{align}
		
		The proof of Theorem \ref{d031} is complete.	
	\end{proof}

	In the following, we prove that the above two models 
	of the $g$-equivariant differential K$^0$-ring are isomophic as rings.
	
	For an equivariant geometric triple $\underline{E}$ over $B$,
	we can consider it as an equivariant geometric family over $B$ 
	in which the fiber is a point.
	
	\begin{thm}\label{lemma:1.03}
		The map
		\begin{align}\label{eq:1.09}
			J:\widehat{K}_{FL,g}^0(B)\to \widehat{K}_{g}^0(B),
			\quad [\underline{E},\phi]\mapsto [\underline{E},\phi],
		\end{align}
		is well-defined. Moreover, it is an isomorphism of rings.
	\end{thm}
	\begin{proof}
		The well-definedness of $J$ follows directly from Corollary 
		\ref{cor:2.09}.
		
		By (\ref{eq:2.084}) and (\ref{eq:1.12}), $J$ is a ring 
		homomorphism. From (\ref{eq:1.11}) and (\ref{d034}),  
		we have the commutative diagram
		\begin{equation}\label{eq:3.052} 
			\begin{split}
				\xymatrix{
					K_{G}^{1}(B) \ar[d]^{\Id} \ar[r]^-{\ch_g}
					&(\Omega^{\mathrm{odd}}(B^{g}, \C)/\Im\, d
					\ar[d]^{\Id} \ar[r]^-{a}
					&\widehat{K}_{FL,g}^0(B)\ar[d]^{J} 
					\ar[r]^-{} &K_{G}^0(B)\ar[d]^{\Id} \ar[r] &0\\
					K_{G}^{1}(B) \ar[r]^-{\ch_g}
					&\Omega^{\mathrm{odd}}(B^g, \C)/\Im\, d
					\ar[r]^-{a} 
					&\widehat{K}_g^0(B)\ar[r]^-{I} 
					&K_{G}^0(B)\ar[r] &0.
				}
			\end{split}
		\end{equation}
		Thus by the five lemma, we see that $J$ is an isomorpism.
		
		The proof of Theorem \ref{lemma:1.03} is complete.
	\end{proof}  
	
	For a finite dimensional complex $G$-representation $M$,
	we can also consider $\underline{M}$ in (\ref{eq:2.43}) 
	as an equivariant geometric
	family over $B$ in which the fiber is a point. Then 
	\begin{align}\label{eq:2.37a}
		(\mF, \phi)\mapsto (\underline{M}\times_B\mF,   
		\chi_{M}(g)\cdot\phi)
	\end{align}
	makes $\widehat{K}_{g}^*(B)$ an 
	$R(G)$-module. From (\ref{eq:2.09}) and (\ref{eq:2.37a}),
	the map $J$ in (\ref{eq:1.09}) is also an isomorphism
	of $R(G)$-modules.
	
	
	\subsection{Pushforward: submersion}\label{s0203}

	Let $\pi_Z: W\rightarrow B$ be an equivariant  submersion of 
	closed $G$-manifolds with closed orientable fiber $Z$. 
	We assume that $TZ$ has an equivariant spin$^c$ structure 
	with associated equivariant line bundle $L_Z$. 
	Set $\underline{\pi_Z}=(W, o_Z, T_{\pi_Z}^HW, g^{TZ}, 
	\underline{L_Z})$ as in (\ref{eq:0.02}).
	
	Let
	\begin{align}\label{d155}
		\mF_X=(V, o_X, T_{\pi_X}^HV, g^{TX}, \underline{L_X},\underline{E})
	\end{align}
	be a $G$-equivariant geometric family over $W$
	with equivariant submersion $\pi_X:V\rightarrow W$
	with closed orientable fiber $X$.
	Then $\pi_Y:=\pi_Z\circ \pi_X: V\rightarrow B$ is a smooth
	equivariant
	submersion with closed orientable fiber $Y$. 
	We have the diagram of smooth submersions:
	
	\begin{center}\label{e02001}
	\begin{tikzpicture}[>=angle 90]
			\matrix(a)[matrix of math nodes,
			row sep=2em, column sep=2.5em,
			text height=1.5ex, text depth=0.25ex]
			{X&Y&V\\
				&Z&W&B.\\};
		\path[->](a-1-1) edge (a-1-2);
		\path[->](a-1-2) edge node[left]{\footnotesize{}} (a-2-2);
		\path[->](a-1-2) edge (a-1-3);
		\path[->](a-2-2) edge (a-2-3);
		\path[->](a-1-3) edge node[left]{\footnotesize{$\pi_X$}} (a-2-3);
		\path[->](a-2-3) edge node[above]{\footnotesize{$\pi_Z$}} (a-2-4);
		\path[->](a-1-3) edge node[above]{\footnotesize{$\pi_Y$}} (a-2-4);
	\end{tikzpicture}
	\end{center}

	Set $T_{\pi_X}^HY:=T_{\pi_X}^HV\cap TY$. Then we have the splitting 
	of smooth vector bundles over $V$,
	\begin{align}\label{e02003}
		TY=T_{\pi_X}^HY\oplus TX,
	\end{align}
	where
	\begin{align}\label{e02002}
		T_{\pi_X}^HY\cong \pi_X^*TZ.
	\end{align}
	Let $\nabla^{TX}$, $\nabla^{TY}$ and $\nabla^{TZ}$ be the
	connections defined in (\ref{eq:1.15}). We set
	$\nabla^{TZ,TX}:=\pi_Z^*\nabla^{TZ}\oplus \nabla^{TX}$.
	Let $o_Y:=\pi_X^*o_Z\cup o_X\in H^0(V,\Z_2)$.
	Since $TZ$ and $TX$ both have equivariant spin$^c$ structures, so 
	is $TY$. The equivariant complex line bundle associated with this 
	spin$^c$ structure is
	\begin{align}\label{d156}
		L_Y:=\pi_X^*L_Z\otimes L_X.
	\end{align}
	Take geometric data 
	$(T_{\pi_Y}^HV, g^{TY}, h^{L_Y}, \nabla^{L_Y})$ of $\pi_Y$ such 
	that $T_{\pi_Y}^HV\subset T_{\pi_X}^HV$, $g^{TY}=\pi_X^*g^{TZ}\oplus 
	g^{TX}$, $h^{L_Y}=\pi_X^*h^{L_Z}\otimes h^{L_X}$ and $\nabla^{L_Y}=
	\pi_X^*\nabla^{L_Z}\otimes 1+1\otimes \nabla^{L_X}$. We obtain a new 
	equivariant geometric family over $B$,
	\begin{align}\label{d073}
		\mF_{Y}:=(V, o_Y, T_{\pi_Y}^HV, g^{TY}, \underline{L_Y},
		\underline{E}).
	\end{align}
	The following definition is the $g$-version of \cite[(4.24)]{Liu21} 
	(cf. also \cite[\S 3.2.3]{Bunke2009}, \cite[\S 3.2.3]{BunkeSchick13}).
	
	\begin{defn}
	The pushforward map
		\begin{align}
		\widehat{\pi_Z}_{!}:\widehat{K}_g^{\bullet}(W)\rightarrow \widehat{K}_g^{\bullet}(B)
		\end{align}
		associated with $\pi_Z$
		is defined by 
		\begin{multline}\label{d082}
			\widehat{\pi_Z}_{!}([\mF_X, \phi])=\left[\mF_Y, \int_{Z^g} 
			\td_g(\nabla^{TZ}, \nabla^{L_Z})\wedge \phi\right.
			\\
			+\left.\widetilde{\td}_g\left(\nabla^{TZ,TX}, \nabla^{L_Y}, 
			\nabla^{TY}, \nabla^{L_Y}\right)\ch_g(\underline{E})
			\right]\in \widehat{K}_g^{\bullet}(B).
		\end{multline}
	\end{defn}

	The well-definedness of the pushforward map $\widehat{\pi_Z}_{!}$
	follows directly from
	\cite[Th. 4.7]{Liu21}, where we use the functoriality
	of equivariant $\eta$-forms.
	
	Remark that in \cite{Bunke2009,BunkeSchick13,Liu21}, 
	there is a differential K-orientation term in the definition of
	the pushforward map. Here for simplicity, we omit the 
	differential K-orientation part in (\ref{d082}), which is 
	enough for our use.
	The definition in (\ref{d082}) corresponds to \cite[(4.24)]{Liu21}
	for $\sigma_Z=0$.
	
	If the fiber $X$ is a point,
	from (\ref{d082}), we have
	\begin{align}\label{eq:2.41}
		\widehat{\pi_Z}_{!}([\underline{E}, \phi])=\left[(\underline{\pi_Z},
		\underline{E}), \int_{Z^g} \td_g(\nabla^{TZ}, 
		\nabla^{L_Z})\wedge \phi\right]\in \widehat{K}_g^{\bullet}(B).
	\end{align}

	For $G=\mathbb{T}^m$,  $g\in \mathbb{T}^m$ defines a prime ideal 
	$I(g)$ in $R(\mathbb{T}^m)$, 
	namely all 
	characters of $\mathbb{T}^m$ which vanish at $g$. 
	For any $R(\mathbb{T}^m)$-module 
	$\mathcal{M}$, 
	we denote by $\mathcal{M}_{I(g)}$ the module obtained from 
	$\mathcal{M}$ by localizing at this prime ideal. 
	In fact, 
	we can write an element of 
	$\mathcal{M}_{I(g)}$ by 
	\begin{align}\label{eq:2.086} 
		u/\chi\quad  \text{with}\ u\in \mM,\ \chi\in R(\mathbb{T}^m),\ 
		\chi(g)\neq 0.
	\end{align}
	Two fractions $u/\chi$ and $u'/\chi'$ represent the same element
	of $\mathcal{M}_{I(g)}$ if there exists $\chi''\in R(\mathbb{T}^m)$ 
	with $\chi''(g)\neq 0$ such that
	\begin{align}\label{eq:2.39}
		\chi''\chi'\,u=\chi''\chi\, u'\in \mM.
	\end{align}
	Moreover, $\mathcal{M}_{I(g)}$ is a module over 
	the local ring $R(\mathbb{T}^m)_{I(g)}$.
	
	For a map $\mathcal{M}\to \mathcal{N}$ of $R(\mathbb{T}^m)$-modules,
	we can naturally define the localized map
	$\mathcal{M}_{I(g)}\to \mathcal{N}_{I(g)}$ of 
	$R(\mathbb{T}^m)_{I(g)}$-modules. For example,
	it is easy to see that 
	the pushforward map $\widehat{\pi_Z}_!$ in (\ref{d082})
	is a map of $R(\mathbb{T}^m)$-modules.
	So when $G=\mathbb{T}^m$, for $g\in \mathbb{T}^m$,
	we have a well-defined localized pushforward map
	\begin{align}\label{eq:2.37}
		\widehat{\pi_Z}_{!}:\widehat{K}_g^{\bullet}(W)_{I(g)}\rightarrow 
		\widehat{K}_g^{\bullet}(B)_{I(g)}.
	\end{align}
	For the embedding $\iota:W^{\mathbb{T}}\to W$, the restriction
	on $W^{\mathbb{T}}$ also induces a natural map
	$\hat{\iota}^*: \widehat{K}_g^{\bullet}(W)\to \widehat{K}_g^{\bullet}(W^{\mathbb{T}})$
	of $R(\mathbb{T}^m)$-modules. Then its localized version 
	\begin{align}\label{eq:2.41a}
		\hat{\iota}^*: \widehat{K}_g^{\bullet}(W)_{I(g)}\to 
		\widehat{K}_g^{\bullet}(W^{\mathbb{T}})_{I(g)}
	\end{align}
	is also well-defined as a map of $R(\mathbb{T}^m)_{I(g)}$-modules.
	By Proposition \ref{prop:1.04} and Theorem \ref{lemma:1.03}, for
	$g\notin A_{\mathbb{T}}$, we can define the pushforward map 
	$\hat{\iota}_!:\widehat{K}^0_g(W^{\mathbb{T}})\rightarrow 
	\widehat{K}^0_g(W)$ through $J$ in Theorem \ref{lemma:1.03} as a 
	well-defined $R(G)$-module homomorphism,
	which induces the well-defined
	localized pushforward map
	\begin{align}\label{eq:1.30}
		\hat{\iota}_!:\widehat{K}^0_g(W^{\mathbb{T}})_{I(g)}\rightarrow 
		\widehat{K}^0_g(W)_{I(g)}.
	\end{align}
	
	On the other hand, since localization preserves exact sequences,
	from (\ref{d034}), we have the exact sequence
	\begin{align}\label{eq:2.46}
		K_G^{1/0}(B)_{I(g)}\overset{\footnotesize{\ch_{g}}}{\rightarrow} 
		\left(\Omega^{\mathrm{odd/even}}(B^g, \C)/\Im\, d\right)_{I(g)}
		\overset{\footnotesize{a}}{\rightarrow} \widehat{K}_g^{0/1}(B)_{I(g)} 
		\overset{\footnotesize{I}}{\rightarrow} K_G^{0/1}(B)_{I(g)}
		\rightarrow 0.
	\end{align}	
	If $G=S^1$ and $B=\mathrm{pt}$ is a point, we have 
	$K_G^{0}(\mathrm{pt})_{I(g)}=R(S^1)_{I(g)} $, 
	$\left(\Omega^{\mathrm{even}}(\mathrm{pt}^g, \C)/\Im\, d\right)_{I(g)}
	=\C$ and $K_G^{1}(\mathrm{pt})_{I(g)}=0$. So from (\ref{eq:2.46}), 
	we have
	\begin{align}\label{eq:2.63}
		\widehat{K}_g^{1}(\mathrm{pt})_{I(g)}=\C/\Q_g.
	\end{align}
	Here $\Q_g$ is the image of $\ch_g$ on $R(S^1)_{I(g)}$
	which is originally defined in \cite[(0.13)]{LM20}.

	
	\section{Localization}\label{s03}

	In this section, we assume $G=\mathbb{T}^m$ and establish our main 
	result: a localization formula for $\eta$-forms.
	In Section \ref{s0301}, we develop a localization formula
	in the geometric family version of the $g$-equivariant differential 
	K-theory. In Section \ref{s0302}, we establish a localization formula 
	for $\eta$-forms. In Section \ref{s0303}, we explain how to obtain 
	Corollaries \ref{thm:0.03} and \ref{cor:0.5} from Theorem \ref{thm:1.9}.  In Section \ref{s0304},
	we prove Theorem \ref{thm:0.7}.


	\subsection{Localization theorem for $g$-equivariant differential 
		K-theory}\label{s0301}
	
	Let $W$ be a closed manifold with smooth $\mathbb{T}^m$-action.
	The following theorem is the differential K-theory version of 
	\cite[Th. 1.1]{ASegal68}, and is the generalization of
	\cite[Th. 3.2]{LM20}. We use the notations in Section \ref{s0101}.
	
	\begin{thm}[Localization Theorem]\label{thm:1.06} 
		For the embedding $\iota: W^{\mathbb{T}}\to W$ and $g\in 
		\mathbb{T}^m\backslash A_{\mathbb{T}}$, 
		the restriction map
		$\hat{\iota}^*:\widehat{K}_g^{\bullet}(W)_{I(g)}\rightarrow 
		\widehat{K}_g^{\bullet}(W^{\mathbb{T}})_{I(g)}$ in (\ref{eq:2.41a})
		is an $R(\mathbb{T}^m)_{I(g)}$-module isomorphism.
	\end{thm} 
	\begin{proof}
		If $g\in \mathbb{T}^m\backslash A_{\mathbb{T}}$, $W^g=W^{
			\mathbb{T}}$.
		From (\ref{eq:2.46}), 
		we have the commutative diagram
		\begin{equation}\label{eq:1.29} 
		\begin{gathered}
			\xymatrix@C=0.9em{
			K_{\mathbb{T}^m}^{1/0}(W)_{I(g)} \ar[d]^{\iota^*} \ar[r]^-{\ch_g}
			&\left(\Omega^{\mathrm{odd/even}}(W^{\mathbb{T}}, \C)/\Im\, 
			d\right)_{I(g)} \ar[d]^{\Id} \ar[r]^-{a}
			&\widehat{K}_g^{0/1}(W)_{I(g)}\ar[d]^{\hat{\iota}^*} 
			\ar[r]^-{I} &K_{\mathbb{T}^m}^{0/1}(W)_{I(g)}\ar[d]^{\iota^*} 
			\ar[r] &0     \\
			K_{\mathbb{T}^m}^{1/0}(W^{\mathbb{T}})_{I(g)}  \ar[r]^-{\ch_g}
			&\left(\Omega^{\mathrm{odd/even}}(W^{\mathbb{T}}, \C)/\Im\, 
			d\right)_{I(g)} \ar[r]^-{a} 
			&\widehat{K}_g^{0/1}(W^{\mathbb{T}})_{I(g)}\ar[r]^-{I} 
			&K_{\mathbb{T}^m}^{0/1}(W^{\mathbb{T}})_{I(g)}\ar[r] &0.
			}
		\end{gathered}		
		\end{equation}
		Here $\iota^*:K_{\mathbb{T}^m}^{0/1}(W)_{I(g)}
		\rightarrow K_{\mathbb{T}^m}^{0/1}(W^{T})_{I(g)}$
		is the localized restriction map of 
		$R(\mathbb{T}^m)_{I(g)}$-modules induced by $\iota$.
		Using the localization in topological K-theory 
		\cite[Th. 1.1]{ASegal68}, $\iota^*$ is 
		an $R(\mathbb{T}^m)_{I(g)}$-module isomorphism.
		Then by the five lemma on (\ref{eq:1.29}), $\hat{\iota}^*$ in 
		(\ref{eq:1.29}) is an 
		isomorphism of $R(\mathbb{T}^m)_{I(g)}$-modules.
		
		The proof of Theorem \ref{thm:1.06} is complete.
	\end{proof}

	For an ungraded equivariant geometric triple $\underline{E}$, we denote by
	\begin{align}\label{eq:2.058} 
		\lambda_t(\underline{E})=\sum_{j\geq 0}\Lambda^j(\underline{E})t^j.
	\end{align}
	By Theorem \ref{lemma:1.03} and (\ref{eq:2.058}),
	$\left[\underline{\lambda_{-1}(N^*)},0\right] \in
	\widehat{K}_g^0(W^{\mathbb{T}})$.
	The following proposition is an extension of 
	\cite[Th. 2.17]{LM20}. Here we first introduce a
	short but not constructive proof and then generalize the
	constructive proof of \cite[Th. 2.17]{LM20} to the 
	torus case.
	
	\begin{prop}\label{prop:2.03} 
		For $g\in \mathbb{T}^m\backslash A_{\mathbb{T}}$,
		$\left[\underline{\lambda_{-1}(N^*)}, 0\right]$ is invertible in 
		$\widehat{K}_{g}^0(W^{\mathbb{T}})_{I(g)}$. Therefore there exist 
		equivariant geometric 
		triples $\underline{E_1}$, $\underline{E_2}$ over $W^{\mathbb{T}}$
		and $\chi\in \Z[x_1,\cdots,x_m]$
		such that 
		\begin{align}\label{eq:3.02} 
			\left[\underline{\lambda_{-1}(N^*)},0\right]^{-1}=
			\left(\left[\underline{E_1},0\right]-
			\left[\underline{E_2},0\right]\right)/\chi\in 
			\widehat{K}_{g}^0(W_{}^{\mathbb{T}})_{I(g)}.
		\end{align}
	\end{prop}
	\begin{proof}

		Fix $y\in W^{\mathbb{T}}$. Then $\lambda_{-1}(N^*)_y$ is a 
		representation of $\mathbb{T}^m$ and the character 
		$\chi_{\lambda_{-1}(N^*)_y}(g)\neq 0$ 
		(see e.g., \cite[(1.37)]{LM20}). From the proof of 
		\cite[Lemma 2.7]{ASegal68}, we know that $\lambda_{-1}(N^*)$
		is invertible in $K_{\mathbb{T}^m}^0(W^{\mathbb{T}})_{I(g)}$.
		By the localization of (\ref{eq:3.032}), we can define
		\begin{align}
			\lambda_{-1}(N^*)\,\cup : K_{\mathbb{T}^m}^1(W^{\mathbb{T}
			})_{I(g)}\to
			K_{\mathbb{T}^m}^1(W^{\mathbb{T}})_{I(g)},
		\end{align}
	which is also an $R(\mathbb{T}^m)_{I(g)}$-module isomorphism. 
	Then from (\ref{eq:1.28}), (\ref{eq:2.084}) and (\ref{eq:1.12}), 
	we have a commutative diagram
	\begin{equation}
		\begin{gathered}
			\xymatrix@C=0.9em{
			K_{\mathbb{T}^m}^{1}(W^{\mathbb{T}})_{I(g)} \ar[d]^{\lambda_{-1}
				(N^*)\cup} \ar[r]^-{\ch_g}
	   &\left(\Omega^{\mathrm{odd}}(W^{\mathbb{T}}, \C)/\Im\, d\right)_{I(g)}
		\ar[d]^{\ch_g(\underline{\lambda_{-1}(N^*)})\wedge} \ar[r]^-{a}
	   &\widehat{K}_g^{0}(W^{\mathbb{T}})_{I(g)}\ar[d]^{\left[\underline{
				\lambda_{-1}(N^*)}, 0\right]\cup} 
		\ar[r]^-{I} &K_{\mathbb{T}^m}^{0}(W^{\mathbb{T}})_{I(g)}\ar[d]^{
			\lambda_{-1}(N^*)\cup} \ar[r] &0\\
	K_{\mathbb{T}^m}^{1}(W^{\mathbb{T}})_{I(g)}  \ar[r]^-{\ch_g}
	&\left(\Omega^{\mathrm{odd}}(W^{\mathbb{T}}, \C)/\Im\, d\right)_{I(g)}
		\ar[r]^-{a} 
	&\widehat{K}_g^{0}(W^{\mathbb{T}})_{I(g)}\ar[r]^-{I} 
	&K_{\mathbb{T}^m}^{0}(W^{\mathbb{T}})_{I(g)}\ar[r] &0.
				}
			\end{gathered}
		\end{equation}
		By the five lemma, the map
		\begin{align}
			\left[\underline{
				\lambda_{-1}(N^*)}, 0\right]\,\cup:\widehat{K}_g^0(W^{
				\mathbb{T}})_{I(g)}
			\to \widehat{K}_g^0(W^{\mathbb{T}})_{I(g)}
		\end{align}
		is an $R(\mathbb{T}^m)_{I(g)}$-module isomorphism. 
		So $\left[\underline{\lambda_{-1}(N^*)}, 0\right]^{-1}\in 
		\widehat{K}_g^0(W^{\mathbb{T}})_{I(g)}$.
		
		By Theorem \ref{lemma:1.03} and (\ref{eq:2.086}),
		there exist geometric triples $\underline{E_1'}$, $\underline{E_2'}$ 
		over $W^{\mathbb{T}}$ and 
		$\chi'\in R(\mathbb{T}^m)=\Z[x_1,\cdots, x_m,x_1^{-1},\cdots, 
		x_m^{-1}]$, such that
		$\left[\underline{\lambda_{-1}(N^*)},0\right]^{-1}=
		\left(\left[\underline{E_1'},0\right]-
		\left[\underline{E_2'},0\right]\right)/\chi'$.
		Write $\chi'=\chi/\chi''$ such that
		$\chi, \chi''\in \Z[x_1,\cdots, x_m]$. 
		Let $\underline{M_{\chi''}}$ denote the equivariant
		geometric triple in (\ref{eq:2.43}) associated with the 
		character $\chi''$. For $\underline{E_i}=\underline{M_{\chi''}}
		\otimes \underline{E_i'}$, $i=1,2$, we have
		$\left[\underline{\lambda_{-1}(N^*)},0\right]^{-1}=
		\left(\left[\underline{E_1},0\right]-
		\left[\underline{E_2},0\right]\right)/\chi$.
		
		The proof of Proposition \ref{prop:2.03} is complete.
	\end{proof} 
	
	In fact, as in \cite[Th. 2.17]{LM20},
	we can also construct equivariant geometric 
	triples $\underline{E_1}$, $\underline{E_2}$ over $W^{\mathbb{T}}$
	and $\chi\in \mathbb{Z}[x_1,\cdots,x_m]$ in (\ref{eq:3.02}) 
	explicitly in torus case.
	
	We first recall the $\lambda$-ring structure on 
	$\widehat{K}^0(W^{\mathbb{T}})$ for non-equivariant case,
	where $\widehat{K}^0(W^{\mathbb{T}})$ is the non-equivariant 
	differential K$^0$-group in \cite[Def. 2.6]{FreedLott10}. In our notation, 
	$\widehat{K}^0(W^{\mathbb{T}})$ is just 
	$\widehat{K}^0_{FL,g}(W^{\mathbb{T}})$ 
	for trivial group $G$ and unity $g$.
	For $k\in \N$ and a cycle $(\underline{F}, \phi)$ of 
	$\widehat{K}^0(W^{\mathbb{T}})$, from  
	\cite[Th. 2.6]{LM20}, there exists $\phi_k\in 
	\Omega^{\mathrm{odd}}(W^{\mathbb{T}}, \R)/\mathrm{Im}\,d$,
	such that the map $(\underline{F}, \phi)\mapsto 
	(\Lambda^k(\underline{F}), \phi_k)$ forms a well-defined
	map
	\begin{align}\label{eq:3.07}
	\lambda^k: \widehat{K}^0(W^{\mathbb{T}})\to 
	\widehat{K}^0(W^{\mathbb{T}}).
	\end{align}
	Moreover, $\phi=0$ implies $\phi_k=0$ and the map
	(\ref{eq:3.07}) forms a pre-$\lambda$-ring structure
	on $\widehat{K}^0(W^{\mathbb{T}})$.
	Recently, the authors proved that this pre-$\lambda$-ring
	structure is in fact a $\lambda$-ring structure \cite{LM26}.
	
	In the following, if $\underline{E}$ is an ungraded
	 equivariant
	geometric triple, 
	we will denote by $\underline{'E}$ and $'E$ the 
	corresponding ungraded geometric triple and vector 
	bundle forgetting 
	the group action respectively. 
	
		The $\gamma$-operations are defined by
	\begin{align}\label{eq:2.055} 
		\gamma_t(\underline{'E})=\sum_{j\geq 0}\gamma^j(\underline{'E})t^j
		:=\lambda_{\frac{t}{1-t}}(\underline{'E}).
	\end{align}	
	So 
	$\gamma^i(\underline{'E})$, defined as in \eqref{eq:2.055}, 
	is a finite dimensional virtual Hermitian vector 
	bundle with induced metric and connection. We also denote it by 
	$\underline{\gamma^i('E)}$.
	
	Let $r$ be the rank of $'E$. 
	Let $\underline{\mathrm{rk} 'E}$ be the
	trivial complex vector bundle of rank $r$ with trivial metric
	and connection. 
	From \cite[(2.65), (2.66)]{LM20},
	formally, we can express
	\begin{align}\label{eq:2.064} 
		\lambda_{t}(\underline{'E })^{-1}
		=(1+t)^{-r}\left(1+\sum_{k=1}^{\infty}t^k(1+t)^{-k}
		\big(P_{k,+}(\underline{'E})-P_{k,-}(\underline{'E})\big) \right),
	\end{align} 
	where $P_{k,\pm}(\underline{'E})$ are finite dimensional
	Hermitian vector 
	bundles with induced metrics and connections
	obtained from 
	\begin{align}
		P_{k,+}(\underline{'E})-P_{k,-}(\underline{'E})
		=\sum_{\substack{(n_1,\cdots,n_r)\in \N^r,\\ 
				\sum_{i=1}^ri\cdot n_i=k}}(-1)^{\sum_{i=1}^r n_i}
		\frac{(\sum_{i=1}^r n_i)!}{\prod_{i=1}^rn_i!} 
		\prod_{i=1}^r\left(\gamma^{i}(\underline{'E }
		-\underline{\rank
			'E })\right)^{n_i}.
	\end{align}
	We denote by  $P_{k,+}(\underline{'E})=P_{k,-}(\underline{'E})=0$
	for $k<0$ and $P_{0,+}(\underline{'E})=\underline{\C}$,
	$P_{0,-}(\underline{'E})=0$. 
	By \cite[Corollary 2.13]{LM20}, there exists $\mN>0$ such that
	for any $k>\mN$, 
	\begin{align}\label{eq:2.081}
		[P_{k,+}(\underline{'E}),0]=[P_{k,-}(\underline{'E}),0]
		\in\widehat{K}^0(W^{\mathbb{T}}).
	\end{align}
	
	Recall that we denote by $N_{\alpha}$ the normal bundle of
	a connected component $W_{\alpha}^{\mathbb{T}}$ in $W$, which 
	has a complex structure.
	From (\ref{eq:1.04a}),
	we write
	\begin{align}
		N_{\alpha}=\sum_{j=1}^k\beta_j\cdot \,'N_{\alpha}(\beta_j),
		\quad N_{\alpha}(h)=\sum_{j=1}^k\beta_j(h)\cdot 
		\,'N_{\alpha}(\beta_j),
		\quad \forall h\in \mathbb{T}^m.
	\end{align}
	So for $h\in \mathbb{T}^m$, 
	\begin{multline}
		\lambda_{-1}(N_{\alpha}^{*})(h)=\bigotimes_{j=1}^k\left(
		1+\sum_{l=1}^{\rank N_{\alpha}(\beta_j)}(-1)^l
		\beta_j(h)^{-l}\cdot\Lambda^l\left('N_{\alpha}(\beta_j
		)^{*}\right)\right)
		\\
		=\bigotimes_{j=1}^k\lambda_{-\beta_j(h)^{-1}}
		\left('N_{\alpha}(\beta_j
		)^{*}\right).
	\end{multline}
	
	For equivariant geometric triple, 
	we write 
	\begin{align}\label{eq:2.088}
		\underline{\lambda_{-1}(N_{\alpha}^{*})}(h)=\bigotimes_{j=1}^{k}
		\lambda_{-\beta_j(h)^{-1}}
		\left(\underline{'N_{\alpha}(\beta_j
			)^{*}}\right)
	\end{align}
	Set 
	\begin{align}\label{eq:2.089} 
		r_{\alpha,j}=\rank N_{\alpha}(\beta_j),\quad m_{\alpha}=\dim 
		W_{\alpha}^{\mathbb{T}}.
	\end{align}
	By (\ref{eq:2.064}), formally we have
	\begin{multline}\label{eq:2.090}
		\left(\lambda_{-\beta_j(h)^{-1}}
		\left(\underline{'N_{\alpha}(\beta_j
			)^{*}}\right)\right)^{-1}
		\\
		=(1-\beta_j(h)^{-1})^{-r_{\alpha,j}}
		\left(1+\sum_{l=1}^{\infty}
		\left(\frac{-\beta_j(h)^{-1}}{1-\beta_j(h)^{-1}}\right)^l
		\left(P_{l,+}\left(\underline{'N_{\alpha}(\beta_j
			)^{*}}\right)-P_{l,-}
		\left(\underline{'N_{\alpha}(\beta_j
			)^{*}}\right)\right)\right)
		\\
		=\frac{\beta_j(h)^{r_{\alpha,j}}}{(\beta_j(h)-1)^{r_{\alpha,j}}}
		\left(1+\sum_{l=1}^{\infty}
		\frac{(-1)^l}{(\beta_j(h)-1)^{l}}
		\left(P_{l,+}\left(\underline{'N_{\alpha}(\beta_j
			)^{*}}\right)-P_{l,-}
		\left(\underline{'N_{\alpha}(\beta_j
			)^{*}}\right)\right)\right).
	\end{multline}
	By Theorem \ref{lemma:1.03} and (\ref{eq:2.081}), there exists 
	$\mN_{r_{\alpha,j}, m_{\alpha}}>0$ 
	such that for any $l> 
	\mN_{r_{\alpha,j}, m_{\alpha}}$, 
	\begin{align}\label{eq:2.091}
		\left[P_{l,+}\left(\underline{'N_{\alpha}(\beta_j
			)^{*}}\right)-P_{l,-}
		\left(\underline{'N_{\alpha}(\beta_j
			)^{*}}\right),0\right]=0
		\in \widehat{K}^0(W_{\alpha}^{\mathbb{T}}).
	\end{align}
	
	For $\mN\in \N$,
	we set
	\begin{multline}\label{eq:2.092} 
		\left(\lambda_{-\beta_j(h)^{-1}}
		\left(\underline{'N_{\alpha}(\beta_j
			)^{*}}\right)\right)^{-1}_{\mN}
		\\
		:=\frac{\beta_j(h)^{r_{\alpha,j}}}{(\beta_j(h)-1)^{r_{\alpha,j}}}
		\left(1+\sum_{l=1}^{\mN}
		\frac{(-1)^l}{(\beta_j(h)-1)^{l}}
		\left(P_{l,+}\left(\underline{'N_{\alpha}(\beta_j
			)^{*}}\right)-P_{l,-}
		\left(\underline{'N_{\alpha}(\beta_j
			)^{*}}\right)\right)\right),
	\end{multline}
	and
	\begin{align}\label{eq:3.14}
		\underline{\lambda_{-1}(N_{\alpha}^*)^{-1}_{\mN}}(h)
		:=\bigotimes_{j=1}^k
		\left(\lambda_{-\beta_j(h)^{-1}}
		\left(\underline{'N_{\alpha}(\beta_j
			)^{*}}\right)\right)^{-1}_{\mN}.
	\end{align}
	
	For $g\in \mathbb{T}^m\backslash A_{\mathbb{T}}$,
	$W^g=W^{\mathbb{T}}$, thus $\beta_j(g)-1\neq 0$.
	From (\ref{eq:2.092}) and (\ref{eq:3.14}), we can construct 
	equivariant geometric triples $\underline{E_{\alpha,1}^{\mN}}$
	, $\underline{E_{\alpha,2}^{\mN}}$ and $\chi_{\alpha}^{\mN}\in 
	\mathbb{Z}[x_1,\cdots, x_m]$ explicitly such that
	\begin{align}\label{eq:3.15}
		\left[\underline{\lambda_{-1}(N_{\alpha}^*)^{-1}_{\mN}},0\right]
		=\left(\left[\underline{E_{\alpha,1}^{\mN}},0\right]-
		\left[\underline{E_{\alpha,2}^{\mN}},0\right]\right)/
		\chi_{\alpha}^{\mN}
		\in \widehat{K}_{g}^0(W_{\alpha}^{\mathbb{T}})_{I(g)}.
	\end{align}
	
	Set
	\begin{align}\label{eq:3.16}
		\mathcal{N}_0:=\sup_{\alpha,j}\mN_{r_{\alpha,j}, m_{\alpha}}.
	\end{align}
	From (\ref{eq:2.091})-(\ref{eq:3.14}), we see that for any 
	$\mN, \mN'>\mN_0$,
	\begin{align}\label{eq:2.093} 
		\left[\underline{\lambda_{-1}(N_{\alpha}^*)^{-1}_{\mN}},0\right]
		=\left[\underline{\lambda_{-1}(N_{\alpha}^*)^{-1}_{\mN'}},0\right]
		\in\widehat{K}_{g}^0(W_{\alpha}^{\mathbb{T}})_{I(g)}.
	\end{align}
	Moreover from (\ref{eq:2.084}), (\ref{eq:2.088})-(\ref{eq:2.093}), 
	for any $\mN>\mN_0$, we have
	\begin{align}\label{eq:2.094} 
	\left[\underline{\lambda_{-1}(N_{\alpha}^*)},0\right]\cup
	\left[\underline{\lambda_{-1}(N_{\alpha}^*)^{-1}_{\mN}},0\right]
	=1\in \widehat{K}_{g}^0(W_{\alpha}^{\mathbb{T}})_{I(g)}.
	\end{align}
	Thus 
	\begin{align}\label{eq:2.095} 
		\left[\underline{\lambda_{-1}(N_{\alpha}^*)},0\right]^{-1}=
		\left[\underline{\lambda_{-1}(N_{\alpha}^*)^{-1}_{\mN}},0\right]
		=\left(\left[\underline{E_{\alpha,1}^{\mN}}, 0 \right] 
		-\left[\underline{E_{\alpha,2}^{\mN}}, 0 \right]\right)/
		\chi_{\alpha}^{\mN}\in 
		\widehat{K}_{g}^0(W_{\alpha}^{\mathbb{T}})_{I(g)}
	\end{align}
	for any $\mN>\mN_0$.
	As in \cite[Corollary 2.18]{LM20}, 
	for any $\mN>\mN_0$,
	$g\in \mathbb{T}^m\backslash A_{\mathbb{T}}$,
	\begin{align}\label{eq:2.098} 
		\ch_g\left(\underline{\lambda_{-1}(N_{\alpha}^*)}\right)\cdot
		\ch_g\left(\underline{\lambda_{-1}(N_{\alpha}^*)^{-1}_{\mN}}\right)
		=1\in \Omega^{*}(W_{\alpha}^{\mathbb{T}},\C),
	\end{align}
	where 
	\begin{align}
		\ch_g\left(\underline{\lambda_{-1}(N_{\alpha}^*)^{-1}_{\mN}}\right)
		=\chi_{\alpha}^{\mN}(g)^{-1}\left(\ch_g\left(
		\underline{E_{\alpha,1}^{\mN}}\right)
		-\ch_g\left(\underline{E_{\alpha,2}^{\mN}}\right) \right).
	\end{align}
	Remark that if $G=S^1$, $\underline{\lambda_{-1}(
		N_{\alpha}^*)^{-1}_{\mN}}$
	is the same
	as in \cite[(2.100)]{LM20}.
	
	The proof of the following theorem is the same as that of 
	\cite[Th. 3.4]{LM20}.
	
	\begin{thm}\label{thm:3.12} 
		For $g\in \mathbb{T}^m\backslash A_{\mathbb{T}}$,
		the localized pushforward map
		$\hat{\iota}_!$ in (\ref{eq:1.30})
		is a $R(\mathbb{T}^m)_{I(g)}$-module isomorphism and 
		\begin{align}\label{eq:3.055} 
			\hat{\iota}^*\circ 
			\hat{\iota}_!=\left[\underline{\lambda_{-1}(N^*)},0\right]\cup:
			\widehat{K}_g^0(W^{\mathbb{T}})_{I(g)}
			\overset{{\sim}}{\longrightarrow}
			\widehat{K}_g^0(W^{\mathbb{T}})_{I(g)}.
		\end{align}
		Thus the inverse map of $\hat{\iota}^*:
		\widehat{K}^0_g(W)_{I(g)}\rightarrow 
		\widehat{K}^0_g(W^{\mathbb{T}})_{I(g)}$ is given by
		$\hat{\iota}_!\circ \left[\underline{\lambda_{-1}(N^*)},
		0\right]^{-1}\cup$.
	\end{thm}

	
	\subsection{Localization for geometric families}\label{s0302}
	
	In this subsection, we will prove the main result of this paper.
	We use the notations in Sections \ref{s0102}, \ref{s0104} and
	\ref{s0203}. 
	
	Let $B$ be a closed manifold with trivial $\mathbb{T}^m$-action.
	Let $\mathcal{F}$ be a $G$-equivariant geometric family in 
	(\ref{eq:1.04}) with $G=\mathbb{T}^m$. 
	In order to obtain a nontrivial result, we may assume that
	$W$ and $B$ are connected and
	$\dim Z>0$.
	Let $\underline{\pi_{Z^{\mathbb{T}}}}$ be the family of geometric 
	data in (\ref{eq:1.38b}). Then it induces a localized
	pushforward map $\widehat{\pi_{Z^{\mathbb{T}}}}_!: 
	\widehat{K}^0_g(W^{\mathbb{T}})_{I(g)}
	\to \widehat{K}^{\bullet}_g(B)_{I(g)}$ as in (\ref{eq:2.37})
	for $g\in \mathbb{T}^m$, which is also a map of 
	$R(\mathbb{T}^m)_{I(g)}$-modules.

	In the following, we prove the commutativity of the diagram 
	in (\ref{eq:0.10}).

	\begin{thm}\label{thm:1.08} 
		For any $g\in \mathbb{T}^m\backslash A_{\mathbb{T}}$, 
		the following diagram commutes
		\begin{equation}\label{eq:0.14} 
			\begin{split}
				\xymatrix{
					\widehat{K}^0_g(W^{\mathbb{T}})_{I(g)} 
					\ar[rd]_{\widehat{\pi_{Z^{\mathbb{T}}}}_!}
					&&\widehat{K}^0_g(W)_{I(g)}
					\ar[ll]_{
						\left[\underline{\lambda_{-1}(N^*)},0\right]^{-1}
						\cup\,\hat{\iota}^*} 
					\ar[ld]^{\widehat{\pi_Z}_{!}}\\
					&\widehat{K}^{\bullet}_g(B)_{I(g)}.&
				}
			\end{split}
		\end{equation}	
	\end{thm}
	\begin{proof}
		By Theorem \ref{thm:3.12}, since $g\in \mathbb{T}^m\backslash 
		A_{\mathbb{T}}$, we only need to prove the commutativity of the 
		following diagram,
		\begin{equation}\label{eq:1.35} 
			\begin{split}
				\xymatrix{
					\widehat{K}^0_g(W^{g})_{I(g)} 
					\ar[rd]_{\widehat{\pi_{Z^{\mathbb{T}}}}_!}
					\ar[rr]^{
						\hat{\iota}_{!}} 
					&&\widehat{K}^0_g(W)_{I(g)}	
					\ar[ld]^{\widehat{\pi_Z}_!}\\
					&\widehat{K}^{\bullet}_g(B)_{I(g)}.&
				}
			\end{split}
		\end{equation}
		
		By \cite[(1.38)]{LM20}, (\ref{eq:1.026}) and (\ref{eq:2.41}), 
		for $[\underline{\mu}, \phi]/\chi\in \widehat{K}^0_g(W^{g})_{I(g)}$, 
		we have
	\begin{align}\label{eq:1.36}
		\widehat{\pi_{Z^{\mathbb{T}}}}_!([\underline{\mu}, \phi]/\chi)
		=\left[(
		\underline{\pi_{Z^{\mathbb{T}}}},\underline{\mu}), \int_{Z^g} 
		\td_g(\nabla^{TZ^g}, \nabla^{L_{Z^g}})\wedge \phi
		\right]/\chi,
	\end{align}
	and
	\begin{multline}\label{eq:1.37}
	\widehat{\pi_Z}_{!}\circ \hat{\iota}_{!}([\underline{\mu}, \phi]/\chi)
	=\widehat{\pi_Z}_{!}\left([\underline{\xi}, \ch_g(
	\underline{\Lambda^{}(N^*)})\wedge \phi ]/\chi\right)
		\\
	=\left[(\underline{\pi_Z}, \underline{\xi}), \int_{Z^g} 
	\td_g(\nabla^{TZ}, 
	\nabla^{L_{Z}})\wedge \ch_g(\underline{\Lambda^{}(N^*)})
	\wedge\phi \right]/\chi
		\\
	=\left[(\underline{\pi_Z}, \underline{\xi}), \int_{Z^g} 
	\td_g(\nabla^{TZ^g}, \nabla^{L_{Z^g}})\wedge \phi
	\right]/\chi.
	\end{multline}
		
		We consider the embedding $\iota':W^{g}\sqcup W^{g}
		\to W\sqcup W^{g}$ such that restricted to the first part 
		$\iota'=\iota$ and to the second part $\iota'=\mathrm{Id}$.
		Note that $(\underline{\mu},\underline{\mu^{\mathrm{op}}})$ is 
		an equivariant geometric triple over $W^{g}\sqcup W^{g}$.
		Then we have
		\begin{align}
		\widehat{\iota'}_*(\underline{\mu},\underline{\mu^{\mathrm{op}}})
		=(\underline{\xi},\underline{\mu^{\mathrm{op}}}).
		\end{align}
		
		Recall that we denote by $\mF_{W}:=(\underline{\pi_Z}, \underline{\xi})$
		and $\mF_{W^{\mathbb{T}}}:=(
		\underline{\pi_{Z^{\mathbb{T}}}},\underline{\mu})$
		in (\ref{eq:1.80}).
		Since $\mathrm{Ind}(\mF_{W^{\mathbb{T}}}
		+(\mF_{W^{\mathbb{T}}})^{\mathrm{op}})=0\in
		K_{\mathbb{T}^m}^{\bullet}(B)$, from \cite[Prop. 3.6]{Liu19}, we have
		$\mathrm{Ind}(\mF_{W}+(\mF_{W^{\mathbb{T}}})^{\mathrm{op}})
		=0\in K_{\mathbb{T}^m}^{\bullet}(B)$. By (\ref{eq:2.23}),
		there exists a perturbation operator $\mA'$ associated with 
		$\mF_{W^{\mathbb{T}}}+(\mF_{W^{\mathbb{T}}})^{\mathrm{op}}$ 
		such that 
$\widetilde{\eta}_g(\mF_{W^{\mathbb{T}}}
+(\mF_{W^{\mathbb{T}}})^{\mathrm{op}}, \mA')=0$.
		So from Proposition \ref{prop:1.10}, there exist 
		perturbation operator $\mA$ associated with 
		$\mF_{W}+(\mF_{W^{\mathbb{T}}})^{\mathrm{op}}$ and 
		$x\in K_{\mathbb{T}^m}^{\bullet}(B)$ such that
		\begin{align}\label{eq:3.28}
		\widetilde{\eta}_g(\mF_{W}+(\mF_{W^{\mathbb{T}}})^{\mathrm{op}}, 
		\mA)=\ch_g(x).
		\end{align}
		From Definition \ref{defn:2.03} and (\ref{eq:3.28}), 
		$(\mF_{W^{\mathbb{T}}}, 0)$ is paired with $(\mF_{W},\ch_g(x))$. 
		Note that from Theorem \ref{d031} (4), 
		$[\emptyset, \ch_g(x)]=0\in \widehat{K}^{\bullet}_g(B)$.
		We have $[\mF_{W^{\mathbb{T}}}, 0]
		=[\mF_{W}, 0]\in \widehat{K}^{\bullet}_g(B)$.
		By (\ref{eq:1.36}) and (\ref{eq:1.37}), the proof of 
		Theorem \ref{thm:1.08} is complete.
		
	\end{proof}

	Now we take $[\underline{E},0]\in \widehat{K}_g^0(W)_{I(g)}$
	for equivariant geometric triple $\underline{E}$ over $W$.
	From (\ref{eq:2.41}),
	\begin{align}
		\widehat{\pi_Z}_!([\underline{E},0])=[(\underline{\pi_Z}, 
		\underline{E}), 0].
	\end{align}
	By (\ref{eq:2.41}), (\ref{eq:2.41a}) and (\ref{eq:2.095}), 
	for $\mN>\mN_0$,
	\begin{multline}
		\widehat{\pi_{Z^{\mathbb{T}}}}_!\circ \left[\underline{
			\lambda_{-1}(N^*)},0\right]^{-1}
		\cup\,\hat{\iota}^*([\underline{E},0])
		=\sum_{\alpha}\left[\left(\underline{\pi_{Z^{\mathbb{T}}_{\alpha}}}, 
		\underline{\lambda_{-1}(N_{\alpha}^*)^{-1}_{\mN}}\widehat{\otimes} 
		\underline{E}|_{W^{\mathbb{T}}_{\alpha}} \right), 0 \right]
		\\
		=\sum_{\alpha}\left[\left(\underline{\pi_{Z^{\mathbb{T}}_{\alpha}}}, 
		\left(\underline{E_{\alpha,1}^{\mN}}\oplus\underline{\left(E_{
		\alpha,2}^{\mN} \right)^{\mathrm{op}}}\right)\widehat{\otimes} 
		\underline{E}|_{W^{\mathbb{T}}_{\alpha}}
		\right), 0\right]/\chi_{\alpha}^{\mN}.
	\end{multline}
	Denote by
	\begin{align}
		\mF:=(\underline{\pi_Z}, \underline{E}),\quad 
		\mF_{\alpha,\mN}:=\left(\underline{\pi_{Z^{\mathbb{T}}_{\alpha}}}, 
		\left(\underline{E_{\alpha,1}^{\mN}}\oplus\underline{\left(
		E_{\alpha,2}^{\mN} \right)^{\mathrm{op}}}\right)\widehat{\otimes} 
		\underline{E}|_{W^{\mathbb{T}}_{\alpha}}
		\right).
	\end{align}
	From Theorem \ref{thm:1.08}, we obtain that 
	\begin{align}\label{eq:3.34}
		[\mF, 0]=
		\sum_{\alpha}\left[\mF_{\alpha,\mN}, 0\right]/\chi_{\alpha}^{\mN}
		\in \widehat{K}^{\bullet}_g(B)_{I(g)}
	\end{align}
	for any $\mN>\mN_0$. 
	Using the index map $I$ in (\ref{d054}) on each side of
	(\ref{eq:3.34}), we obtain the following 
	localization formula for the equivariant family index.
	
	\begin{cor}\label{cor:3.5}
For $\mN>\mN_0$ and $g\in \mathbb{T}^m\backslash A_{\mathbb{T}}$, we have
\begin{align}
\ind D(\mF)=\sum_{\alpha\in \mathfrak{B}}
\ind (D(\mF_{\alpha,\mN}))/\chi_{\alpha}^{\mN}\in K_{\mathbb{T}^m}^{\bullet}(B)_{I(g)}.
\end{align}
	\end{cor}

	Furthermore, set 
	\begin{align}\label{eq:3.35}
		\chi_0^{\mN}:=\prod_{\alpha\in \mathfrak{B}}\chi_{\alpha}^{\mN}.
	\end{align}
	Then by (\ref{eq:2.37a}), (\ref{eq:2.39}) and (\ref{eq:3.34}), 
	there exists $\chi'\in R(G)$ such that 
	\begin{align}\label{eq:3.36}
	\left[\underline{M_{\chi_0^{\mN}\chi'}}\times_B\mF^E, 0\right]=\left[ 
	\sum_{\alpha} \underline{M_{\chi_0^{\mN}\chi'/\chi_{\alpha}^{\mN}}}
	\times_B \mF_{\alpha,\mN}, 0\right]\in \widehat{K}^{\bullet}_g(B).
	\end{align}
	Denote by
	\begin{align}\label{eq:3.36a}
		\mF_M:=\underline{M_{\chi_0^{\mN}\chi'}}\times_B\mF^E, 
		\quad \mF_{M,\mN}^{\mathbb{T}}:=\sum_{\alpha} 
		\underline{M_{\chi_0^{\mN}
		\chi'/\chi_{\alpha}^{\mN}}}\times_B \mF_{\alpha,\mN}.
	\end{align}
	Then by Corollary \ref{cor:3.5}, we have
	$\ind(D(\mF_M))=\ind(D(\mF_{M,\mN}^{\mathbb{T}}))$. Thus by Proposition
	\ref{prop:1.4}, there exists a perturbation operator $\mA$
	associated with $\mF_M+(\mF_{M,\mN}^{\mathbb{T}})^{\mathrm{op}}$.
	
	By Lemma \ref{lemma:2.09}, there exists family cycle
	$(\widetilde{\mF},0)$ such that	$(\mF_M+\widetilde{\mF},0)$ 
	is paired with $(\mF_{M,\mN}^{\mathbb{T}}+\widetilde{\mF},0)$. 
Thus by Definition \ref{defn:2.03},
	there exists perturbation operator
	$\mA'$ such that modulo exact forms on $B$,
	\begin{align}\label{eq:1.40}
		\tilde{\eta}_g(\mF_M+\widetilde{\mF}+\widetilde{\mF}^{\mathrm{op}}
		+(\mF_{M,\mN}^{\mathbb{T}})^{\mathrm{op}},\mA')=0.
	\end{align}
	From (\ref{eq:2.23}), there exists
	perturbation operator $\mA''$ associated 
	$\widetilde{\mF}+\widetilde{\mF}^{\mathrm{op}}$ such that
	$\tilde{\eta}_g(\widetilde{\mF}+\widetilde{\mF}^{\mathrm{op}},\mA'')=0$. 
	So by Proposition \ref{prop:2.01} and (\ref{eq:1.40}), 
	there exists $x\in K_G^{\bullet}(B)$, such that
	\begin{multline}\label{eq:2.45a}
		\ch_g(x)=\tilde{\eta}_g\left((\mF_M
		+(\mF_{M,\mN}^{\mathbb{T}})^{\mathrm{op}})
		+(\widetilde{\mF}+\widetilde{\mF}^{\mathrm{op}}),
		\mA\oplus \mA'' \right)
		\\
		-\tilde{\eta}_g(\mF_M+\widetilde{\mF}+\widetilde{\mF}^{\mathrm{op}}
		+(\mF_{M,\mN}^{\mathbb{T}})^{\mathrm{op}},\mA')
		\\
		=\tilde{\eta}_g\left(\mF_M+(\mF_{M,\mN}^{\mathbb{T}})^{\mathrm{op}},
		\mA \right)+\tilde{\eta}_g(\widetilde{\mF}
		+\widetilde{\mF}^{\mathrm{op}},\mA'')
		-\tilde{\eta}_g(\mF_M+\widetilde{\mF}+\widetilde{\mF}^{\mathrm{op}}
		+(\mF_{M,\mN}^{\mathbb{T}})^{\mathrm{op}},\mA')
		\\
		=\tilde{\eta}_g\left(\mF_M+(\mF_{M,\mN}^{\mathbb{T}})^{\mathrm{op}},
		\mA \right).
	\end{multline}
	
	Since $\mA'$ depends on $g\in \mathbb{T}^m\backslash 
	A_{\mathbb{T}}$,
	a priori, $x$ depends on $g\in \mathbb{T}^m\backslash 
	A_{\mathbb{T}}$. However, as in \cite[Th. 0.5]{LM20}, 
	we will prove in the following that
	it is independent of the choice of $g\in \mathbb{T}^m\backslash 
	A_{\mathbb{T}}$.   It is Theorem \ref{thm:0.03a} in the introduction
	and the main result of this paper.
	In the introduction, we simply denote by
	\begin{align}\label{eq:3.43}
		\chi_{\alpha}:=\chi_{\alpha}^{\mN},\quad
		\chi_0:=\chi_0^{\mN},\quad
		\underline{\mE_{\alpha}}=\underline{E_{\alpha, 1}^{\mN}}\oplus
		\underline{\left(E_{\alpha, 2}^{\mN}\right)^{\mathrm{op}}},\quad
	\mF_{\alpha}:=\mF_{\alpha,\mN},\quad \mF_{M}^{\mathbb{T}}:=
	\mF_{M,\mN}^{\mathbb{T}}
	\end{align}
	for $\mN>\mN_0$.
	
	\begin{thm}\label{thm:1.9}
	For $\mN>\mN_0$ and any perturbation operator $\mathcal{A}$ 
	associated with 
		$\mF_M+(\mF_{M,\mN}^{\mathbb{T}})^{\mathrm{op}}$, there exists $x\in 
		K_{\mathbb{T}^m}^{\bullet}(B)$ such that for any $g\in \mathbb{T}^m
		\backslash A_{\mathbb{T}}$, modulo exact forms on $B$, we have
	\begin{align}
		\tilde{\eta}_g(\mF_M+(\mF_{M,\mN}^{\mathbb{T}})^{\mathrm{op}},
		\mathcal{A})=\ch_g(x).
	\end{align}
	\end{thm}
	\begin{proof}
	Since the group action on $B$ is  trivial, from (\ref{eq:1.005}) and
	(\ref{eq:2.45a}),
	for $g=(g_1,\cdots, g_m)\in  \mathbb{T}^m\backslash A_{\mathbb{T}}$, 
	we can write
		\begin{align}\label{eq:3.40}
		\tilde{\eta}_g(\mF_M+(\mF_{M,\mN}^{\mathbb{T}})^{\mathrm{op}},
		\mathcal{A})
		=\sum_{\max_{k}|l_k|\leq N(g)} a_{l_1\cdots l_m}(g) g_1^{l_{1}}
			\cdots g_m^{l_{m}},
		\end{align}
		where $a_{l_1\cdots l_m}(g)\in \ch(K^{\bullet}(B))$ 
		and  $N(g)\in \mathbb{N}$.
		We only need to prove that $a_{l_1\cdots l_m}(g)$ and 
		 $N(g)$ are independent 
		of $g\in \mathbb{T}^m\backslash A_{\mathbb{T}}$.
		
		For  $K\in \mathfrak{t}^m$, we denote by
		\begin{align}
		f_g(zK):=\tilde{\eta}_{ge^{zK}}(\mF_M
		+(\mF_{M,\mN}^{\mathbb{T}})^{\mathrm{op}},\mathcal{A})
		\end{align}
	for $z\in \C$ and $|z|$ small enough.
	From Proposition \ref{prop:1.09}, since $g\in \mathbb{T}^m
	\backslash A_{\mathbb{T}}$,
	$f_g(zK)=\tilde{\eta}_{g, zK}(\mF_M
	+(\mF_{M,\mN}^{\mathbb{T}})^{\mathrm{op}},\mathcal{A})$.
	By Proposition \ref{prop:1.08}, there exists $\beta_g>0$
	depending on $g$ such that 
	$\frac{\partial f_g(zK)}{\partial \bar{z}}=0$ for $|zK|<\beta_g$. 
	Thus for any compact submanifold $X$ of $B$,
	$\frac{\partial}{\partial \bar{z}}\int_X f_g(zK) =0$
	for $|zK|<\beta_g$.
	So $\int_X f_g(tK)$ is real analytic in $t\in \R$ for $|tK|<\beta_g$.

	Set
	$$\mathscr{G}:=\{g=(g_1,\cdots,g_m)\in \mathbb{T}^m :
	(g_1^{\gamma_1}\cdots g_m^{\gamma_m})_{(
	\gamma_1,\cdots, \gamma_m)\in \Z^m}\ 
	\text{is linearly independent over}\ \Q\}.$$
	Easy to see that $g\in \mathscr{G}$ if and only if
	$g_1,\cdots, g_m\in S^1$ are algebraically independent over $\Q$.
	We claim that \textbf{$\mathscr{G}$ is dense in $\mathbb{T}^m$}.		
	In fact, for any $h=(h_1, \cdots,h_m)\in \mathbb{T}^m$
	and any $\varepsilon>0$, we can construct $g\in \mathscr{G}$
	such that $|h_i-g_i|<\varepsilon$ for any $1\leq i\leq m$
	as follows. For any countable field $\mathscr{K}\subset \C$, denote by
	$A(\mathscr{K})$ the set of complex numbers which are algebraic 
	over $\mathscr{K}$.
	 It means that $x_0\in A(\mathscr{K})$ if and only if there exists
	nonzero $f\in \mathscr{K}[x]$ such that $f(x_0)=0$.
	If $\mathscr{K}$ is countable, a classical result for algebraic number
	theory (see e.g., \cite[Cor. B-2.41]{Ro15}) 
	says that $A(\mathscr{K})$ is countable. In particular,
	$A(\Q)$ is countable. Thus we can take $g_1\in S^1\backslash A(\Q)$ 
	such that $|h_1-g_1|<\varepsilon$. For fixed $g_1$, 
	set $\mathscr{K}_1:=\Q(g_1)$. 
	Then $A(\mathscr{K}_1)$ is countable since $\mathscr{K}_1$ is countable. 
	We take $g_2\in S^1\backslash A(\mathscr{K}_1)$ such that 
	$|h_2-g_2|<\varepsilon$. In general, we can 
	choose $g_i\in S^1\backslash A(\mathscr{K}_{i-1})$, 
	$\mathscr{K}_{i-1}:=\Q(g_1,\cdots,g_{i-1})$, by steps such that
	$|h_i-g_i|<\varepsilon$ for any $1\leq i\leq m$.  
	From the construction we see that such $g=(g_1,\cdots,g_m)$ 
	lies in $\mathscr{G}$.
		
	Fix $\tilde{g}\in \mathscr{G}\cap(\mathbb{T}^m
	\backslash A_{\mathbb{T}})$. There exists a small neighborhood
	 $U_{\tilde{g}}$ of $0\in \mathfrak{t}^m$
	depending on $\tilde{g}$ such that
	$\exp_{\tilde{g}}(U_{\tilde{g}})\subset \mathbb{T}^m
	\backslash A_{\mathbb{T}}$ and for any compact 
	submanifold $X$ of $B$, $\int_Xf_{\tilde{g}}(tK)$
	is real analytic in $t\in \R$ for $tK\in U_{\tilde{g}}$.	
	For any fixed $K\in U_{\tilde{g}}$, we denote by
	$V_{\tilde{g},K}$ the connected component of 
	$\{tK: t\in \R \}\cap U_{\tilde{g}}$ containing $0$.
	Thus $\int_Xf_{\tilde{g}}(tK)$ in analytic in $t$
	for $tK\in V_{\tilde{g},K}$.
		
	Let $T_N=\{g\in \mathbb{T}^m\backslash A_{\mathbb{T}} : N(g)\leq N \}$
	 with $N(g)$ in (\ref{eq:3.40}). Then 
	 $\cup_{N=1}^{\infty}(T_N\cap \exp_{\tilde{g}}(V_{\tilde{g},K})) 
	 =\exp_{\tilde{g}}(V_{\tilde{g},K})$.
	Since $V_{\tilde{g},K}$ is an uncountable set, there exists
	$N_{\tilde{g},K}\in \N$ such that
	 $T_{N_{\tilde{g},K}}\cap \exp_{\tilde{g}}(V_{\tilde{g},K})$
	is an uncountable set.
	Note that $\langle\ch(K^*(B)), [X]\rangle\subset \Q$
	is a countable set, where $[X]$ is the homology
	class represented by $X$ and $\langle \,,\,\rangle$
	is the pair between the cohomology classes and the homology classes.
   Since $a_{l_1\cdots l_m}(g)\in \ch(K^*(B))$, 		
	following the same argument around \cite[(3.61)]{LM20}, 
there exists $(a_{l_1\cdots l_m}^X(K))_{\max_{k}|l_k|\leq N_{\tilde{g},K}}
\in \Q^{m(2N_{\tilde{g},K}+1)}$ 
such that there are uncountable $t\in \R$ such that
		$tK\in V_{\tilde{g},K}$ and 
		\begin{align}\label{eq:3.41a}
\int_Xf_{\tilde{g}}(tK)=\sum_{\max_{k}|l_k|\leq N_{\tilde{g},K}} 
a_{l_1\cdots l_m}^X(K) \tilde{g}_1^{l_{1}}
		\cdots \tilde{g}_m^{l_{m}} e^{t(l_1K_1+\cdots l_mK_m)},
		\end{align}
		where $\tilde{g}=(\tilde{g}_1,\cdots, \tilde{g}_m)$
		and $K=(K_1,\cdots, K_m)\in \mathfrak{t}^m$.
	We claim that $N_{\tilde{g},K}$ and $a_{l_1\cdots l_m}^X(K)$ are also 
	independent of $K$.	In fact, for any $K'\in \mathfrak{t}^m$,
	by (\ref{eq:3.41a}) at $t=0$, we have
		\begin{align}
		\sum_{\max_{k}|l_k|\leq N_{\tilde{g},K}} a_{l_1\cdots l_m}^X(K) 
		\tilde{g}_1^{l_{1}}
		\cdots \tilde{g}_m^{l_{m}}=\int_Xf_{\tilde{g}}(0)
		=\sum_{\max_{k}|l_k|\leq N_{\tilde{g},K'}} a_{l_1\cdots l_m}^X(K') 
		\tilde{g}_1^{l_{1}}\cdots \tilde{g}_m^{l_{m}}.
		\end{align}
		Since $\tilde{g}\in \mathscr{G}$, we have
		$N_{\tilde{g}, K}=N_{\tilde{g}, K'}$ and $a_{l_1\cdots l_m}^X(K)
		=a_{l_1\cdots l_m}^X(K')$.
		Therefore, for any $\tilde{g}\in \mathscr{G}\cap(\mathbb{T}^m
		\backslash A_{\mathbb{T}})$, there exists a neighborhood
		$U_{\tilde{g}}$ of $\tilde{g}$ in $\mathbb{T}^m
		\backslash A_{\mathbb{T}}$ such that (\ref{eq:3.40})
		holds and $N(g)$ and $\langle a_{l_1\cdots l_m}(g),
		[X]\ra$ are independent of $g\in U_{\tilde{g}}$.

		From \cite[Cor. II.30]{T54}, since 
		$B$ is oriented, the rational homology
		group of $B$ have a basis consisting of elements realisable
		by compact submanifolds $X_1,\cdots, X_{s}$.
		Then $\int_{X_i} a_{l_1\cdots l_m}(g)$ is constant
		on $g\in U_{\tilde{g}}$ for any $1\leq i\leq s$.
		Since the pair between the rational cohomology group and the 
		rational homology group is non-degenerate,
		$a_{l_1\cdots l_m}(g)$ is constant
		on $g\in U_{\tilde{g}}$.
		Note that $N(g)$ is constant on $g\in U_{\tilde{g}}$.
		Since $\tilde{g}\in \mathscr{G}\cap(\mathbb{T}^m
		\backslash A_{\mathbb{T}})$
		is arbitrary chosen and $\mathscr{G}$ is dense in $\mathbb{T}^m$,
		$N(g)$ and $a_{l_1\cdots l_m}(g)$
		are constant on each connected component of 
		$\mathbb{T}^m\backslash A_{\mathbb{T}}$.
		
	For $g\in A_{\mathbb{T}}$ and for $t\neq 0$ small, it follows from 
	Proposition \ref{prop:1.09} that
	\begin{align}\label{eq:3.41}
	\tilde{\eta}_{ge^{tK}}(\mF_M+(\mF_{M,\mN}^{\mathbb{T}})^{\mathrm{op}},
	\mathcal{A})
	=\tilde{\eta}_{g,tK}(\mF_M+(\mF_{M,\mN}^{\mathbb{T}})^{\mathrm{op}},
	\mathcal{A})     -\mM_{g,tK},
	\end{align}
	is a meromorphic function for $0\neq t\in \R$, $|t|$ small.
	By Proposition \ref{prop:1.05}, we can choose $K\in \mathfrak{t}^m$
	such that $ge^{tK}$ intersects  transversally with $A_{\mathbb{T}}$
	at $g\in A_{\mathbb{T}}$. From the argument before (\ref{eq:3.41}),
	we know that for $t>0$ small, if we write
	$K=(K_1,\cdots, K_m)\in \mathfrak{t}^m$, $ge^{tK}=(g_{1}e^{tK_1},
	\cdots, g_{m}e^{tK_m})\in \mathbb{T}^m\backslash A_{\mathbb{T}}$, 
	we have
	\begin{align}\label{eq:3.42}
	\tilde{\eta}_{ge^{tK}}(\mF_M+(\mF_{M,\mN}^{\mathbb{T}})^{\mathrm{op}},
	\mathcal{A})
	=\sum_{} a_{l_1\cdots l_m}g_{1}^{l_{1}}\cdots g_{m}^{l_{m}}
		e^{t(l_1K_1+\cdots+l_mK_m)},
	\end{align}
	where $a_{l_1\cdots l_m}\in \ch(K^*(B))$ are independent on $t$ 
	for $t>0$ small. Similarly, for $t<0$ small enough, there exist 
	$a_{l_1\cdots l_m}'\in \ch(K^*(B))$, which do not depend on $t$, 
	such that
	\begin{align}\label{eq:3.43a}
	\tilde{\eta}_{ge^{tK}}(\mF_M+(\mF_{M,\mN}^{\mathbb{T}})^{\mathrm{op}},
	\mathcal{A})
		=\sum_{} a_{l_1\cdots l_m}'g_{1}^{l_{1}}\cdots g_{m}^{l_{m}}
			e^{t(l_1K_1+\cdots+l_mK_m)}.
	\end{align}
	From (\ref{eq:3.41})-(\ref{eq:3.43a}), when $t\to 0$,
	$\lim_{t \to 0}\tilde{\eta}_{ge^{tK}}(\mF_M
	+(\mF_{M,\mN}^{\mathbb{T}})^{\mathrm{op}},\mathcal{A})$ exists. 
	Combining with Proposition \ref{prop:1.09} and (\ref{eq:3.41}), 
	$\tilde{\eta}_{ge^{tK}}(\mF_M+(\mF_{M,\mN}^{\mathbb{T}})^{\mathrm{op}},
	\mathcal{A})$ is analytic in $t$ for $|t|$ small. 
	Set 
	\begin{align}
			f(t)=\sum_{} a_{l_1\cdots l_m}g_{1}^{l_{1}}\cdots 
			g_{m}^{l_{m}}e^{t(l_1K_1+\cdots+l_mK_m)}.
	\end{align}
	Then $f(t)$ is analytic in $t$ for $|t|$ small. By (\ref{eq:3.42}), 
	$f(t)=\tilde{\eta}_{ge^{tK}}(\mF_M
	+(\mF_{M,\mN}^{\mathbb{T}})^{\mathrm{op}},\mathcal{A})$ 
	for $t>0$ small. From the uniqueness of the
	analytic function,  $f(t)=\tilde{\eta}_{ge^{tK}}(\mF_M+
	(\mF_{M,\mN}^{\mathbb{T}})^{\mathrm{op}},\mathcal{A})$ also holds
	 for $t<0$ small. So $a_{l_1\cdots l_m}=a_{l_1\cdots l_m}'$
	for any $(l_1,\cdots, l_m)\in \mathbb{Z}^m$.
	Thus $a_{l_1\cdots l_m}(g)$
	and $N(g)$ in (\ref{eq:3.40}) are constant on different connected 
	components of $\mathbb{T}^m\backslash A_{\mathbb{T}}$.
		
	The proof of this theorem is complete.
	\end{proof}

	\subsection{Proof of Corollaries \ref{thm:0.03} 
	and \ref{cor:0.5}}\label{s0303}

	Take $\mN>\mN_0$.
	If $\ker D(\mF)$ and $\ker D(\mF_{\alpha,\mN})$ form vector bundles 
	over $B$ for any $\alpha$, from Definition \ref{defn:1.01} and 
	(\ref{eq:3.36a}),
	\begin{align}\label{eq:3.45}
	\tilde{\eta}_g(\mF_M)=\chi_0^{\mN}(g)\chi'(g)\tilde{\eta}_g(\mF),\quad
	 \tilde{\eta}_g(\mF_{M,\mN}^{\mathbb{T}})=\chi_0^{\mN}(g)\chi'(g)
	\sum_{\alpha} \chi_{\alpha}^{\mN}(g)^{-1}
	\tilde{\eta}_g(\mF_{\alpha,\mN}).
	\end{align}
	By Proposition \ref{prop:1.5}, for any perturbation
	operator $\mA$ associated with $\mF_M+
	(\mF_{M,\mN}^{\mathbb{T}})^{\mathrm{op}}$, there exists 
	$y\in K_G^{\bullet}(B)$ such that
	for any $g\in \mathbb{T}^m\backslash A_{\mathbb{T}}$,
	 modulo exact forms on $B$,
	\begin{align}\label{eq:3.55}
	\tilde{\eta}_g(\mF_M)-\tilde{\eta}_g(\mF_{M,\mN}^{\mathbb{T}})
	=\tilde{\eta}_g(\mF_M+
	(\mF_{M,\mN}^{\mathbb{T}})^{\mathrm{op}})
	=\tilde{\eta}_g(\mF_M+
	(\mF_{M,\mN}^{\mathbb{T}})^{\mathrm{op}},\mA)+\frac{1}{2}\ch_g(y).
	\end{align}
	Then (\ref{eq:0.07}) follows directly from 
	Theorem \ref{thm:1.9}, (\ref{eq:3.45}) and (\ref{eq:3.55}).
	
	If there exist perturbation operators $\mA$ and $\mA_{\alpha}$
	associated with $\mF$ and $\mF_{\alpha,\mN}$ respectively, 
	from Definition \ref{defn:1.02} and (\ref{eq:3.36a}),
	\begin{align}\label{eq:3.46}
		\begin{split}
		\tilde{\eta}_g(\mF_M, 1\widehat{\otimes}\mA)&=\chi_0^{\mN}(g)
		\chi'(g)\tilde{\eta}_g(\mF,\mA),
		\\
	\tilde{\eta}_g(\mF_{M,\mN}^{\mathbb{T}}, 1\widehat{\otimes}
	(\oplus_{\alpha} \mA_{\alpha}))
	&=\chi_0^{\mN}(g)\chi'(g)\sum_{\alpha}\chi_{\alpha}^{\mN}(g)^{-1}
			\tilde{\eta}_g(\mF_{\alpha,\mN},\mA_{\alpha}).
		\end{split}
	\end{align}
	From Proposition \ref{prop:2.01}, for any perturbation
	operator $\mA'$ associated with $\mF_M+
	(\mF_{M,\mN}^{\mathbb{T}})^{\mathrm{op}}$, there exists
	 $y\in K_G^{\bullet}(B)$ such that
	for any $g\in \mathbb{T}^m\backslash A_{\mathbb{T}}$,
	 modulo exact forms on $B$,
	\begin{multline}\label{eq:3.57}
	\tilde{\eta}_g(\mF_M, 1\widehat{\otimes}\mA)
	-\tilde{\eta}_g(\mF_{M,\mN}^{\mathbb{T}}, 1\widehat{\otimes}
	(\oplus_{\alpha} \mA_{\alpha}))
	=\tilde{\eta}_g(\mF_M+
	(\mF_{M,\mN}^{\mathbb{T}})^{\mathrm{op}}, 1\widehat{\otimes}\mA 
	\oplus 1\widehat{\otimes}(\oplus_{\alpha}
	\mA_{\alpha}^{\mathrm{op}}) )
	\\
	=\tilde{\eta}_g(\mF_M+
	(\mF_{M,\mN}^{\mathbb{T}})^{\mathrm{op}},\mA')+\ch_g(y).
	\end{multline}
	Since the equivariant geometric family does note change,
	the integration part in Proposition \ref{prop:2.01}
	does not appear here.
	Then (\ref{eq:0.08}) follows from
	Theorem \ref{thm:1.9}, (\ref{eq:3.46}) and (\ref{eq:3.57}).
	
	The proof of Corollary \ref{thm:0.03} is complete.

	If the dimension of all fibers of $\mF$ is odd and $\ker D(\mF)$ 
	forms a vector bundle over $B$, by Proposition \ref{prop:1.5}, 
	for any perturbation operator $\mA_1$ associated with $\mF$,
	there exists $x_1\in K_G^{\bullet}(B)$ such that
	for any $g\in \mathbb{T}^m\backslash A_{\mathbb{T}}$,
	modulo exact forms on $B$,
	\begin{align}\label{eq:3.58}
	\tilde{\eta}_g(\mF_M)
	=\tilde{\eta}_g(\mF_M,\mA_1)+\frac{1}{2}\ch_g(x_1).
	\end{align}
	As in (\ref{eq:3.57}), by Proposition \ref{prop:2.01}, for any 
	perturbation operators $\mA_2$ and $\mA$ associated with $\mF_M+
	(\mF_{M,\mN}^{\mathbb{T}})^{\mathrm{op}}$ and 
	$\mF_{M,\mN}^{\mathbb{T}}$ respectively,
	there exists $x_2\in K_G^{\bullet}(B)$ such that
	for any $g\in \mathbb{T}^m\backslash A_{\mathbb{T}}$,
	modulo exact forms on $B$,
	\begin{align}\label{eq:3.59}
	\tilde{\eta}_g(\mF_M,\mA_1)-\tilde{\eta}_g(\mF_{M,\mN}^{\mathbb{T}},\mA)
	=\tilde{\eta}_g(\mF_M+
	(\mF_{M,\mN}^{\mathbb{T}})^{\mathrm{op}},\mA_2)+\ch_g(x_2).
	\end{align}
	Then (\ref{eq:0.07b}) follows from Theorem \ref{thm:1.9},
	the first formula
	of (\ref{eq:3.45}), (\ref{eq:3.58}) and (\ref{eq:3.59}).

	The proof of Corollary \ref{cor:0.5} is complete.
	
	
	\subsection{Proof of Theorem \ref{thm:0.7}}\label{s0304}
	
In this subsection, we assume that $W^{\mathbb{T}}=\emptyset$.

Let $\pi_{Z!}:K_{\mathbb{T}^m}^0(W)\to K_{\mathbb{T}^m}^{\bullet}(B)$, 
$[E]\mapsto \ind(D(\mF))$ be the pushforward map associated with 
the equivariant topological K-groups (cf. e.g., \cite[Th. 2.9]{Liu21}). 
Then it induced a well-defined localized map
$\pi_{Z!}:K_{\mathbb{T}^m}^0(W)_{I(g)}\to 
K_{\mathbb{T}^m}^{\bullet}(B)_{I(g)}$.
From \cite[Prop. 1.5]{ASegal68}, if 
$g\in \mathbb{T}^m\backslash A_{\mathbb{T}}$, 
$K_{\mathbb{T}^m}^0(W)_{I(g)}=0$.
Thus we have $\ind(D(\mF))=0\in K_{\mathbb{T}^m}^{\bullet}(B)_{I(g)}$. 
So there exists $\chi\in R(\mathbb{T}^m)$ such that 
$\ind(D(\underline{M_{\chi}}\otimes\mF))=0\in K_{\mathbb{T}^m}^{\bullet}(B)$.

From (\ref{d034}), since $W^{\mathbb{T}}=\emptyset$,
if $g\in \mathbb{T}^m\backslash A_{\mathbb{T}}$,
$I:\widehat{K}_g^{\bullet}(W)\to K_{\mathbb{T}^m}^{\bullet}(W)$
is an isomorphism. By (\ref{d054}), we have 
$[\pi_Z^*(\underline{M_{\chi}})\otimes E , 0 ]=0\in \widehat{K}_g^{0}(W)$.
Note that by \cite[Th. 4.7]{Liu21}, $\widehat{\pi_Z}_!:
\widehat{K}_g^{0}(W)\to \widehat{K}_g^{\bullet}(B)$
in (\ref{d082}) is
also well-defined when $W^{g}=\emptyset$.
Thus 
\begin{align}\label{eq:3.60}
\widehat{\pi_Z}_![\pi_Z^*(\underline{M_{\chi}})\otimes E , 0 ]
=[\underline{M_{\chi}}\otimes\mF, 0]=0\in 
\widehat{K}_g^{\bullet}(B).
\end{align}
From Definition \ref{defn:2.03}, for any perturbation operator $\mA$ 
associated with $\underline{M_{\chi}}\otimes\mF$,
$(\underline{M_{\chi}}\otimes\mF, 0)$ is paired with
$(0, \widetilde{\eta}_g(\underline{M_{\chi}}\times_B\mF, \mA))$. 
So from (\ref{eq:3.60}), 
$[0, \widetilde{\eta}_g(\underline{M_{\chi}}\times_B\mF, \mA)]
=a(-\widetilde{\eta}_g(\underline{M_{\chi}}\times_B\mF, \mA))=0
\in \widehat{K}_g^{\bullet}(B)$.
By (\ref{d034}), there exists $x\in K_G^{\bullet}(B)$
such that modulo exact forms on $B$,
\begin{align}
\widetilde{\eta}_g(\underline{M_{\chi}}\times_B\mF, \mA)=\ch_g(x).
\end{align}	
Since Propositions \ref{prop:2.01} and \ref{prop:1.09} 
still hold if $W^{\mathbb{T}}=\emptyset$, following the same process 
as in Theorem \ref{thm:1.9}, we obtain that 
$x$ is independent of $g$, which implies (\ref{eq:0.08a}).
	
If $\ker D(\mF)$ forms a vector bundle over $B$,
by Proposition \ref{prop:1.5} and (\ref{eq:0.08a}), there exists 
$x_1\in K_G^{\bullet}(B)$ such that
for any $g\in \mathbb{T}^m\backslash A_{\mathbb{T}}$,
 modulo exact forms on $B$,
\begin{align}
\widetilde{\eta}_g(\underline{M_{\chi}}\times_B\mF)
=\widetilde{\eta}_g(\underline{M_{\chi}}\times_B\mF, \mA)
+\frac{1}{2}\ch_g(x_1)=\ch_g(x)+\frac{1}{2}\ch_g(x_1).
\end{align}	
Since $\widetilde{\eta}_g(\underline{M_{\chi}}\times_B\mF)
=\chi(g)\widetilde{\eta}_g(\mF)$, we obtain (\ref{eq:0.08b}).

The proof of Theorem \ref{thm:0.7} is complete.

	\appendix

	\section{Characterization of Chern-Simons classes}\label{appendix}
	
	In this appendix, we introduce a unique characterization
	for equivariant Chern-Simons classes used in this paper,
	which is a parallel version of \cite[Th. 1.29]{BGSI}
	and \cite[Th. 2.10]{Bi98}.
	We assume that $G$ is any compact Lie group.
	
	\begin{lemma}\label{lemma:2.06}
		For $g\in G$, there exists a unique way to attach to
		two equivariant ungraded geometric triples
		$\underline{E_1}=(E, h_1^E, \nabla_1^E)$ and 
		$\underline{E_2}=(E, h_2^E, \nabla_2^E)$ on a compact
		$G$-manifold $B$ with the same $G$-equivariant vector
		bundle $E$
		a class $\hat{\eta}_g
		(\underline{E_1}, \underline{E_2})$ in $\Omega^{\mathrm{odd}}
		(B^g, \mathbb{R})/\mathrm{Im}\,d$ such that
	\begin{enumerate}
	\item $d \hat{\eta}_g
		(\underline{E_1}, \underline{E_2})=\mathrm{ch}_g(\underline{E_2})
		-\mathrm{ch}_g(\underline{E_1})$;
		\item for every $G$-equivariant smooth map $f: B'\to B$, 
		$\hat{\eta}_g
		(f^*\underline{E_1}, f^*\underline{E_2}) =f^* \hat{\eta}_g
		(\underline{E_1}, \underline{E_2})$;
	\item if $\underline{E_1}=\underline{E_2}$, 
		$\hat{\eta}_g(\underline{E_1}, \underline{E_2})=0$.
	\end{enumerate}
	In this case, the unique class $\hat{\eta}_g
	(\underline{E_1}, \underline{E_2})$ is just the 
	equivariant Chern-Simons class $\widetilde{\ch}_g\left(\underline{E_1},
	\underline{E_2}\right) $ in $\Omega^{\mathrm{odd}}
	(B^g, \mathbb{R})/\mathrm{Im}\,d$.
	\end{lemma}
	\begin{proof}
		Since $\widetilde{\mathrm{ch}}_g(\underline{E_1}, \underline{E_2})$
		satisfies
		all conditions, the existence is obvious. 
		For the uniqueness,
		if $\hat{\eta}_g$ verifies (1), (2) and (3),
		we consider the equivariant geometric triples
		over $B\times [0,1]$, where the $G$-action on $[0,1]$ is trivial.
		Let $\underline{\widetilde{E}_1}$
	and $\underline{\widetilde{E}_2}$ be two equivariant geometric triples
	over $B\times [0,1]$ such that $\widetilde{E}_1|_{s\in [0,1]}=
	\widetilde{E}_2|_{s\in [0,1]}=E$, $h_1^{\widetilde{E}}|_{s=0}=h_1^E$,
		$h_1^{\widetilde{E}}|_{s=1}=h_2^E$,
		$\nabla_1^{\widetilde{E}}|_{s=0}=\nabla_1^E$,
		$\nabla_1^{\widetilde{E}}|_{s=1}=\nabla_2^E$,
		$h_2^{\widetilde{E}}|_{s\in [0,1]}=h_2^E$,
		$\nabla_2^{\widetilde{E}}|_{s\in [0,1]}=\nabla_2^E$.
	Then $\hat{\eta}_g
	(\underline{\widetilde{E}_1}, \underline{\widetilde{E}_2})\in
	\Omega^{\mathrm{odd}}(B^g\times [0,1], \R)/\mathrm{Im}\,d$ satisfying 
    the conditions (1)-(3). From the condition (1), we have
		\begin{align}\label{eq:2.16}
			d \hat{\eta}_g
			(\underline{\widetilde{E}_1}, \underline{\widetilde{E}_2})
			=\mathrm{ch}_g(\underline{\widetilde{E}_2})
			-\mathrm{ch}_g(\underline{\widetilde{E}_1}).
		\end{align}
		Taking the fiber integration of (\ref{eq:2.16}) along $[0,1]$,
		by \cite[(1.29)]{LM20},	we have
		\begin{align}\label{eq:2.26}
			\int_0^1\frac{\partial}{\partial s}
			\hat{\eta}_g
			(\underline{\widetilde{E}_1}, \underline{\widetilde{E}_2})ds
			=-\int_0^1\{\mathrm{ch}_g(\underline{\widetilde{E}_1}) \}^{ds}ds
			=-\widetilde{\mathrm{ch}}_g(\underline{E_1}, \underline{E_2}).
		\end{align} 
		From the conditions (2) and (3), 
		\begin{multline}\label{eq:2.27}
			\int_0^1\frac{\partial}{\partial s}
			\hat{\eta}_g
			(\underline{\widetilde{E}_1}, \underline{\widetilde{E}_2})ds
			=i_1^*\hat{\eta}_g
			(\underline{\widetilde{E}_1}, \underline{\widetilde{E}_2})
			-i_0^*\hat{\eta}_g
			(\underline{\widetilde{E}_1}, \underline{\widetilde{E}_2})
			\\
			=\hat{\eta}_g
			(\underline{E_2}, \underline{E_2})-\hat{\eta}_g
			(\underline{E_1}, \underline{E_2})=-\hat{\eta}_g
			(\underline{E_1}, \underline{E_2}),
		\end{multline}
		where $i_0:B\times\{0\}\to B\times [0,1]$ and 
		$i_1:B\times\{1\}\to B\times [0,1]$ are natural embeddings.
		Thus from (\ref{eq:2.26}) and (\ref{eq:2.27}), 
		we obtain the uniqueness.
	\end{proof}

	\begin{cor}\label{cor:2.09}
		For equivariant ungraded geometric triples $\underline{E_i}=
		(E, h^E_i, \nabla^E_i)$ over $B$, $i=1,2$, denote by $\mathcal{F}_i
		=\underline{E_i}$ the equivariant geometric family over $B$ in which 
		the fiber is a point. Then there exists a perturbation operator
		$\mathcal{A}$ associated with $\mathcal{F}_1+
		\mathcal{F}_2^{\mathrm{op}}$ such that
		\begin{align}\label{eq:2.28}
			\widetilde{\eta}_g((\mF_1+\mF_2^{\mathrm{op}})|_{B^g},
			\mA|_{B^g})=\widetilde{\ch}_g\left(\underline{E_1},
			\underline{E_2}\right).
		\end{align}
	\end{cor}
	\begin{proof}
		We only need to verify that the equivariant eta form
		$\widetilde{\eta}_g((\mF_1+\mF_2^{\mathrm{op}})|_{B^g},
		\mA|_{B^g})$ satisfies
		the conditions in Lemma \ref{lemma:2.06}.
		By Definition \ref{defn:1.03}, as in (\ref{eq:2.34}),
		\begin{align}\label{eq:a.06}
			\mathcal{A}=	
			\left(\begin{array}{cc}
				0 & \mathrm{Id} \\
				\mathrm{Id} & 0
			\end{array}
			\right)
		\end{align}
		is a perturbation operator associated with $\mathcal{F}_1+
		\mathcal{F}_2^{\mathrm{op}}$. 
		For $G$-equivariant smooth map $f:B'\to B$, $f^*\mathcal{A}$
		is also a perturbation operator associated with $f^*
		(\mathcal{F}_1+
		\mathcal{F}_2^{\mathrm{op}}) $. So from Definition \ref{defn:1.02}
		(cf. \cite[(4.38)]{Liu21}), we have
		\begin{align}\label{eq:2.31}
			\widetilde{\eta}_g(f^*(\mF_1+\mF_2^{\mathrm{op}})|_{B^g},
			f^*\mA|_{B^g})=f^*\widetilde{\eta}_g(
			(\mF_1+\mF_2^{\mathrm{op}})|_{B^g},
			\mA|_{B^g}).
		\end{align}
		If $\underline{E_1}
		=\underline{E_2}$, from (\ref{eq:1.16a}) and 
		(\ref{eq:1.23a}), we have
		\begin{align}\label{eq:2.30}
		(\mathbb{B}_u+\sqrt{u}\chi(\sqrt{u})\mathcal{A})^2
		=\left(\begin{array}{cc}
				(\nabla^E_1)^2+u\chi^2(\sqrt{u})\mathrm{Id} & 0 \\
				0 & (\nabla^E_2)^2+u\chi^2(\sqrt{u})\mathrm{Id}
			\end{array}
			\right).
		\end{align}
		Thus $\mathrm{Tr}_s[g\mathcal{A}\exp(-(\mathbb{B}_u
		+\sqrt{u}\chi(\sqrt{u})	\mathcal{A})^2)]=0$.
		From Definition \ref{defn:1.02}, we have 
		$\widetilde{\eta}_g((\mF_1+\mF_2^{\mathrm{op}})|_{B^g},
		\mA|_{B^g})=0$.
		So from (\ref{eq:1.25a}) and Lemma \ref{lemma:2.06}, we obtain 
		(\ref{eq:2.28}).
	\end{proof}
	

\end{document}